\documentclass{article}
\usepackage{amsmath}
\usepackage{amsfonts}
\usepackage{graphicx}
\newcounter{remark}
\newtheorem{proposition}{Proposition}
\def\theremark {\arabic{remark}}
\newenvironment{remark}{\refstepcounter{remark}\par\noindent{\bf Remark\ \theremark}\ }{\par}
\newenvironment{proof}{\noindent{\it Proof.\/}\rm}{\hfill $\Box$\medskip }
\newcounter{example}
\def\theexample{\arabic{example}}
\newenvironment{example}{\refstepcounter{example}\par\noindent{\bf Example\ \theexample}\ }{\par}

\title{Shishkin mesh simulation: A new stabilization
technique for convection-diffusion problems and complementary experiments}
\author{Bosco
Garc\'{\i}a-Archilla\thanks{Departamento de Matem\'{a}tica Aplicada
II, Universidad de Sevilla, Sevilla, Spain. Research supported by
Spanish MEC under grant MTM2009-07849 (bosco@esi.us.es)}
}
\begin{document}
\maketitle

\begin{abstract}
We present together in this document the paper B. Garc\'{\i}a-Archilla, Shiskin mesh simulation:
A new stabilization
technique for convection-diffusion problems, {\em Comput. Methods Appl. Mech. Engrg.\/}, 256 (2013),
1--16,
and the manuscript  B. Garc\'{\i}a-Archilla, Shiskin mesh simulation: Complementary experiments,
which contains some of the experiments that, for the sake of brevity were
not included~in the first one. {\it This second paper with complementary experiments is neither self-contained nor intended to be published
in any major journal\/}, but is intended to be accessible to anyone wishing to learn more on the
performance of the SMS method. Following a principle of reproducible computational research,
the source codes of the experiments in the present papers are available from the author on request or
on 
"{\tt http://personal.us.es/bga/bga\_files/software\_bga.html}". Excluded are the codes
corresponding to Example 7 which is subject of current resarch
\end{abstract}

\newpage
\noindent{}
\vskip1.2truecm
\setcounter{page}{1}

\begin{center}
{\Large
Shishkin mesh simulation: A new stabilization}
\medskip

{\Large technique for convection-diffusion problems}
\bigskip

{\large Bosco
Garc\'{\i}a-Archilla}\footnote{Departamento de Matem\'{a}tica Aplicada
II, Universidad de Sevilla, 41092 Sevilla, Spain. Research supported by
Spanish MEC under grant MTM2009-07849 (bosco@etsi.us.es)}
\medskip

{\large
}
\vskip 1truecm

\begin{abstract}
A new stabilization procedure is presented. 
It is based on a simulation of the interaction between the coarse and fine parts of a
Shishkin mesh, but can be applied on coarse and irregular meshes and
on domains with nontrivial geometries.
The technique, which does not require adjusting any parameter, can be applied to different
stabilized and non stabilized methods.
Numerical experiments show it to obtain oscillation-free approximations
on problems with boundary and internal layers,
on uniform and nonuniform meshes and on domains with curved boundaries.

\end{abstract}

 Finite-element methods. Galerkin method. SUPG method. 

\end{center}
\section{Introduction}

The numerical solution of convection-diffusion problems when convection dominates is, despite
more than 30 years of research, a challenging problem nowadays. Standard finite-element  or
finite-difference methods typically suffer from unphysical or spurious oscillations unless meshes are
taken so fine that are useless for all practical purposes. The reason is the presence of
layers or thin regions where solutions change fast. Modification of standard methods,
known as stabilized methods have been proposed in the literature, from upwind
methods 35 years ago~\cite{Tabata}, to strongly-consistent stabilized methods like the
streamline upwind/Petrov-Galerkin (SUPG) method~\cite{SUPG}, also known as the streamline diffusion
finite element method (SDFEM), or the Galerkin least-squares (GALS) method~\cite{GALS}.
More recently, local projection stabilization (LPS) methods, \cite{LP1}, \cite{LP2}, \cite{TobiskaLP},
continuous interior
penalty (CIP) methods~\cite{CIP1}, or discontinuous Galerkin (DG) methods \cite{HSS02},
\cite{Ri92} have been introduced, to cite a few of the many
techniques proposed (see~\cite{Roos-Stynes-Tobiska}, \cite{Stynes-Acta} for a survey of methods).
It must be noticed, however, that computational studies (see e.g., \cite{Augustin-et-al-2011}, \cite{John-Knobloch-2007}) find it hard to put a particular method above the others. It must be also mentioned
that most of these methods depend on at least one parameter about which there is no unanimous
agreement on its optimal choice in practical problems~\cite{John-Knobloch-Savescu}.

A different approach is to use layer-adapted meshes. Among these we cite
Shishkin meshes (described below)~\cite{Shishkin2}, \cite{Shishkin}, which
have received considerable attention in recent years \cite{Franz-Kellogg-Stynes},
\cite{Franz-Linss-Roos}, \cite{Franz-Matthies}, \cite{Kopteva-Oriordan}, 
\cite{Linss}, \cite{Linss-Stynes}
\cite{Oriordan-Stynes}, \cite{Stynes-Tobiska} and~\cite{Zhang}. However, it is generally acknowledged
that the main drawback of Shishkin meshes is the difficulty to design them on domains with nontrivial
geometries, although some works overcoming this difficulty can be found in the literature
\cite{Xenophontos-Fulton}, \cite{Kopteva-Oriordan}.

The method we propose, however, does not suffer from the above indicated drawbacks: It does
not depend on parameters and, although it is based on the idea of simulating a Shishkin mesh,
the experiments we present show it produces excellent results on domains with nontrivial geometries.

We consider the problem
\begin{align}
&-\varepsilon \Delta u+b\cdot \nabla u+ cu=f,\quad {\rm in}\quad \Omega,\label{eq:model}\\
&u=g_1,\quad {\rm in}\ \partial \Omega_D,\quad \frac{{\partial
u}}{{\partial n}}=g_2,\quad {\rm in}\ \partial
\Omega_N.\label{eq:modelbc}
\end{align}
Here, $\Omega$ is a bounded
domain in ${\mathbb R}^d$, $d=1,2,3$,
its boundary $\partial\Omega$ being the disjoint union of~$\Gamma_D$ and~$\Gamma_N$,
$b$ and $c$ are given functions and $\varepsilon> 0$ is a
constant diffusion coefficient. We
assume that $\Gamma^{-}\subset
\partial \Omega_D$, $\Gamma^{-}$ being the inflow boundary of
$\Omega\subset \mathbb{R}^d$, i.e., the set of points $x\in \partial \Omega$ such that
$b(x)\cdot n(x)<0$.

It is well-known if $\varepsilon \ll \sup\{ \left|b(x)\right|\mid x\in\Omega\}$ ($\left|\cdot\right|$ being
the euclidean norm) boundary layers are likely to develop along $\partial\Omega\backslash
\Gamma^{-}$, although they have different structure on~$\Gamma^{0}=
\{ x\in\partial \Omega\mid b(x)\cdot n(x)=0\}$ and~$\Gamma^{+}=
\{ x\in\partial \Omega\mid b(x)\cdot n(x)>0\}$. As already mentioned, these boundary layers,
when present, are
responsible of the spurious oscillations that pollute the numerical approximations obtained
with standard methods unless extremely fine meshes are used. For uniform meshes,
oscillations typically disappear when the {\it mesh P\'eclet number}
$$
\hbox{\rm Pe}=\frac{\left\|b\right\|_{L^\infty(\Omega)^2} h}{2\varepsilon}
$$
($h$ being the mesh size) is of the order of~1.

Let us briefly  describe now the idea of the method we propose in the following simple problem:
\begin{align}
\label{simple1d}
 L(u)&\equiv -\varepsilon u''(x)+b(x)u'(x)+c(x)=f(x),\qquad 0<x<1,
\\
\label{simple1dbc}
u(0)&=u(1)=0.
\end{align}
In~(\ref{simple1d}) we assume that $b$, $c$ and~$f$ are sufficiently smooth functions,
and that
\begin{equation}
0<\beta<\min_{x\in[0,1]} b(x),\qquad 0\le \min_{x\in[0,1]} c(x).
\label{cond1d}
\end{equation}
The standard Galerkin linear finite-element method for~(\ref{simple1d}-\ref{simple1dbc})
on a partition or mesh $0=x_0<x_1<\ldots< x_J=1$ of~$[0,1]$ obtains a continuous
piecewise linear approximation
$u_h(x)$ to~$u$. As it is customary, $h$ denotes the mesh diameter, $h=\max_{1\le j\le J}h_j$,
where $h_j=x_j-x_{j-1}$, for $j=1,\ldots,J$.
The approximation
can be expressed as $U(x)=u_1\varphi_1(x)+\cdots+u_{J-1}\varphi_{J-1}(x)$, where the
$\varphi_j(x)$ are the basis or hat (piecewise linear) functions taking value 1 at the node
$x_j$ and $0$ in the rest of the nodes of the partition (thus, $U(x_j)=u_j$). The values $u_j$,
$j=1,\ldots,J-1$, are obtained by solving the linear system of equations
\begin{equation}
\label{simple1dgal}
a(u_h,\varphi_i)=(f,\varphi_i), \qquad i=1,\ldots,J-1,
\end{equation}
where, $a$ is the bilinear form associated with~(\ref{simple1d}), which is given by
$$
a(v,w)=\varepsilon(v',w')+(bv'+c,w),
$$
$(\cdot,\cdot)$ being the standard inner product in~$L^2(0,1)$,
$$
(f,g)=\int_0^1f(x)g(x)\,\hbox{\rm d}x.
$$

The Shishkin mesh with $J=2N$ nodes is composed of two uniform meshes with $N$ subintervals
on each side of the transition point $x_N=1-\sigma$, where
$$
\sigma=\min(\frac{1}{2},\frac{2}{\beta}\varepsilon\log N),
$$
for an adequate constant $\beta$, 
that is, $x_j=j(1-\sigma)/N$, for $j=0,\ldots,N$, and $x_{N+j}=x_N+j\sigma/N$, for
$j=1,\ldots,N$.   Let us consider the coarse and fine grid
parts of the Galerkin approximation given by
\begin{align*}
U_c(x)&=u_1\varphi_1(x)+\cdots+u_{N-1}\varphi_{N-1}(x),\\
U_f(x)&=u_{N+1}\varphi_{N+1}(x)+\cdots+u_{2N-1}\varphi_{2N-1}(x),
\end{align*}
so that $U_c+u_N\varphi_N+U_f$ is the Galerkin approximation on the Shishkin mesh.
Since for $i=1,\ldots,J-1$, the support of the basis function
$\varphi_i$ is~$[x_{i-1},x_{i+1}]$, we have
$a(U_c,\varphi_{N+j})=0$ and~$a(U_f,\varphi_{j})=0$, for $j=1,\ldots,N-1$.
Consequently the system~(\ref{simple1dgal}) on the Shishkin mesh
can be rewritten as
\begin{align}
a(U_c,\varphi_i)\,\hphantom{{}+u_Na(\varphi_N,\varphi_i)}
\hphantom{{}+a(U_f,\varphi_N)}\,\,&=(f,\varphi_i),\quad i=1,\ldots,N-2,
\label{sis1d1a}
\\
a(U_c,\varphi_i)\,{}+u_Na(\varphi_N,\varphi_i)
\hphantom{{}+a(U_f,\varphi_N)}\,\,&=(f,\varphi_i),\quad i=N-1,
\label{sis1d1b}
\\
\vphantom{\bigg|}
a(U_c,\varphi_N)+u_Na(\varphi_N,\varphi_N)
+a(U_f,\varphi_N)&=(f,\varphi_N),
\label{sis1d2a}
\\
\hphantom{a(U_c,\varphi_i)+{}}u_Na(\varphi_N,\varphi_i)\,
+\,a(U_f,\varphi_i)\,\,&=(f,\varphi_i),\quad i=N+1,\ldots, 2N-1.
\label{sis1d2b}
\end{align}

\begin{figure}
\begin{center}
\includegraphics[height=4.5cm]{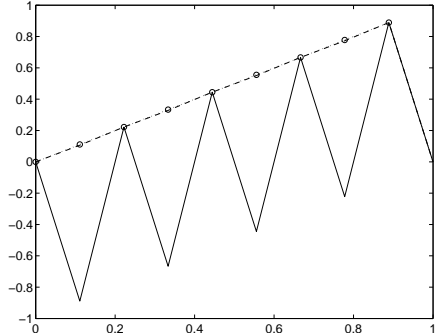}
\end{center}
\vspace{-0.5cm}
\caption{Galerkin approximation on a uniform mesh with $N=9$  (continuous line) to the solution
of~(\ref{simple1d2}--\ref{simple1d2bc}) with $\varepsilon=10^{-8}$, $\sigma=4\varepsilon\log(2*N)$,
$b=f=1$, and $c=0$.
The~$U_c$ part of the
Galerkin approximation on a Shishkin mesh with $J=18$ (broken line) for same $\varepsilon$
and~$f$.  Circles are the
values of the true solution on the nodes.}
\label{Fig_ex1}
\end{figure}

We notice that were it not for the presence of the $u_Na(\varphi_N,\varphi_i)$
in~(\ref{sis1d1b}), the system~(\ref{sis1d1a}--\ref{sis1d1b}) would be the
equations
\begin{equation}
\label{eq:galsigma}
a(U,\varphi_i)=(f,\varphi_i),\qquad i=1,\ldots,N-1,
\end{equation}
of the Galerkin approximation~$U=U_1\varphi_1+\cdots+U_{N-1}\varphi_{N-1}$
for the problem
\begin{align}
\label{simple1d2}
 -\varepsilon u''(x)+b(x)u'(x)+c(x)&=f(x),\qquad 0<x<1-\sigma,
\\
\label{simple1d2bc}
u(0)=u(1-\sigma)&=0.
\end{align}
The Galerkin approximation~$U$ for this problem, unless $\varepsilon N>1/2$,
is likely to have spurious oscillations of large amplitude
as we show in Fig.~\ref{Fig_ex1} for $\varepsilon=10^{-8}$, $\sigma=4\varepsilon\log(J)$, $b(x)=f(x)=1$, $c=0$ and $N=9$.
It is however the presence
of~$u_Na(\varphi_N,\varphi_i)$ in equation~(\ref{sis1d1b}) that suppresses the oscillations, as we
can see in~Fig.~\ref{Fig_ex1}, where the component $U_c$ of the Galerkin approximation
on a Shishkin grid with $J=2N=18$ is also shown (discontinuous line) together with the true
solution at the nodes of the coarse part of the mesh.

It is remarkable that just by adding the value
\begin{equation}
\alpha^*=u_Na(\varphi_{N-1},\varphi_N)
\label{alpha*}
\end{equation}
to the last equation of the Galerkin method for~(\ref{simple1d2}--\ref{simple1d2bc}) we get the oscillation-free approximation~$U_c$. 
%
Obviously, in order to have the value of $\alpha^*$ we have to solve the whole system (\ref{sis1d1a}--\ref{sis1d2b}).
In the present paper, we introduce a technique to approximate~$\alpha^*$ without the need to
compute the whole approximation on the Shishkin grid. In Fig.~\ref{Fig_ex1}, the approximation
computed with the estimated~$\alpha^*$ is indistinguishable from~$U_c$. Numerical experiments in the present paper
show that, in two-dimensional problems, the oscillation-free approximation on a coarse mesh can
be obtained by this technique at half the computational cost of a Shishkin grid, and a more substantial
gain can be expected in three-dimensional problems.

Furthermore, this technique can be extended
when the grid is no part of any Shishkin grid, while, at the same time, managing to get rid of the
spurious oscillations. This allows to obtain accurate
approximations on domains with non trivial boundaries, where Shishkin meshes may be
difficult to construct. In spite of this, we call the new technique {\it Shishkin mesh simulation\/}
(SMS), since it was derived, as described above, in an attempt to simulate Shishkin grids.

We must mention, however, that in the present paper we only consider the case of dominant convection,
both in the analysis and in the numerical experiments. The question of how to modify the method
(if necessary) when
the mesh P\'eclet number~$\hbox{\rm Pe}$ tends to one will be addressed elsewhere.

It is well-known that the Galerkin method is a far from ideal method in convection-diffusion problems.
Let us also notice that despite the good properties of stabilized methods developed
in recent years, the SUPG method is still considered the standard
approach~\cite{John-Knobloch-Savescu}. For this reason, in the numerical experiments
presented below, we compare the new method with the SUPG method. 

The rest of the paper is as follows. In Section~\ref{Se:main} we describe the SMS method
in detail. In Section~\ref{Se:analisis} we present a limited analysis of the new technique.
Section~\ref{Se:experiments} contains the numerical experiments and Section~\ref{Se:conclusions}
the conclusions. 

\section{The new technique: Shishkin mesh simulation}
\label{Se:main}

\subsection{The one-dimensional case}

We consider~(\ref{simple1d}--\ref{simple1dbc}) satisfying~(\ref{cond1d}). Given 
a partition
$0=x_0<x_1<\ldots< x_J=1$, of $[0,1]$,  we denote by $X_{h}$
the space of continuous piecewise
linear functions, and by~$V_{h}$ the subspace of~$X_{h}$ of
functions taking zero values at $x=0$ and $x=1$, so that we can express
$$
X_h=\hbox{\rm span}\{\varphi_0\}\oplus V_h\oplus\hbox{\rm span}\{\varphi_J\} . 
$$
We consider the operator
$L_h$ given by
$$
L_h(v_h)=bv_h'+cv_h,\qquad v_h\in V_h.
$$
(See also Remark~\ref{re:higher_order} below).
We denote by~$u_h\in V_h$ the standard Galerkin linear finite-element approximation, and by~$\tilde u_h$ the
SMS approximation. This is found as the solution of the least-squares problem
\begin{equation}
\min_{\tilde u_h\in V_h ,\alpha\in{\mathbb R}} \left\| L_h( \tilde u_h)-f\right\|_{L^2(0,x_{J-1})},
\label{sms1d1}
\end{equation}
subject to the restriction
\begin{equation}
a(\tilde u_h,\varphi_h)+\alpha \varphi_h(x_{J-1})=(f,\varphi),\qquad\varphi_h\in V_h.
\label{sms1d2}
\end{equation}
Observe that since $\varphi_h(x_{J-1})=0$ for $\varphi_h$ not proportional to~$\varphi_{J-1}$,
and recalling how we defined~$\alpha^*$ in~(\ref{alpha*}), equation
(\ref{sms1d2}) is similar to~(\ref{sis1d1a}--\ref{sis1d1b}), and then, the least-squares
problem~(\ref{sms1d1}) is the way of finding the value $\alpha$ hopefully close to~$\alpha^*$.



Notice also that the restriction~(\ref{sms1d2}) is, in fact, a set of as many
independent restrictions
as interior nodes or nodal basis functions in~$V_h$. Consequently, in the optimality
conditions, there must be a Lagrange multiplier for every interior node. We
gather all this multipliers in a function~$z_h\in V_h$ whose value at every
node is that of the corresponding Lagrange multiplier. Thus, 
the optimality conditions of the SMS approximation can be written as the following linear problem:
find $\tilde u_h,z_h\in V_h$ and~$\alpha\in{\mathbb R}$ such that
\begin{align}
\label{opt1d1}
\bigl(L_h(\tilde u_h),L_h(\varphi_h)\bigr)_{L^2(0,x_{J-1})} - a_h(\varphi_h,z_h) &=
(f,L_h(\varphi_h))_{L^2(0,x_{J-1})},\qquad \varphi_h\in V_h,\\
z_h(x_{J-1})& =0,\label{opt2}\\
a(\tilde u_h,\varphi_h)+\alpha \varphi_h(x_{J-1})&=(f,\varphi),\qquad\varphi_h\in V_h,
\label{opt1d3}
\end{align}
where here and in the sequel~$(\cdot,\cdot)_{L^2(I)}$ denotes the standard inner product
in~$L^2(I)$.


\begin{remark}\label{re:higher_order}
Observe that for linear elements the operator $L_h$ coincides with
\begin{equation}
L_h(v_h)=\sum_{j=1}^{J} (-\varepsilon v_h''+bv_h'+cv_h)_{|_{(x_{j-1},x_j)}},
\label{L_h_higher}
\end{equation}
that is, $L$ applied element by element. This expression of~$L_h$ is better suited to the
SMS method for higher-order elements, a topic that will be studied elsewhere.
\end{remark}

\medskip
\begin{remark}\label{re:GALS} The fact that the approximation $\tilde u_h$ is found by solving
the least-squares problem~(\ref{sms1d1}) may suggest a possible relation with the
GALS method (compare for example least-squares problem~(\ref{sms1d1}) with that
in~\cite[p.~327]{Roos-Stynes-Tobiska}). However, they are very different methods, since
for the examples in~Section~\ref{Se:experiments} in the present paper, SUPG and GALS methods
are identical (see e.g.~\cite{GALS}), and, as shown
in~Section~\ref{Se:experiments}, the SMS method and the SUPG method produce
markedly different results.
\end{remark}

\subsection{The multidimensional case}
\label{Se:multidimensional}

We will assume that $b(x)\ne 0$ for $x\in\Omega$, and that every characteristic (i.e., solution
of $dx/dt=b(x)$)  in~$\Omega$ enters and leaves~$\Omega$ in finite time. Also
for simplicity, we will assume that the domain~$\Omega$ in~(\ref{eq:model}--\ref{eq:modelbc}) has a
polygonal or polyhedral boundary.

Let~${\cal T}_h$ a triangulation of it, that is,  a partition
of the closure~$\overline\Omega$ of~$\Omega$ in $n$-simplices
with only faces and vertices in common. For every $\tau\in {\cal T}_h$, let ${\cal N}(\tau)$ be the
set of its vertices, and let~${\cal N}_h=\cup_{\tau\in{\cal T}_h}{\cal N}(\tau)$ be the set of vertices of~${\cal T}_h$. 

Similarly to the one-dimensional case, let~$X_h$ the space
of continuous piecewise linear polynomials.
We  express
$$
X_h=X_h^{-}\oplus V_h\oplus X_h^{+},
$$
where $\varphi_h=0$ on $\Gamma_D$ if~$\varphi_h\in V_h$, and for $\varphi_h$
in~$X_h^{-}$ (resp.~$X_h^{+}$),
if $\varphi_{h}(x)\ne 0$ for $x\in{\cal N}_h$, then $x\in \Gamma^{-}$ (resp.~$\Gamma_D
\backslash\Gamma^{-}$).
In the standard Galerkin method, first an element~$u_h^D\in X_{h}^{-}\oplus X_h^{+}$ is selected such
that the restriction ${u_h^D}_{\mid_{\Gamma_D}}$ to $\Gamma_D$ is a good approximation
to the Dirichlet data~$g_1$ in~(\ref{eq:modelbc}). This restriction is typically the interpolant or the
$L^2(\Gamma_D)$-projection onto the restrictions to~$\Gamma_D$ of~functions in~$X_h$.
Then, the Galerkin
approximation~$u_h\in u_h^D+V_h$ satisfies
\begin{equation}
a(u_h,\varphi_h)=(f,\varphi_h)+\varepsilon\langle g_2,\varphi_h\rangle_{\Gamma_N},
\qquad\forall\varphi_h\in V_h,
\label{gal_2d}
\end{equation}
where here and in the sequel, $\langle\cdot,\cdot\rangle_{\Gamma}$ denotes the~$L^2$ inner product
on $\Gamma\subseteq\partial\Omega$, and
$$
a(v,w)=\varepsilon (\nabla v,\nabla w)+(b\cdot\nabla v+cv,w),\qquad v,w\in H^1(\Omega).
$$

\begin{figure}[h]
$$
\includegraphics[height=3truecm]{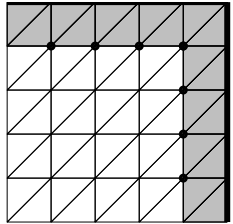}
$$
\caption{A triangulation of the unit square with $\Gamma^{0}\cup\Gamma^{+}$ marked with
a thicker line, the set $\Omega_h^{+}$ in grey and the the points of~${\cal N}_\delta$ with
circles.}
\label{fig_firstgammaplus}
\end{figure}

For the new method, similarly to the one-dimensional case, we consider
$$
L_h(v_h)=b\cdot\nabla v_h +c v_h.
$$
In order to describe the new approximation, we must set up the multidimensional version of the last
interval~$(x_{J-1},x_J)$ in the one-dimensional case.
For this purpose
we denote
$$
\Gamma_D^{0+}=(\Gamma^{+}\cup
\Gamma^{0})\cap\Gamma_D.
$$
We consider a  suitable $\Omega_h^{+}$ with $\Gamma_D^{0+}\subset\partial\Omega_h^{+}$ (to be specified below)
and let us denote by ${\cal N}_\delta$ and~${\mathbb N}_\delta$ the set of vertices
in~$\partial\Omega_h^{+}\backslash\Gamma_D$ and their indices, respectively, that is
$$
{\cal N}_\delta= {\cal N}\bigl(\partial\Omega_h^{+}\backslash\Gamma_D),\qquad
{\mathbb N}_\delta=\{j\in{\mathbb N} \mid x_j\in{\cal N}_\delta\},
$$
and by ${\mathbb R}^{{\mathbb N}_\delta}$ we will refer to the set of real vectors of the form
$(t_j)_{j\in{\mathbb N}_\delta}$.
The points in~${\cal N}_\delta$ and the set~$\Omega_h^{+}$ will play a role similar to that of~$x_{J-1}$ and~$(x_{J-1},x_J)$ in the one-dimensional
case. This can be seen in~Fig.~\ref{fig_firstgammaplus}, where,  shadowed in grey,
we show the set~$\Omega_h^{+}$
for a triangulation of the unit square for $b=[1,1]^T$, with $\Gamma_D=\partial\Omega$, so that $\Gamma^0\cup\Gamma^{+}$ is
consists of the sides $y=1$ and $x=1$. We show the points in~${\cal N}_\delta$
marked with circles.


The approximation $\tilde u_h\in u_h^D+V_h$, is then found by solving
the least-squares problem
\begin{equation}
\min_{\tilde u_h\in u_h^D+V_h, t\in{\mathbb R}^{{\mathbb N}_\delta}}\bigl\| L_h(\tilde u_h)-f
\bigr\|_{L^2(\Omega\backslash\Omega_h^{+})},\label{sms2d1}
\end{equation}
subject to the restriction
\begin{equation}
\label{sms2d2}
a(\tilde u_h,\varphi_h)+\sum_{j\in{\mathbb  N}_\delta} t_j \varphi_h(x_{j})=(f,\varphi_h)
+\varepsilon\langle g_2,\varphi_h\rangle_{\Gamma_N},
\qquad\varphi_h\in V_h.
\end{equation}
The optimality conditions of this problem can be written as the following linear problem:
find $\tilde u_h\in u_h^D+V_h$, $z_h\in V_h$ and~$t\in {\mathbb R}^{{\mathbb N}_\delta}$ such that
\begin{align}
\label{opt2d1}
\!\!\bigl(L_h(\tilde u_h),L_h(\varphi_h)\bigr)_{L^2(\Omega\backslash\Omega_h^{+})}
\! - a_h(\varphi_h,z_h) &=(L_h(\varphi_h),f)_{L^2(\Omega\backslash\Omega_h^{+})},\quad \varphi_h\in V_h,\\
z_h(x_j)& =0,\qquad j\in{\mathbb N}_\delta,
\label{opt2d2}\\
a(\tilde u_h,\varphi_h)+\sum_{j\in{\mathbb N}_\delta}t_j \varphi_h(x_{j})&=(f,\varphi_h)
+\varepsilon\langle g_2,\varphi_h\rangle_{\Gamma_N},
\quad\varphi_h\in V_h,
\label{opt2d3}
\end{align}
where, as in the one-dimensional case, $z_h$ is the Lagrange multiplier
of~restriction~(\ref{sms2d2}).



Let us specify now the set~$\Omega_h^{+}$. 
The obvious choice of setting~$\Omega_h^{+}=B_h$ where
$$
B_h=\bigcup_{\tau\cap(\Gamma_D^{0+})\ne \emptyset} \tau,
$$
as in Fig.~\ref{fig_firstgammaplus}, may lead to an unstable method 
if there are nodes $x\in {\cal N}_h$ interior to~$B_h$,  as it can be
easily checked in a numerical experiment. 
\begin{figure}[h]
$$
\includegraphics[height=3truecm]{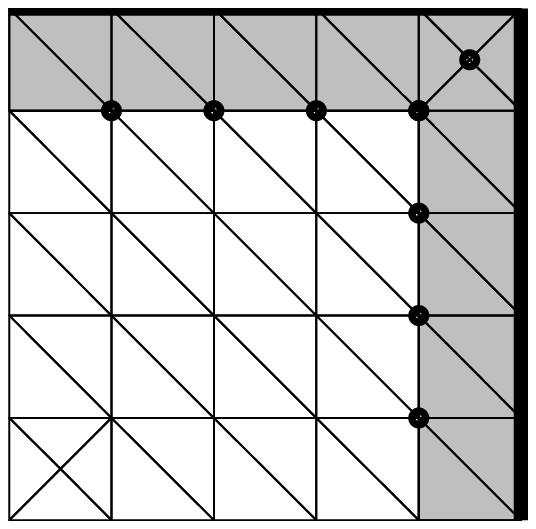}\qquad\qquad
\includegraphics[height=3truecm]{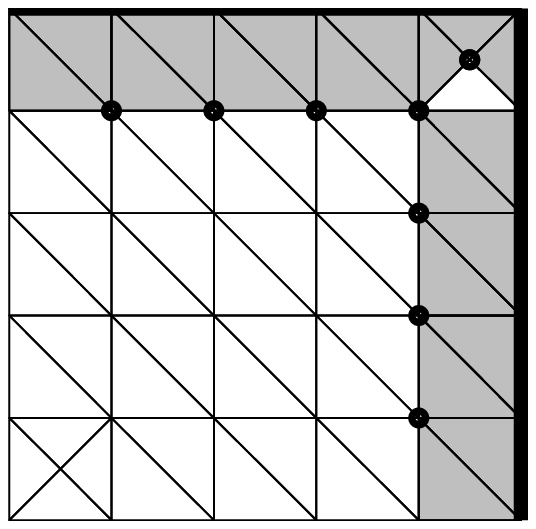}
$$
\caption{A triangulation of the unit square with $\Gamma^{0}\cup\Gamma^{+}$ marked with
a thicker line and  the points of~${\cal N}_\delta$ with
circles. In grey, the sets $B_h$ (left plot) and $\Omega_h^{+}$ (right plot) for $b=[1,1]^T$.}
\label{fig_secondgammaplus}
\end{figure}
To understand why in such a case the SMS method may be unstable, consider
the limit case~$\varepsilon=0$, $b$ constant and~$c=0$ (so that the bilinear form
$a$ is skew-symmetric) and consider a mesh as that depicted
in~Fig.~\ref{fig_secondgammaplus}, where there is one
node~$x_i$ interior to~$B_h$ who has only one neighbour node~$x_k
\in {\cal N}_\delta$. Then, on the
one hand, the basis function $\varphi_i$ of the node~$x_i$ satisfies that
$a_h(\varphi_i,\varphi_h)=0$ for all~$\varphi_h\in V_h$, except when
$\varphi_h$ is a multiple of~the basis function~$\varphi_k$ of the node~$x_k$.
Furthermore, since the support of~$\varphi_i$ is contained in~$B_h$, we have that 
$\left\| L_h(\varphi_i)\right\|_{L^2(\Omega\backslash B_h)}=0$. Consequently,
taking $\tilde v_h=\varphi_i$, $z_h=0$,  $t_l=0$ for $l\ne k$ and $t_k=-a(\varphi_i,\varphi_k)$
we have a nontrivial solution of
\begin{align}
\label{opt2d1_hom}
\bigl(L_h(\tilde v_h),L_h(\varphi_h)\bigr)_{L^2(\Omega\backslash\Omega_h^{+})}
 - a_h(\varphi_h,z_h) &=0,\qquad \varphi_h\in V_h,\\
z_h(x_j)& =0,\qquad j\in{\mathbb N}_\delta,
\label{opt2d2_hom}\\
a(\tilde v_h,\varphi_h)+\sum_{j\in{\mathbb N}_\delta}t_j \varphi_h(x_{j})&=0,
\qquad\varphi_h\in V_h,
\label{opt2d3_hom}
\end{align}
if we set $\Omega_h^{+}=B_h$. It can be easily checked that this is also the
case if the node~$x_i$ interior to~$B_h$ is connected to more that one
node in~${\cal N}_\delta$. For $\varepsilon >0$, one can expect a unique solution
but, as we have checked in practice, it is a solution where
$\tilde u(x_i)=O(1/\varepsilon)$.

To avoid
this situation we remove from~$B_h$ the upwind triangle of any interior point. More precisely,
for every node~$x_i\in{\cal N}_h$, let us denote by
$\tau_{-}(x_i)$ its upwind triangle, that is
$$
\tau_{-}({x_i})=\{\tau\in {\cal T}_h \mid \{x_i-\lambda b\}\cap \tau \ne \emptyset,\quad \lambda
\rightarrow 0\}.
$$
Then we define
$$
\Omega_h^{+}=B_h\backslash\bigl(
\bigcup_{x_i\in{\cal N}_h\cap\stackrel{\circ}{\widehat{B}_h}}\tau_{-}(x_i)\bigr).
$$
If there are two upwind triangles because $x_j-\lambda b$ is and edge for
$\lambda$ smaller than some~$\lambda_0$, we may select the first one of them
after ordering all the elements.

With this choice of~$\Omega_h^{+}$, as we will show in Section~\ref{Se:elOmega_h},
the few cases were there are nontrivial solutions
of~(\ref{opt2d1_hom}--\ref{opt2d3_hom}) have an easy fix from the computational point of view.
We also note that a slightly different choice of~$\Omega_h^{+}$ is taken
in~Example~\ref{ej:Silvester} in Section~\ref{Se:experiments}, where the vector field~$b$ has closed
integral curves.


\subsection{Further Extensions}

In the same way that the Shishkin meshes are not restricted to the standard Galerkin
method, the new technique, originally motivated by Shishkin meshes, is not restricted
to the Galerkin method either. In this section we comment on how to use it with some
other methods of widespread use, such as  the SUPG, GALS, LPS and CIP methods.

In these methods the
approximation~$u_h\in u_h^D+V_h$, instead of satisfying~(\ref{gal_2d}),
satisfies
\begin{equation}
a_h(u_h,\varphi_h)=(f,\varphi_h)_h
+\varepsilon\langle g_2,\varphi_h\rangle_{\Gamma_N},
\qquad\forall\varphi_h\in V_h,
\label{other_2d}
\end{equation}
where $a_h$ and~$(\cdot,\cdot)_h$ are mesh-dependent bilinear forms. For the
SUPG method they are given by
\begin{align}
\label{SUPG-bilinear}
a_h(v_h,\varphi_h)&=a(v_h,\varphi_h)+\sum_{\tau\in{\cal T}_h} \delta_{\tau}
(L_{\mid_\tau}(v_h),b\cdot\varphi_h)_{\tau}, \\
 (f,\varphi_h)_h&=(f,\varphi_h)+
\sum_{\tau\in{\cal T}_h} \delta_{\tau}
(f,b\cdot\varphi_h)_{\tau},
\label{SUPG-rhs}
\end{align}
where $\delta_{\tau}$ is and adequately chosen parameter and
$L_{\mid_{\tau}}$ denotes $L$ restricted to~$\tau$. In the GALS method
the term~$b\cdot\nabla\varphi_h$ is replaced by~$L_{\mid_{\tau}}(\varphi_h)$
(see e. g., \cite{Roos-Stynes-Tobiska} for a full description of these methods).
The extension of the SMS technique to these methods is as follows:
 we solve the least-squares problem~(\ref{sms2d1}) but we replace
the restriction~(\ref{sms2d2}) by
\begin{equation}
\label{sms2dSUPG}
a_h(\tilde u_h,\varphi_h)+\sum_{j\in{\mathbb  N}_\delta}t_j \varphi_h(x_{j})=(f,\varphi_h)_h
+\varepsilon\langle g_2,\varphi_h\rangle_{\Gamma_N}.
\qquad\varphi_h\in V_h.
\end{equation}

In the numerical experiments in the present paper, we only consider the SMS for the SUPG method.
To distinguish it from the SMS method in the previous section, we will refer to them
as {\it Galerkin-based} and {\it SUPG-based} SMS methods. 

We could also consider the extension to higher-order finite-element methods. The extension
seems straightforward
since it consists in replacing the approximation space $V_h$ by piecewise quadratics or pieceswise
cubics (together with $L_h$ in~(\ref{L_h_higher})). Numerical experiments (not reported here)
indicate that this straightforward extension does not give as good results as the piecewise
linear elements, at least for one-dimensional problems. The reason seems to be the need to redefine
the set~$\Omega_{h}^+$ for higher-order elements. This will be subject of future studies.


\section{Analysis of the SMS method}
\label{Se:analisis}

\subsection{The one-dimensional case}
\label{Se:analisis1d}
Some understanding of why the new method is able to suppress or, at least, dramatically reduce the
spurious oscillations of standard approximations can be gained by analyzing the
problem~(\ref{simple1d}--\ref{simple1dbc}) when $b$ is a positive constant and $c=0$. We do this
first for the Galerkin method.

\medskip

\noindent{\it The Galerkin method}.

 We first consider the limit case $\varepsilon=0$. In this case, it is easy to check that the Galerkin
approximation satisfies the equations
\begin{equation}
\frac{b}{2}(u_{j+1}-u_{j-1})=f_j \equiv\int_{x_{j-1}}^{x{j+1}} f(x)\varphi_j(x)\,dx, \qquad
j=1,\ldots,J-1.
\label{gal1d_eqs}
\end{equation}
Thus, summing separately odd and even-numbered equations we obtain the following expressions,
\begin{align}
u_{2j}&=u_0+\frac{2}{b}\sum_{i=1}^j f_{2i-1},\qquad j=1,\ldots,(J'-1)/2,
\label{loseven}\\
 u_{2j-1}&= u_{J'}-\frac{2}{b}
\sum_{i=j}^{(J'-1)/2} f_{2i},\qquad
j=1,\ldots,(J'-1)/2,
\label{losodd}
\end{align}
where, here and in the rest of this section
$$
J'=\left\{ \begin{array}{lcl} J,&\quad&\hbox{\rm if~$J$ is odd},\\
J-1,&\quad&\hbox{\rm if~$J$ is even}.\end{array}\right.
$$

We notice that when $J$ is odd,
the expressions in~(\ref{loseven}--\ref{losodd}) are the discrete counterparts
of the problems
\begin{align}
bu_x=f,&\quad 0<x<1,\qquad u(0)=0,\label{avant}\\
 bu_x=f,&\quad  0<x<1,\qquad u(1)=0,\label{patras}
\end{align}
which, unless $f$ has zero mean, have different solutions. Thus, since for sufficiently smooth
$f$ the expressions~(\ref{loseven}--\ref{losodd}) are consistent with~(\ref{avant}--\ref{patras}),
oscillations in the numerical approximation are bound to occur whenever $f$ is not of zero mean.
For $J$ even, it is easy to check that the Galerkin equations~(\ref{gal1d_eqs}) have no solution
unless $f_{1}+f_3+\cdots+f_{J-1}=0$, in which case the solution is not unique. Thus for $J$ even
the Galerkin method is not stable.

However, for $J$ odd, the Galerkin method is stable in the following sense
\begin{equation}
\left\| u_h\right\|_\infty \le \frac{2}{b} \left\|f\right\|_{-h},
\label{stab2}
\end{equation}
where
\begin{equation}
\left\|f\right\|_{-h} = \left\| S_h\right\|_\infty,
\label{fh1}
\end{equation}
and
\begin{equation}
\label{fh2}
S_{2j}=\sum_{i=1}^j f_{2i-1},\qquad S_{2j-1}=\sum_{i=j}^{(J'-1)/2} f_{2i},\qquad
j=1,\ldots,(J'-1)/2.
\end{equation}


\medskip

\noindent{\it The SMS technique}.
Here the equations are
\begin{equation}
\frac{b}{2}(\tilde u_{j+1}-\tilde u_{j-1})+\alpha\delta_{J-1}(j)=f_j \qquad
j=1,\ldots,J-1,
\label{smsan1}
\end{equation}
where $\delta_{J-1}$ denotes Dirac's delta function
$$
\delta_{J-1}(j)=\left\{ \begin{array}{lcl} 0,&\quad & \hbox{\rm if $j\ne J-1$,}\\ 1,&&
\hbox{\rm if $j=J-1$.}\end{array}\right.
$$

As Proposition~\ref{prop_stab} below shows, the SMS method is stable independently of the parity of the grid. To prove stability we will consider  the  auxiliary function $q_h\in V_h$ whose nodal values are
\begin{equation}
q_{j}=-\frac{1-(-1)^j}{b},\qquad j=1,\dots,J-1.
\label{elq1b}
\end{equation}
It is easy to check that when $J$ is odd $q_h$
satisfies
\begin{equation}
\label{laq_h1d}
a(q_h,\varphi_h)=\varphi_h(x_{J-1}),\qquad\varphi_h\in V_h,
\end{equation}
whereas when $J$ is even it satisfies 
\begin{equation}\label{rev1}
a(q_h,\varphi_h)=0,\qquad\varphi_h\in V_h.
\end{equation}

\begin{proposition}\label{prop_stab}
 There exist a positive constant $C$ such that for $\varepsilon=0$ the SMS
approximation
$\tilde u_h$ satisfies
\begin{equation}
\label{stab_sms_prop}
\left\| \tilde u_h\right\|_\infty \le \frac{6}{b}\Bigl( \left\|f\right\|_{-h}+\frac{h}{6J}\left|r\right|
\Bigr),
\end{equation}
where $r=(f,L_h(q_h))_{L^2(0,x_{J-1})}$.
\end{proposition}
\begin{proof} Recall that the SMS approximation~$\tilde u_h$ is part of the solution
$(\tilde u_h,\alpha,z_h)$ of system~(\ref{opt1d1}--\ref{opt1d3}).
Due to~(\ref{laq_h1d}), (\ref{rev1}) and~(\ref{opt2}) we have that
$a(q_h,z_h)=0$, so that taking~$\varphi_h=q_h$ in~(\ref{opt1d1}) we have
\begin{equation}
\label{opt1d1q}
\bigl(L_h(\tilde u_h),L_h(q_h)\bigr)_{L^2(0,x_{J-1})}=
(f,L_h(q_h))_{L^2(0,x_{J-1})}.
\end{equation}
On the other hand,
we notice that $L(q_h)=2(-1)^{j}/h_j$, for~$x\in(x_{j-1},x_j)$ and $j=1,\ldots,J-1$,
so that equation~(\ref{opt1d1q}) becomes
\begin{equation}
\label{smsan1b}
2\sum_{j=1}^{J-1} (-1)^j b\frac{\tilde u_j-\tilde u_{j-1}}{h_j} = 2\sum_{j=1}^{J-1} (-1)^j \frac{1}{h_j}
\int_{x_{j-1}}^{x_j} f(x)\,dx.
\end{equation}

To prove~(\ref{stab_sms_prop}) we treat the cases $J$ odd and $J$ even separately. For
$J$ odd, we notice that
since the equations in (\ref{smsan1}) are those of the Galerkin method except for the last one, 
applying the stability estimate (\ref{stab2}) we have
\begin{equation}
\left\| \tilde u_h\right\|_\infty \le \frac{2}{b}( \left\|f\right\|_{-h}+\left|\alpha\right|).
\label{stab_sms0}
\end{equation}
Also, since $\tilde u_h=u_h-\alpha q_h$, from~(\ref{smsan1b}) it follows that
$$
\biggl( 4\sum_{j=1}^{J-1} \frac{1}{h_j}\biggr)\alpha=
2b\sum_{j=1}^{J-1} (-1)^j \frac{ u_j -u_{j-1}}{h_j}- r,
$$
$r$ being the right-hand side of~(\ref{smsan1b}).
Thus, applying again the stability estimate~(\ref{stab2}), we have
\begin{equation}
\left|\alpha\right| \le 2\left\|f\right\|_{-h}+\frac{h}{4(J-1)}\left|r\right|.
\label{stab1b}
\end{equation}
From this inequality, (\ref{stab_sms0}) and~the fact that $2(J-1)\ge J$, the stability estimate~(\ref{stab_sms_prop}) follows.

For $J$ even, taking $\varphi_h=q_h$ in~(\ref{opt1d3}) and applying~(\ref{rev1})
we have
$$
0=(b(\tilde u_h)_x -f,q_h)+\alpha q_{J-1}=\alpha q_{J-1}+\frac{2}{b}\sum_{j=1}^{J/2} 
f_{2j-1},
$$
so that
\begin{equation}
\alpha=\sum_{j=1}^{J/2} 
f_{2j-1}.
\label{alpha_par}
\end{equation}
Also, since the equations in (\ref{smsan1}) are those of the Galerkin method except for the last one, 
which does not appear in~(\ref{loseven}) and~(\ref{losodd}) for $J$ even,
we have that for $J$ even (\ref{loseven}) and~(\ref{losodd}) also hold with $u_h$ replaced
by~$\tilde u_h$. Thus, it follows that
\begin{equation}
\label{est1c}
\left\| \tilde u_h\right\| \le \frac{2}{b}\left\| f\right\|_{-h}+ \left| \tilde u_{J-1}\right|.
\end{equation}

The value of~$\tilde u_{J-1}$ is then obtained from~(\ref{smsan1b}) by replacing $\tilde u_j$ by
their values as expressed in~(\ref{loseven}--\ref{losodd}), which gives
$$
-\biggl(2b\sum_{j=1}^{J-1}\frac{1}{h_j}\biggr) \tilde u_{J-1}+4 \sum_{j=1}^{J-2}
\Bigl( \frac{1}{h_j}+\frac{1}{h_{j+1}}\Bigr) S_j = r,
$$
and, thus,
\begin{equation}
\left| \tilde u_{J-1}\right| \le \frac{4}{b}\left\|f\right\|_{-h} + \frac{h}{2b(J-1)} \left|r\right|,
\label{est1c2}
\end{equation}
which, together with~(\ref{est1c}) shows that~(\ref{stab_sms_prop}) also
holds for $J$ even.
\end{proof}


Once the stability is proved, we investigate the convergence of the SMS method by considering the basic solution,
$$
bu_x=f,\qquad u(0)=0,
$$
and the error
$$
\tilde e_h=\tilde u_h - I_h(u),
$$
where $I_h(u)$ is the interpolant of~$u$ in~$X_h$. We prove different estimates depending
on whether
\begin{equation}
\label{asymp_uni}
\max_{1\le j\le J-1} \left| h_{j+1}-h_{j-1}\right|= Ch^2,
\end{equation}
holds for some constant $C>0$.

\begin{proposition} There exist a positive constant $C$ such that for $\varepsilon=0$ the error
$e_h=\tilde u_h-I_h(u)$ of SMS
approximation
satisfies the following estimates:
\begin{equation}
\label{prop1_eq1}
\bigl\| \tilde e_h\bigr\|_\infty \le Ch\left\|f'\right\|_{L^1(0,1)},
\end{equation}
and, if~(\ref{asymp_uni})
holds,
\begin{equation}
\label{prop1_eq2}
\bigl\| \tilde e_h\bigr\|_\infty \le Ch^2(\left\|f'\right\|_{\infty}+\left\|f''\right\|_{L^1(0,1)}).
\end{equation}
\end{proposition}

\begin{proof}
We have
\begin{equation}
\frac{b}{2}(\tilde e_{j+1}-\tilde e_{j-1})+\Bigl(\alpha-\frac{b}{2}u(1)\Bigr)\delta_{J-1}(j)=\tau_j,\qquad
j=1,\ldots,J-1,
\label{eq:error1}
\end{equation}
where
$$
\tau_j= \int_{x_{j-1}}^{x_{j+1}}f(x)\Bigl(\varphi_j-\frac{1}{2}\Bigr)\,dx.
$$
Since
$$
\int_{x_{j-1}}^{x_j}(bI_h(u)_x-f)\,dx=0,\qquad j=1,\ldots,J-1,
$$
so that $(b(\tilde q_h)_x,bI_h(u)_x-f)_{L^2(0,x_{J-1})}=0$, it follows that
\begin{equation}
(b(\tilde q_h)_x,b(\tilde e_h)_x)_{L^2(0,x_{J-1})}=0.
\label{eq:error1bb}
\end{equation}
Thus, applying~(\ref{stab_sms_prop}) with $f$ replaced by~$\tau=[\tau_1,\ldots,\tau_{J-1}]^T$ and
$r$ by~$0$, we have
\begin{equation}
\label{err1b}
\left\| \tilde e_h\right\|_\infty \le \frac{6}{b}\left\| \tau_h\right\|_{-h}.
\end{equation}
Integration by parts reveals
$$
\tau_j =- \int_{x_{j-1}}^{x_{j+1}} f'(x)\Bigl(\int_{x_j}^x (\frac{1}{2}-\varphi_j)\,dy\Bigr)\,dx,
$$
and, further,
$$
\tau_j=f'(x_{j+1})\frac{h_{j+1}^2}{12} -f'(x_{j-1})\frac{h_{j}^2}{12} +
 \int_{x_{j-1}}^{x_{j+1}} f''(x)\Bigl(\int_{x_j}^x\int_{x_j}^y (\frac{1}{2}-\varphi_j)\,dz\,dy\Bigr)\,dx.
$$
Then, it is easy to check that
\begin{equation}
\left\| \tau \right\|_{-h} \le Ch\left\|f'\right\|_{L^1(0,1)},
\label{trucsms1da}
\end{equation}
or
\begin{align}
\left\| \tau\right\|_{-h} \le &
 \max_{k\le (J-1)/2}\biggl(
\biggl| \sum_{j=1}^k f'(x_{2j-1})\frac{h_{2j}^2 -h_{2j-1}^2}{12}\biggr|
+
\biggl|\sum_{j=k}^{(J-1)/2}\!\! f'(x_{2j})\frac{h_{2j+1}^2 -h_{2j}^2}{12}\biggr|\biggr)
\nonumber \\
&{}+Ch^2\bigl(\left\|f'\right\|_{L^\infty(0,1)}+\left\|f''\right\|_{L^1(0,1)}\bigr).
\label{trucsms1d}
\end{align}
We notice that~from~(\ref{err1b}) and~(\ref{trucsms1da}) the estimate~(\ref{prop1_eq1}) follows.
On the other hand, if (\ref{asymp_uni}) holds, then the right-hand side of~(\ref{trucsms1d})
is $O(h^2)$, so that~(\ref{err1b})
and~(\ref{trucsms1d}) imply~(\ref{prop1_eq2}).
\end{proof}

Let us finally comment the case $\varepsilon>0$.
In this case, equations~(\ref{smsan1}) should be modified by adding
the term
$$
\varepsilon\frac{\tilde u_{j}-\tilde u_{j+1}}{h_{j+1}}+\varepsilon\frac{\tilde u_{j}-\tilde u_{j-1}}{h_{j}}
$$
to the left hand side. Then, the stability estimate~(\ref{stab_sms_prop}) holds but with
$\left\|f\right\|_{-h}$
on the right-hand side replaced by
$$
\left\|f\right\|_{-h}+4J\varepsilon \max_{1\le j\le J}\frac{1}{h_j}\left\|\tilde u_h\right\|_\infty,
$$
added to the right-hand side. Thus, we have the same stability
estimate~(\ref{stab_sms_prop}) with
the factor~$6/b$ replaced by~$12/b$ whenever
\begin{equation}
\label{epsilon_small}
\varepsilon<\frac{b}{48J}\min_{1\le j\le J} h_j.
\end{equation}
Finally, with respect to the
error, standard asymptotic analysis (see e. g. \cite[Theorem~1.4]{Roos-Stynes-Tobiska})
shows that
the solution $u$ of~(\ref{simple1d}-\ref{simple1dbc}) and the basic solution~$u_0$ (solution
of~(\ref{avant})) satisfy that
$$
\max_{0\le x\le 1}\bigl| u(x)-u_0(x)+u_0(1)e^{b(x-1)/\varepsilon}\bigr| \le C\varepsilon,
$$
for some $C>0$ independent of $\varepsilon$.
Since $e^{b(x-1)/\varepsilon} < \varepsilon$ for $x<1-\varepsilon\left|\log(\varepsilon)\right|/b$,
then, whenever
\begin{equation}
\label{cota_h_J}
h_J\ge \frac{\varepsilon\left|\log(\varepsilon)\right|}{b},
\end{equation}
the bounds~(\ref{prop1_eq1}) and~(\ref{prop1_eq2}) hold if~(\ref{epsilon_small}) also holds.
\medskip

\begin{remark}\label{re:chapuza} As mentioned in the Introduction, in the present paper we are
concerned with large P\'eclet numbers. For this reason and for simplicity, we have not
pursued conditions which allow for larger values of~$\varepsilon$ than~(\ref{epsilon_small}).
Obviously, a more detailed (and lengthy) analysis would allow to obtain more general conditions
on~$\varepsilon$. 
\end{remark}


\subsection{The multidimensional case}
\label{Se:analysis_2d}
\subsubsection{Uniqueness of solutions}
\label{Se:elOmega_h}

In this section we limit ourselves to study the possibility loosing uniqueness in the new method
when $\varepsilon=0$ and
how to overcome this in practice.  These cases are an indication of when to expect the method
to perform poorly for small $\varepsilon>0$. In the rest of the paper, we denote
$$
\hat \Omega_h=\Omega\backslash\Omega_h^{+}.
$$

We start by characterizing the solutions of~(\ref{opt2d1_hom}--\ref{opt2d3_hom}) when $\varepsilon=0$.
As we
now show,  they are
given by those elements $\tilde v_h\in V_h$ satisfying
\begin{equation}
\left\|L_h(\tilde v_h)\right\|_{L^2(\hat\Omega_h)}=0.
\label{Ltildev_h=0}
\end{equation}
To see this, we take
$\varphi_h=\tilde v_h$ in~(\ref{opt2d1_hom}) and get
$$
\left\|L_h(\tilde v_h)\right\|_{L^2(\hat\Omega_h)}^2+a(\tilde v_h,z_h)=0.
$$
But taking  $\varphi_h=z_h$ in~(\ref{opt2d3_hom}) we have
$$
a(\tilde v_h,z_h)=-\sum_{j\in{\mathbb N}_\delta} t_jz_h(x_j)=0,
$$
where in the last equality we have applied~(\ref{opt2d2_hom}). Thus, (\ref{Ltildev_h=0}) follows.
Conversely, if~(\ref{Ltildev_h=0}) holds,  this implies $L_h(\tilde v_h)=0$ on~$\hat \Omega_h$.
But since we are assuming that $\varepsilon=0$, we also have that $a(v_h,w_h)=(L_h(v_h),w_h)$
for any~$v_h,w_h\in V_h$. Thus,
if~for a basis function $\varphi_j$ we have $0\ne a(\tilde v_h,\varphi_j)=(L_h(\tilde v_h),\varphi_j)$,
it follows that its
support cannot be entirely inside~$\hat \Omega_h$, so that $x_j$ must be in~${\cal N}_\delta$.
Then taking~$z_h=0$ and $t_j=-a(\tilde v_h,\varphi_j)$ we have a solution
of~(\ref{opt2d1_hom}--\ref{opt2d3_hom}).

We devote the rest of this section to comment on
a case that may easily arise in practice where  there is a~$\tilde v_h$ not entirely null in~$\Omega$ that satisfies~(\ref{Ltildev_h=0})
(and, hence, the lost of uniqueness in the SMS method when $\varepsilon=0$).
We will also comment
on what possible remedies can
be applied in practice in these cases. In the rest of the section we assume $c=0$, that is,
$L_h(w_h)=b\cdot\nabla w_h$.

\begin{figure}[h]
$$
\includegraphics[height=2.8truecm]{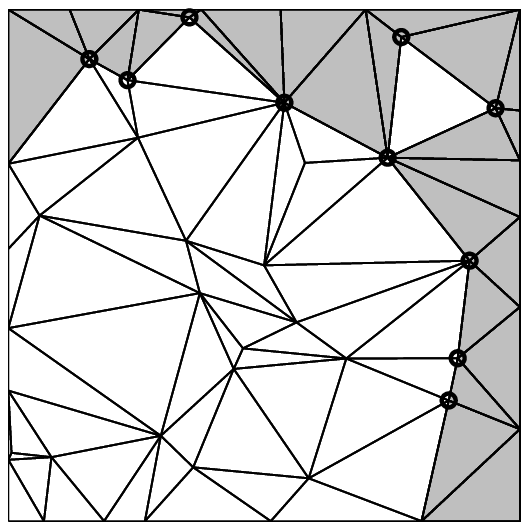}\qquad\qquad
\includegraphics[height=2.8truecm]{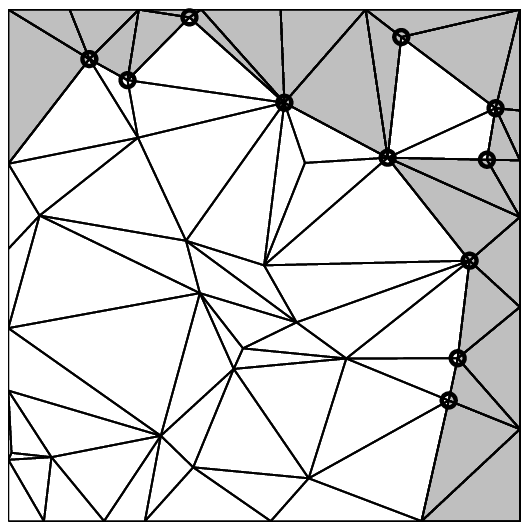}
$$
\caption{Two triangulations of the unit square with $\hat \Omega_h$ of disconnected interior.  The set $\Omega_h^{+}$
for $b=[1,1]^T$ is shadowed in grey  and points of~${\cal N}_\delta$ are marked with
circles.}
\label{fig_isolated}
\end{figure}

Consider an example like those in~Fig.~\ref{fig_isolated}, that is, when the interior set
$\mathring{\hat\Omega{}}_h$ of~$\hat\Omega_h$
is not connected.
We may have $\tilde v_h\in V_h$ not null in~$\Omega$ and yet satisfying~(\ref{Ltildev_h=0}).
To see why, 
take for example the left plot in Fig.~\ref{fig_isolated} and assume that $b$ is constant with strictly positive components for simplicity. Let $x_j$ and~$x_k$ be the most downwind vertices of the isolated
triangle~$\tau$ of~$\hat \Omega_h$. Observe that since their associated basis functions
$\varphi_j$ and~$\varphi_k$ have linearly independent gradients,  it is possible to find real values
$\alpha$ and~$\beta$ so that $\alpha \varphi_j+\beta \varphi_k\ne 0$ but 
$b\cdot(\alpha \nabla \varphi_j+\beta\nabla \varphi_k)=0$. Setting then
$\tilde v_h=\alpha \varphi_j+\beta \varphi_k$ we have that $\tilde v_h\ne 0$ but it
satisfies~(\ref{Ltildev_h=0}).
If $b$ is not constant, then it is possible to find~$\tilde v_h$ such that
$\left\|v_h\right\|_{L^2(\Omega)}=1$
but 
$\left\|L_h( v_h)\right\|_{L^2(\hat\Omega_h)}\approx\hbox{\rm diam}{(\tau)}^2$; there
may not be nontrivial solutions of~(\ref{opt2d1_hom}--\ref{opt2d3_hom}) but the method may not
be very accurate since residuals of size $\hbox{\rm diam}{(\tau)}^2$ allow for errors of
size~$1=\left\| v_h\right\|_{L^2(\Omega)}$. 
Similar arguments show that this is also the case
of the right plot in Fig.~\ref{fig_isolated}. These arguments are easily extended to tetrahedra
in three-dimensional domains. In Section~\ref{Se:appendix} we show some other cases where there
are nontrivial $\tilde v_h$ satisfying~(\ref{Ltildev_h=0}). We also state general conditions to
prevent them.

\begin{figure}[h]
$$
\includegraphics[height=2.5truecm]{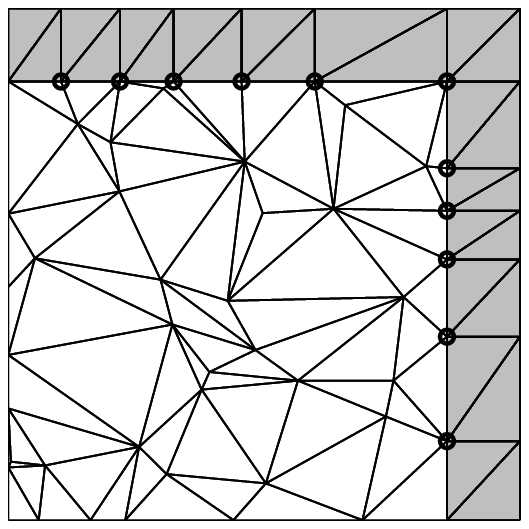}\qquad\quad
\includegraphics[height=2.5truecm]{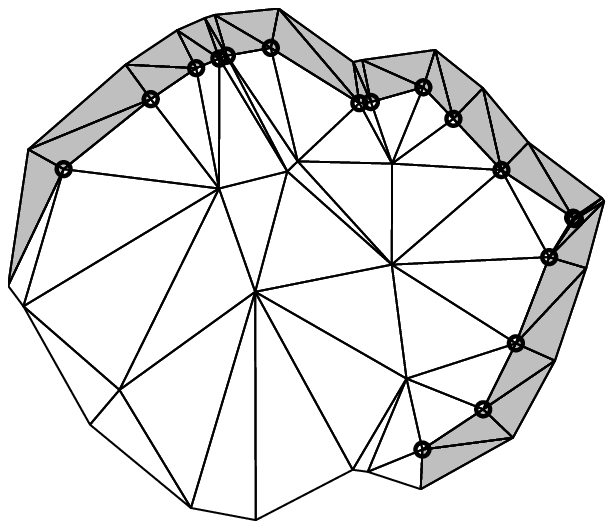}\qquad\quad
\includegraphics[height=2truecm]{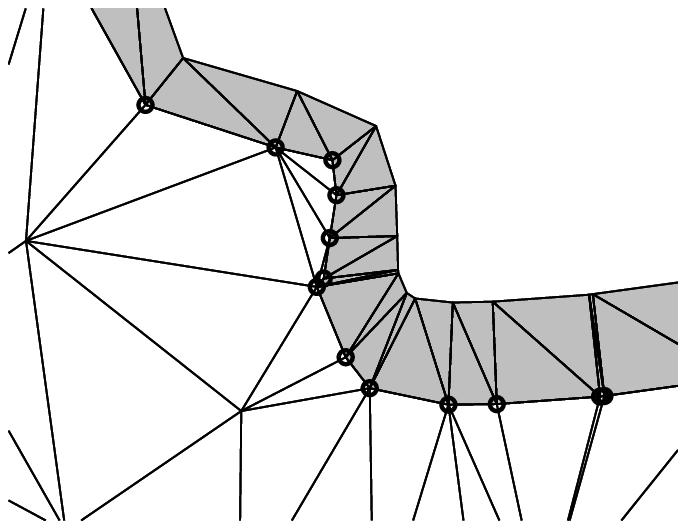}
$$
\caption{Triangulations with strip of elements on~$\Gamma_D^{0+}$ to prevent isolated components
in~$\mathring{\hat\Omega{}}_h$. The set $\Omega_h^{+}$
for $b=[1,1]^T$ is shadowed in grey.}
\label{fig_strip}
\end{figure}
\begin{figure}[h]
$$
\includegraphics[height=2.6truecm]{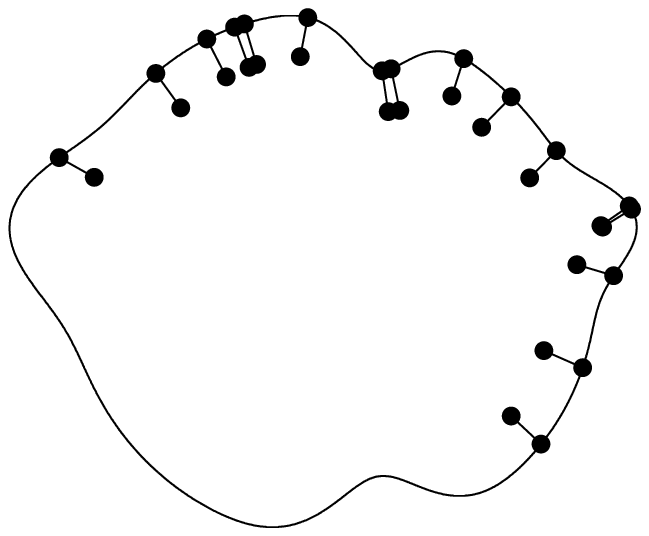}\qquad
\includegraphics[height=2.6truecm]{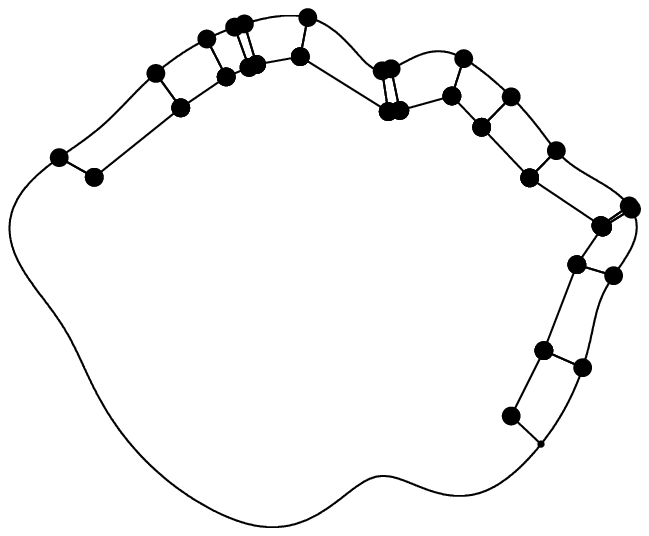}\qquad
\includegraphics[height=2.6truecm]{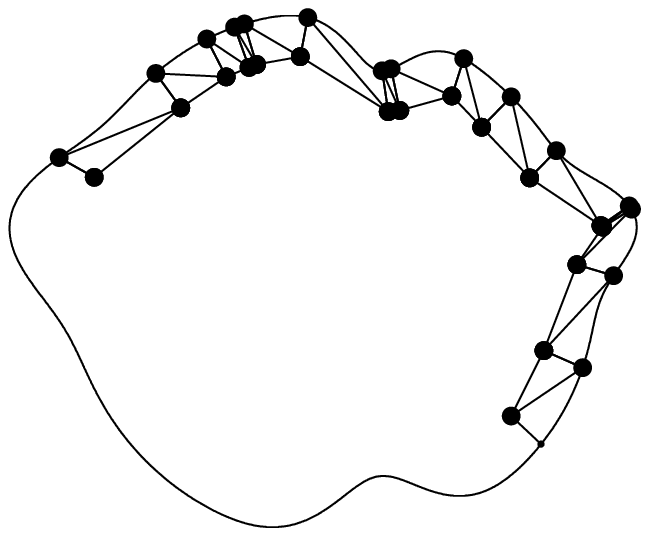}
$$
\caption{The process to build a strip of elements along~$\Gamma_D^{0+}$.}
\label{fig_strip2}
\end{figure}

From the practical point of view, all those conditions are satisfied
if the meshes are designed with a strip of elements on $\Gamma_D^{0+}$ as indicated
in~Fig.~\ref{fig_strip}.
One possible way to form this strip is to take the nodes on the outflow boundary
(which we assume partitioned into segments or faces) and displace them
along the normal vector to obtain the interior nodes of the strip (see Fig.~\ref{fig_strip2}).
Then, sides or faces in the outflow boundary are replicated in the displaced nodes, so that
rectangles are obtained in two-dimensional problems and triangular prisms in three-dimensional ones. In the two-dimensional case, 
the triangulation of the strip is obtained by joining two opposite
vertices of the rectangles.
In three dimensions, tetrahedra can be obtained for example by adapting to triangular prisms
the so-called Kuhn's triangulation of the cube, \cite{Bey}, \cite{Freudenthal}, \cite{Kuhn}.
This is a partition of the unit cube into 6 tetrahedra, all of them having the origin  and the vertex of
coordinates $(1,1,1)$,
as common vertices. 
It is a simple exercise to adapt
Kuhn's triangulation of the cube to triangular prisms.

We notice, however, that 
there are times in practice when one cannot choose the grid because it is part
of a larger problem, or because it is designed to satisfy some other constraints, etc.
In Section~\ref{Se:appendix}, we comment on a different procedure
 to avoid isolated
components in~$\hat\Omega_h$. 

\subsubsection{Further cases of lack of uniqueness}
\label{Se:appendix}
We complete the study of the previous section
on the cases where there are nontrivial~$\tilde v_h\in V_h$
satisfying~(\ref{Ltildev_h=0}). We assume $\varepsilon=0$  and $c=0$
here as well. We start by checking
that only $\tilde v_h=0$ satisfies~(\ref{Ltildev_h=0}) when
$\mathring{\hat\Omega{}}_h$ is disconnected due to an
element upwind of a node interior to~$B_h$, as it was the case
discussed in~Fig.~\ref{fig_secondgammaplus}. Let us take that example again.
We first notice that since we are assuming
that $\varepsilon=0$, then, (\ref{Ltildev_h=0}) imply that $\tilde v_h$ is constant along
characteristics, and  since $\tilde v_h=0$ on~$\Gamma_{-}$, it must vanish on the big component
of~$\hat \Omega_h$. This implies that $\tilde v_h$ also vanish on the most upwind vertex $x_k$ of the
isolated triangle~$\tau$. It also vanish on the vertex on the boundary. So it is only on the remaining
vertex~$x_i$ where it may not vanish. But recall that the triangle~$\tau$ is upwind of $x_i$, so that
$b\cdot \nabla\varphi_i\ne 0$. Thus, $L_h(\tilde v_h)=0$ on~$\tau$ implies that $\tilde v_h=0$
on~$\tau$.

\begin{figure}[h]
$$
\includegraphics[height=3truecm]{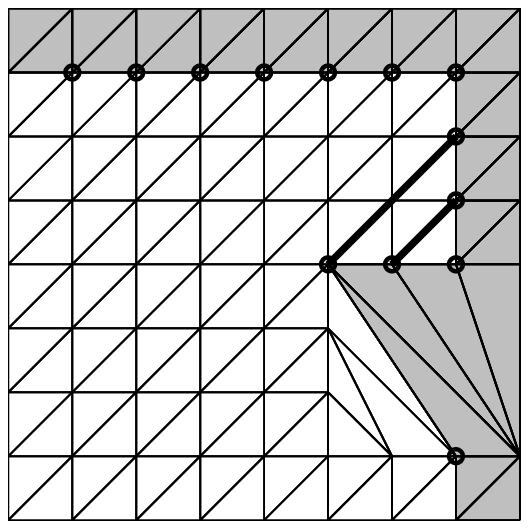}\qquad\qquad
\includegraphics[height=3truecm]{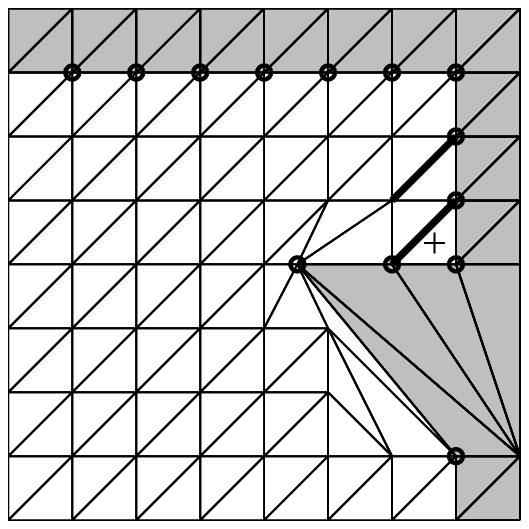}\qquad\qquad
\includegraphics[height=3truecm]{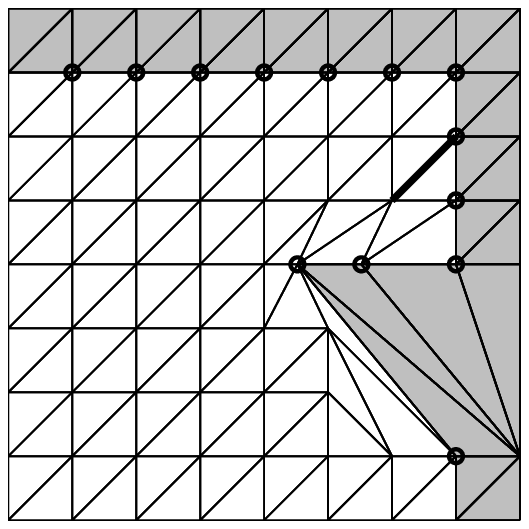}
$$
\caption{Triangulations of the unit square with triangles in~$\hat\Omega_h$ downwind
of~$\Omega_h^{+}$ and with one side (thicker line) parallel to the wind velocity
$b=[1,1]^T$. The set $\Omega_h^{+}$ is shadowed in grey.}
\label{fig_shadowed}
\end{figure}

However, even when~$\hat \Omega_h$ has a connected interior, there
may be nontrivial $\tilde v_h\in V_h$ satisfying~(\ref{Ltildev_h=0}) when there are
triangles (resp.\ tetrahedra) in~$\hat\Omega_h$, downwind of~$\Omega_h^{+}$ and with one side aligned
with the (constant) wind velocity~$b$ (resp.\ one face parallel to~$b$).
This is the case of the triangles with one side plotted
with a thicker line in~the grids in~Fig.~\ref{fig_shadowed}
for $b=[1,1]^T$. 
The basis function of the vertices opposite those sides have their gradients orthogonal
to~$b$. For example, for the triangle $\tau$ marked with a $+$ sign
in the center plot of~Fig.~\ref{fig_shadowed}, let $x_i$ be the vertex opposite to the shadowed side on the triangle,
and $x_j$ and~$x_k$ the other two vertices. Since $b\cdot \nabla \varphi_i=0$ in~$\tau$,
and $\varphi_i=0$
on the rest of~$\hat\Omega_h$, then
$\tilde v_h=\varphi_i$ satisfies (\ref{Ltildev_h=0}). For the grid in the left plot,
we may take $\tilde v_h=\alpha
\varphi_i + \beta(\varphi_j+\varphi_k)$ for
any $\alpha,\beta\in{\mathbb R}$, since $b\cdot \nabla \tilde v_h=0$ in the part of $\hat \Omega_h$
downwind~$\Omega_h^{+}$ and $\tilde v_h$ is
null in the rest of~$\hat\Omega_h$. As Proposition~\ref{prop:shadowed} below shows,
there is no nontrivial $\tilde v_h\in \Omega$ satisfying~(\ref{Ltildev_h=0}) for the triangulation on the
right-hand side of~Fig.~\ref{fig_shadowed}. 

\begin{figure}[h]
$$
\includegraphics[height=2.8truecm]{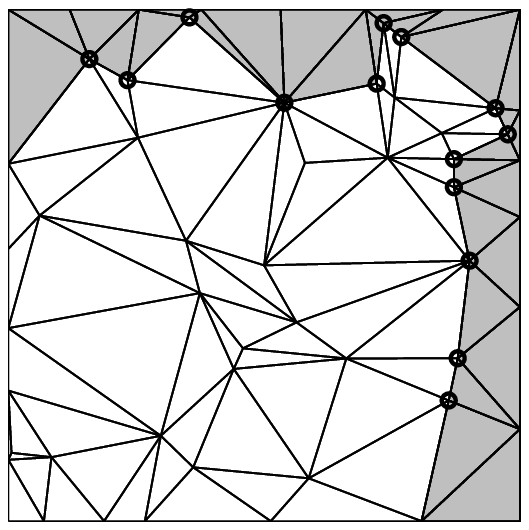}\qquad\qquad
\includegraphics[height=2.8truecm]{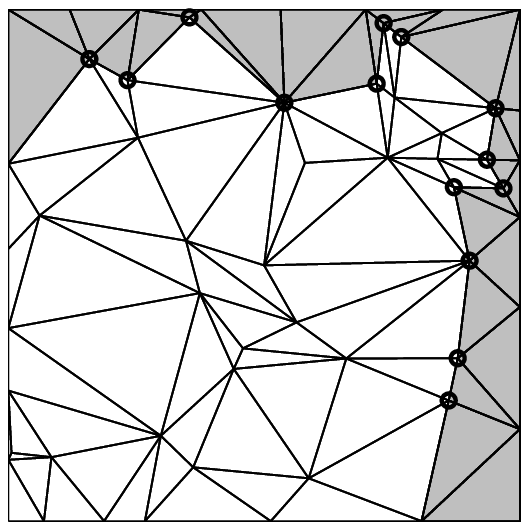}
$$
\caption{The triangulations in Fig.~\ref{fig_isolated} after regular refinement of~upwind triangles
closest to isolated components.  The new set $\Omega_h^{+}$
for $b=[1,1]^T$ is shadowed in grey  and points of the new~${\cal N}_\delta$ are marked with
circles.}
\label{fig_isolated_r}
\end{figure}

\begin{figure}[h]
$$
\includegraphics[height=3truecm]{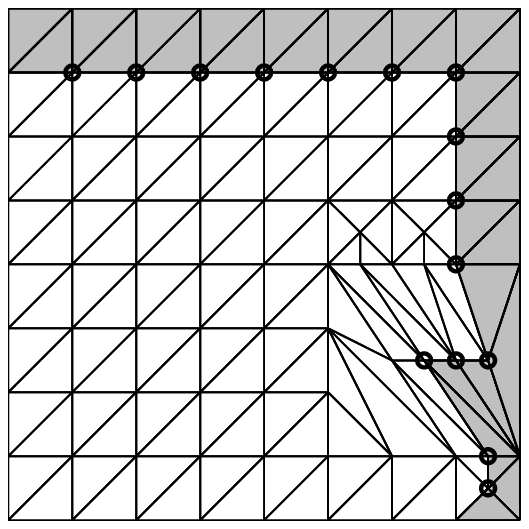}\qquad\qquad
\includegraphics[height=3truecm]{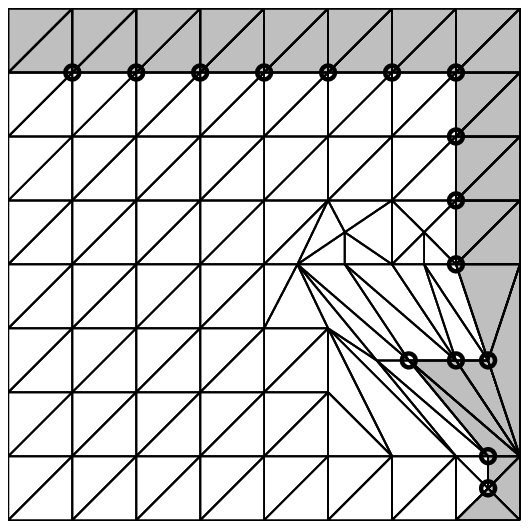}
$$
\caption{The triangulations in Fig.~\ref{fig_shadowed} after regular refinement of triangles
in~$\Omega_h^{+}$
upwind of~$\hat\Omega_h$
The set~$\Omega_h^{+}$ for $b=[1,1]^T$ is shadowed in grey and points of the new~${\cal N}_\delta$
are marked with circles.}
\label{fig_shadowed_r}
\end{figure}

Nevertheless, notice that as it happens with the isolated components in~Fig.~\ref{fig_isolated},
 the cases depicted
in~Fig.~\ref{fig_shadowed} cannot occur if one designs the mesh with a strip of elements
along~$\Gamma_D^{0+}$. 
However, as commented at the end of Section~\ref{Se:elOmega_h},
there are times in practice when one cannot choose the grid so that cases like those
depicted in~Fig.~\ref{fig_isolated} or~in~Fig.~\ref{fig_shadowed} may occur. We now comment
on the possible remedies to avoid nontrivial solutions
of~(\ref{opt2d1_hom}--\ref{opt2d3_hom}). For example, one can connect isolated components
of~$\mathring{\hat\Omega{}}_h$
by enlarging~$\hat \Omega_h$ with triangles from~$\Omega_h^{+}$. We have checked in practice
that this may leave large portions of $\Gamma_{D}^{0+}$ as boundary of~$\hat \Omega_h$, and this
resulted in the presence of spurious oscillations in the computed approximations.
A better alternative is
to enlarge the grid by refining those elements in~$\Omega_h^{+}$
upwind of the isolated component.
In Fig.~\ref{fig_isolated_r} we show the result of dividing into four similar triangles (red or
regular refinement)
the two upwind neighbours of the isolated components in Fig.~\ref{fig_isolated} and using
longest edge bisection in those triangles that inherit a hanging node.
We see that after the refinement process $\mathring{\hat\Omega{}}_h$ is connected.
Also, 
the lack of uniqueness induced by
those triangles
in~$\hat\Omega_h$ with a side parallel to~$b$ downwind of~$\Omega_h^{+}$ can be prevented
by red-refining their upwind triangles in~$\Omega_h^{+}$, as it can be
seen in~Fig.~\ref{fig_shadowed_r},
where no side parallel to the wind velocity~$b$ is now downwind of~$\Omega_h^{+}$.

The following result states general conditions guaranteeing that (\ref{Ltildev_h=0}) implies~$\tilde
v_h=0$.
Consider the set of points that are not downwind of~$\Omega_h^{+}$, that is,
$$
d(\hat\Omega_h)=\{ x\in \mathring{\hat\Omega{}}_h\; \mid \{ x-tb\mid t>0\}\cap \Omega_h^{+}=
\emptyset\}.
$$

\begin{proposition}\label{prop:shadowed} Let $b$ be constant and $c=0$. Then if either
$d(\hat\Omega_h)=\mathring{\hat\Omega{}}_h$ or for any $x\in
\mathring{\hat\Omega{}}_h\backslash\overline{ d(\Omega_h)}$ there is a
path in~$\mathring{\hat\Omega{}}_h$ from $x$
to a point $y\in d(\Omega_h)$ through elements with no edge or face parallel
to~$b$, then the only element $\tilde v_h\in V_h$ satisfying~(\ref{Ltildev_h=0}) is $\tilde v_h=0$.
\end{proposition}

\begin{proof}
 If $\tilde v_h\in V_h$ satisfies~(\ref{Ltildev_h=0}), then $\tilde v_h=0$
on~$d(\hat\Omega_h)$, since
any~$\tilde v_h$ satisfying~(\ref{Ltildev_h=0}) is constant along the characteristics $x+tb$ which,
eventually intersect~$\Gamma^{-}$ where $\tilde v_h=0$. If~$d(\hat\Omega_h)=
\mathring{\hat\Omega{}}_h$, then $\tilde v_h$ can only be nonzero on~$\Omega_h^{+}$. But elements
on~$\Omega_h^{+}$ have vertices either on~$\Gamma_D^{0+}$ where $\tilde v_h=0$
or on~$\hat\Omega_h$.
Thus, $\tilde v_h$ must be zero on every element in~$\Omega_h^+$.

Suppose now that, $d(\hat\Omega_h)\ne \mathring{\hat\Omega{}}_h$. Then, for $x\in \mathring{\hat\Omega{}}_h\backslash \overline{d(\hat\Omega_h)}$, let $\gamma$ be the path 
in~$\mathring{\hat\Omega{}}_h$ connecting $x$ to a
point~$y\in d(\hat\Omega_h)$, not intersecting any element with a side or face parallel to~$b$.
There is no loss of generality in assuming the path to be a polygonal, and except maybe from
the first and last segments, the remaining segments joint the barycenters or arithmetic means of the vertices
of the elements
it intersects. There will be a first element~$\tau$ of those intersected by~$\gamma$
which is not entirely inside~$\mathring{\hat\Omega{}}_h\backslash \overline{d(\hat\Omega_h)}$.
Since $\tau$ has no edge or face parallel to~$b$, then its
interior $\mathring{\tau{}}\;$ is not entirely contained in~$\mathring{\hat\Omega{}}_h\backslash \overline{d(\hat\Omega_h)}$, and thus, $\tilde v_h$ must vanish on~$\tau$. But then, the previous
element is in a similar situation, with $\tilde v_h$ vanishing on the side or face in common with~$\tau$
and this side or face not being parallel to~$b$, so that $\tilde v_h$ is zero in that element too. Repeating
the argument we conclude that $\tilde v_h$ vanish in~$x$.
\end{proof}

\begin{figure}[h]
$$
\includegraphics[height=3truecm]{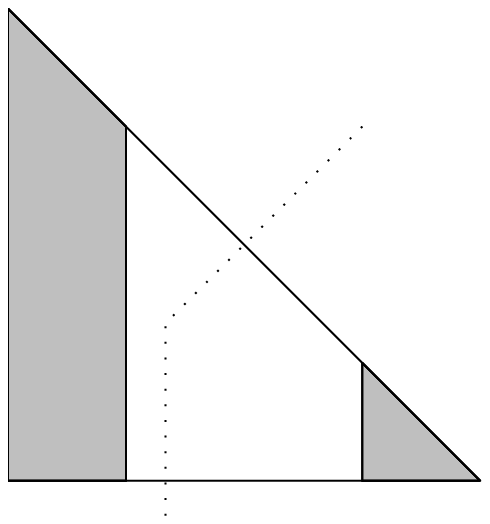}\qquad\qquad
\includegraphics[height=3truecm]{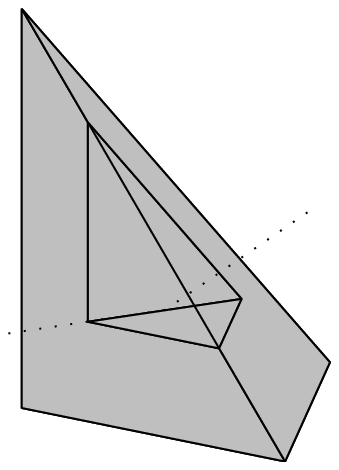}
$$
\caption{Two and three-dimensional elements showing shadowed in grey the part
where $\Omega_h^{+}$ is restricted to be after two refinements where all sides and
edges are bisected. Dotted lines are path entering and leaving the element through
the only sides and faces not on the boundary of~$\Omega$ and passing through the
barycenter.}
\label{fig_casos}
\end{figure}

Finally, for those situations in practice where one has to work with a mesh without a strip of
elements along~$\Gamma_D^{0+}$, we now study how to make robust the technique of
red-refining adequate triangles in~$\Omega_h^{+}$. We first notice that it is not difficult
to design examples 
where
red-refining even {\em all} triangles in~$\Omega_h^{+}$ does not make
$\mathring{\hat\Omega{}}_h$ connected.
However, we now argue that two red-refinements are enough to connect any
isolated component. This is done as in the proof of~Proposition~\ref{prop:shadowed}
by considering a polygonal path in~$\Omega$ joining the isolated component with
with another component 
upwind of it. The segments of
the path may be assumed to join barycenters (resp.\ arithmetic means of the vertices)
of triangles with a common
side (resp.\ tetrahedra with a common face). Consider a refinement strategy where
all sides or edges are bisected, and apply it to all elements in~$\Omega_h^+$ intersecting the
path. It is not difficult to analyze all possible cases, the worst cases being
those where all vertices are in~$\Gamma_D^{0+}$ and all sides (resp.\ faces)
except those two intersecting the path are also on~$\Gamma_D^{0+}$. After two refinements
the elements in~$\Omega_h^+$ are restricted to be in the area or volume induced by
the new vertices of the second refinement closest to the initial vertices,
as we show in~Fig.~\ref{fig_casos}.
Then, in the two-dimensional case,
the original path is completely embedded
in~$\mathring{\hat\Omega{}}_h$, and, in the three-dimensional case, it is on the
border of~$\mathring{\hat\Omega{}}_h$, so that a small change moves it
into~$\mathring{\hat\Omega{}}_h$.

Let us mention that in the case of tetrahedra, the so-called regular refinement is not
unique~\cite{Bey}, and that the arguments above also apply if the refinement is done by
successive bisection of the six edges of the tetrahedron using, for example, the technique
in~\cite{Arnold}, or, in the case of triangles, three applications of longest-edge bisection
to the three sides of a triangle.


\section{Numerical Experiments}
\label{Se:experiments}

In this section we solve~(\ref{eq:model}--\ref{eq:modelbc}) on different domains~$\Omega$
and with different forcing terms and vector fields~$b$. In all cases we take $c=0$.  We used
Matlab 7.13.0 and the backslash command to solve linear systems. Also, for
the SUPG method, we take the streamline diffusion parameter as suggested in
formulae~(5-7) in~\cite{John-Knobloch-2007}, which we reproduce here for the convenience
of the reader. More precisely, in the SUPG method test functions on element~$\tau$ are of the form
$
\varphi_h+\delta_{\tau}b\cdot \nabla\varphi_h$, where $\delta_\tau$ is given by
\begin{align}
\label{delta_tau1}
\delta_{\tau}&=\frac{\hbox{\rm diam}(\tau,b)}{2\left|b\right|}, \qquad \hbox{\rm if Pe}_{\tau}>1,\\
\delta_{\tau}&=\frac{\hbox{\rm diam}(\tau,b)^2}{4\varepsilon}, \qquad \hbox{\rm if Pe}_{\tau}\le1,
\label{delta_tau2}
\end{align}
where $
\hbox{\rm Pe}_{\tau} = \frac{\left|b\right|\hbox{\rm diam}(\tau,b)}{2\varepsilon}$, and,
if  $\varphi_1,\ldots,
\varphi_{d+1}$ are the basis functions in element~$\tau$ (taking value 1 on one of the vertices and 0
on the rest of them) then
$$
\hbox{\rm diam}(\tau,b)=\frac{2\left|b\right|}
{\left|b\cdot\nabla \varphi_1\right|+\cdots+\left|b\cdot\nabla \varphi_{d+1}\right|}.
$$
Here $\left|b\right|$ stands for the euclidean norm of the vector field $b$. If $b$ is not constant, it
is evaluated at the barycenter of element~$\tau$.
The Matlab
codes used in this section are available from the author on request (check also the url
address
{\tt http://personal.us.es/bga/bga\_files/software\_bga.html}). Excluded are the codes
corresponding to Example 7 which is subject of current resarch.
\medskip

\begin{example}\label{ej:comparem}
{\it Simulation of  Shishkin grids}. Since the SMS method was conceived
as a simulation of a Shishkin grid, we now check how well it does it.
We take $\Omega=(0,1)^2$, 
$b=[2,3]^T$,
Dirichlet homogeneous boundary conditions and the forcing term $f$ such that the solution is
$$
u(x,y)=\bigl(x-e^{2(x-1)/\varepsilon}\bigr)\bigl(y^2-e^{3(y-1)/\varepsilon}\bigr).
$$
This example is taken from Example~5.2 in~\cite{Volker1}, but in our case $c=0$. On the
one hand, we solve
this problem with the SUPG method on Shishkin meshes with $J=2N$ subdivision in each
coordinate direction. They are formed as tensor products of one-dimensional Shishkin
meshes, with values $\sigma_x=2\varepsilon\log(N)$ and~$\sigma_y=(3/2)\varepsilon\log(J)$
on the $x$ and $y$ directions, respectively. On the other hand,
we take the coarse part of the Shishkin mesh, which is a triangulation of
\begin{equation}
\Omega_\sigma=[0,1-\sigma_x]\times [0,1-\sigma_y],
\label{Omegasigma}
\end{equation}
with $N$ subdivisions in each coordinate direction, and apply the SMS method
(both Galerkin-based and SUPG based).
For $N=5,10,20,\ldots,320$, we compute the numerical approximations
and measure both the computing time and the $L^\infty$ errors in the interior points of the coarse part of the Shishkin grid
(points which are shared by the three methods)
for $\varepsilon=10^{-4}$ and $\varepsilon=10^{-8}$.
Results are shown in~Fig.~\ref{fig_comparemv}.
\begin{figure}[h]
$$
\includegraphics[height=5.3truecm]{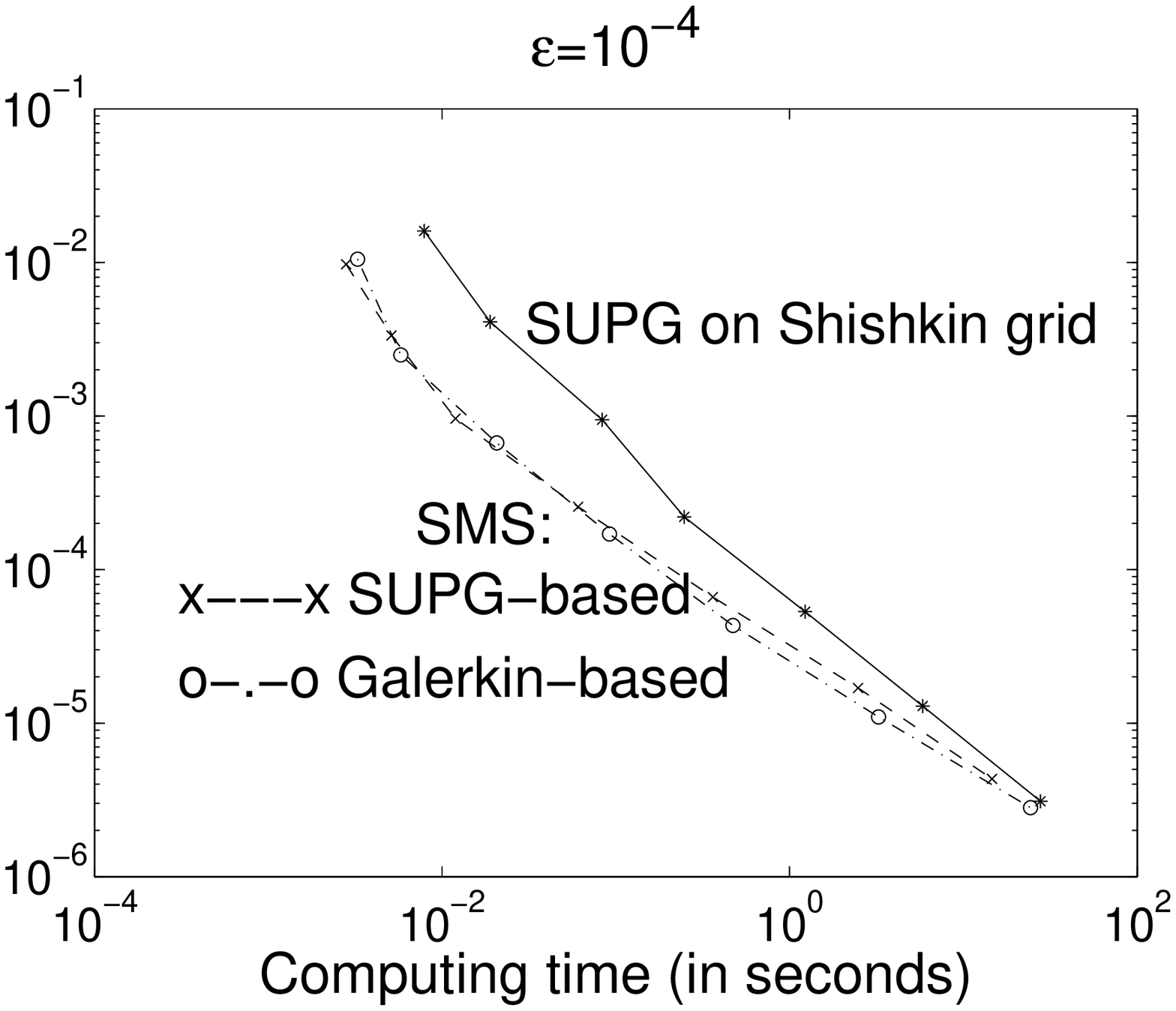}\quad
\includegraphics[height=5.3truecm]{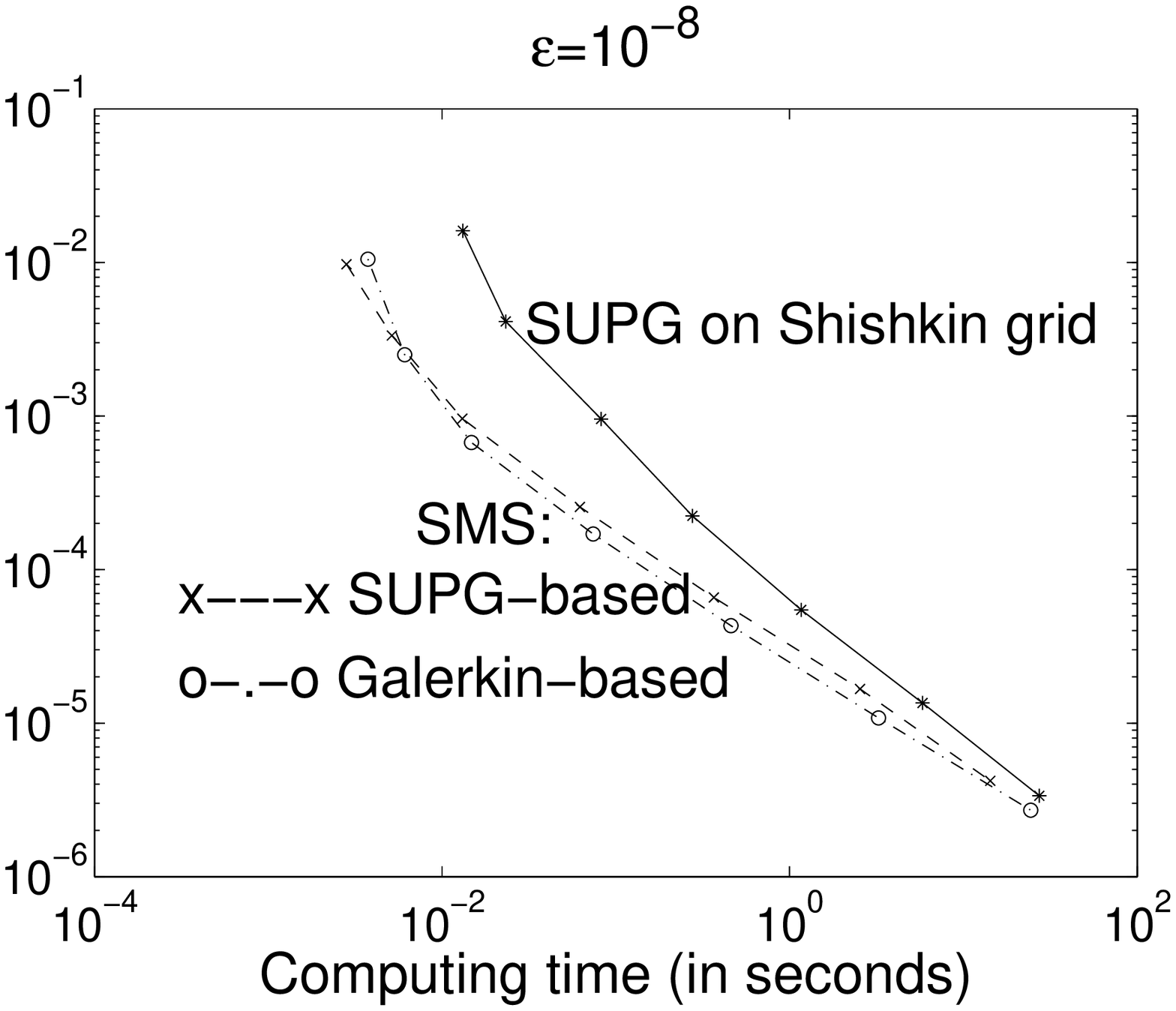}
$$
\caption{Relative efficiency of SMS and SUPG methods in Example~\ref{ej:comparem}.}
\label{fig_comparemv}
\end{figure}
Looking at the errors, and for both values of $\varepsilon$, we notice that the
three methods commit roughly the same errors, being in this example those of the SUPG-based
SMS method slightly worse in general.
In terms of computational efficiency, though, both SMS methods are roughly equally efficient,
and markedly better than
the SUPG method on Shishkin grids, since although the three methods commit roughly the
same errors, the SMS methods compute the approximations between 4 and 2 times faster
than the SUPG on the Shishkin grid (only for $N=320$, the 
SMS methods were only 1.2 times faster than the SUPG method). In the next example, we further
comment on computational costs and the structure of the linear systems to be solved to obtain
the different approximations.
\end{example}
\medskip

\begin{example}\label{ej:comparem_same} {\it Comparisons on the same grid}. In the previous example,
it may be considered unfair to compare the new methods on an $N\times N$ coarse grid with the SUPG on a $2N\times2N$ grid, since this is bound to be more expensive. In the present example
we apply the methods on the same grids (uniform with diagonals running Southwest-Northeast). The convection-problem is the same as in the previous example, and, as before, $L^\infty$ errors
are measured in interior mesh points. The SUPG method in this example gave very poor
results, with errors above $10^{-1}$. For this reason, and following suggestions
in~\cite{Madden-Stynes-96},  we programmed it with some crosswind diffusion.
More precisely, if in the SUPG method test functions on element~$\tau$ are of the form
$
\varphi_h+\delta_{\tau}b\cdot \nabla\varphi_h$, where $\delta_\tau$ is the stabilization
parameter, we used~$\varphi_h+\delta_{\tau}(b\cdot \nabla\varphi_h+\delta_c\partial_x\varphi)$
where the value of~$\delta_c$ was set by trial and error to obtain the
smallest errors. These values were
$\delta_c=0.7701$,    $0.8783$~and    $0.9365$, for $N=10$, $20$ and~$40$ subdivisions
in each coordinate direction. For similar reasons, the value of~$\delta_\tau$
in~formulae~(5-7) in~\cite{John-Knobloch-2007}, was multiplied
by $1.57$, $1.615$ and~$1.64$ for the above-mentioned number of subdivisions, respectively.

The results can be seen in~Fig.~\ref{fig_comparemsame}. Even though the SMS methods can be
up to
twice as costly as the SUPG on the same grid, the errors are, however, between 13 and~60 times
smaller (SUPG-based SMS) and 26 and~130 (Galerkin-based SMS),
so that the SMS method is computatinally much more efficient.
\begin{figure}[h]
$$
\includegraphics[height=5.3truecm]{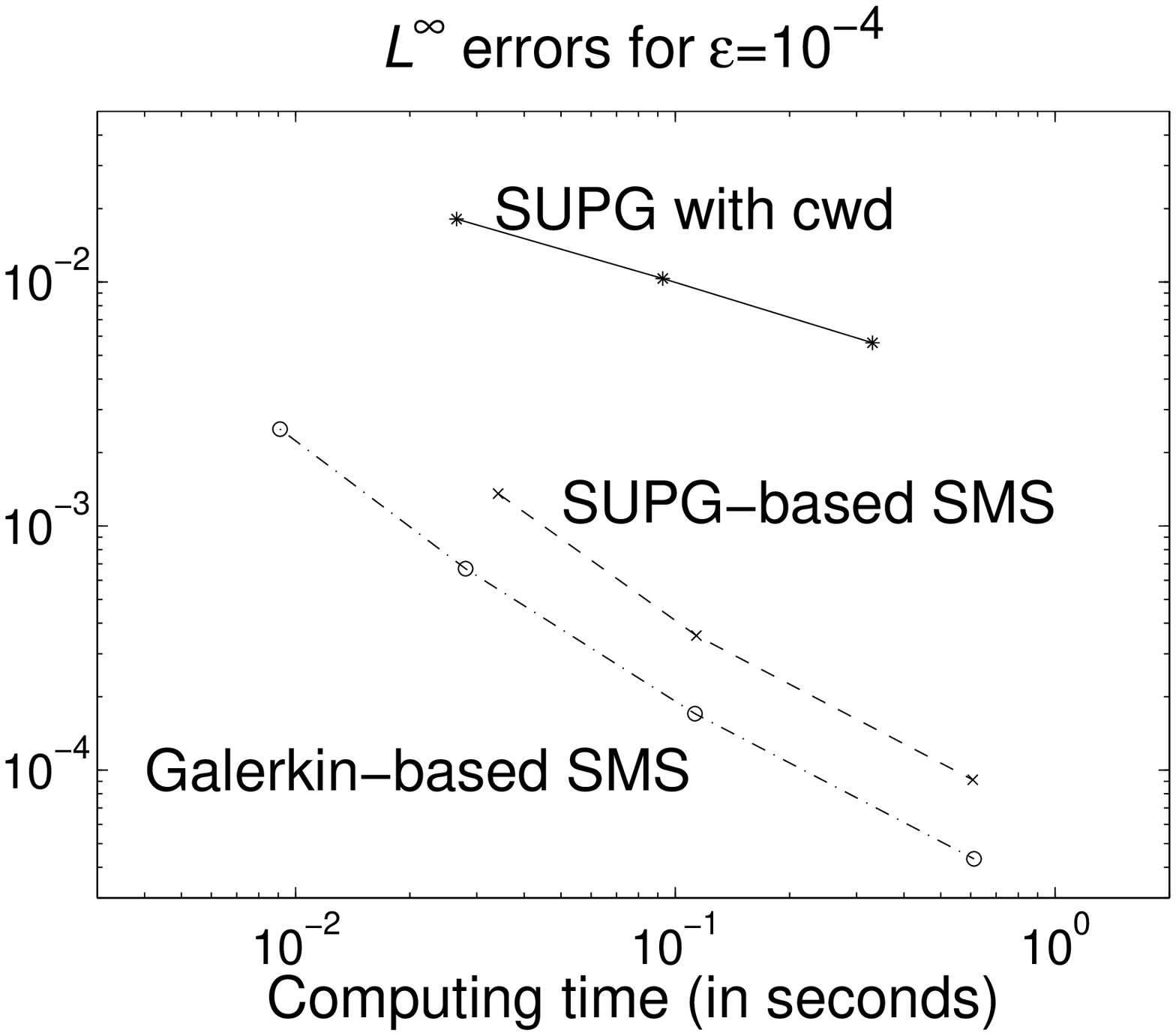}\quad
\includegraphics[height=5.3truecm]{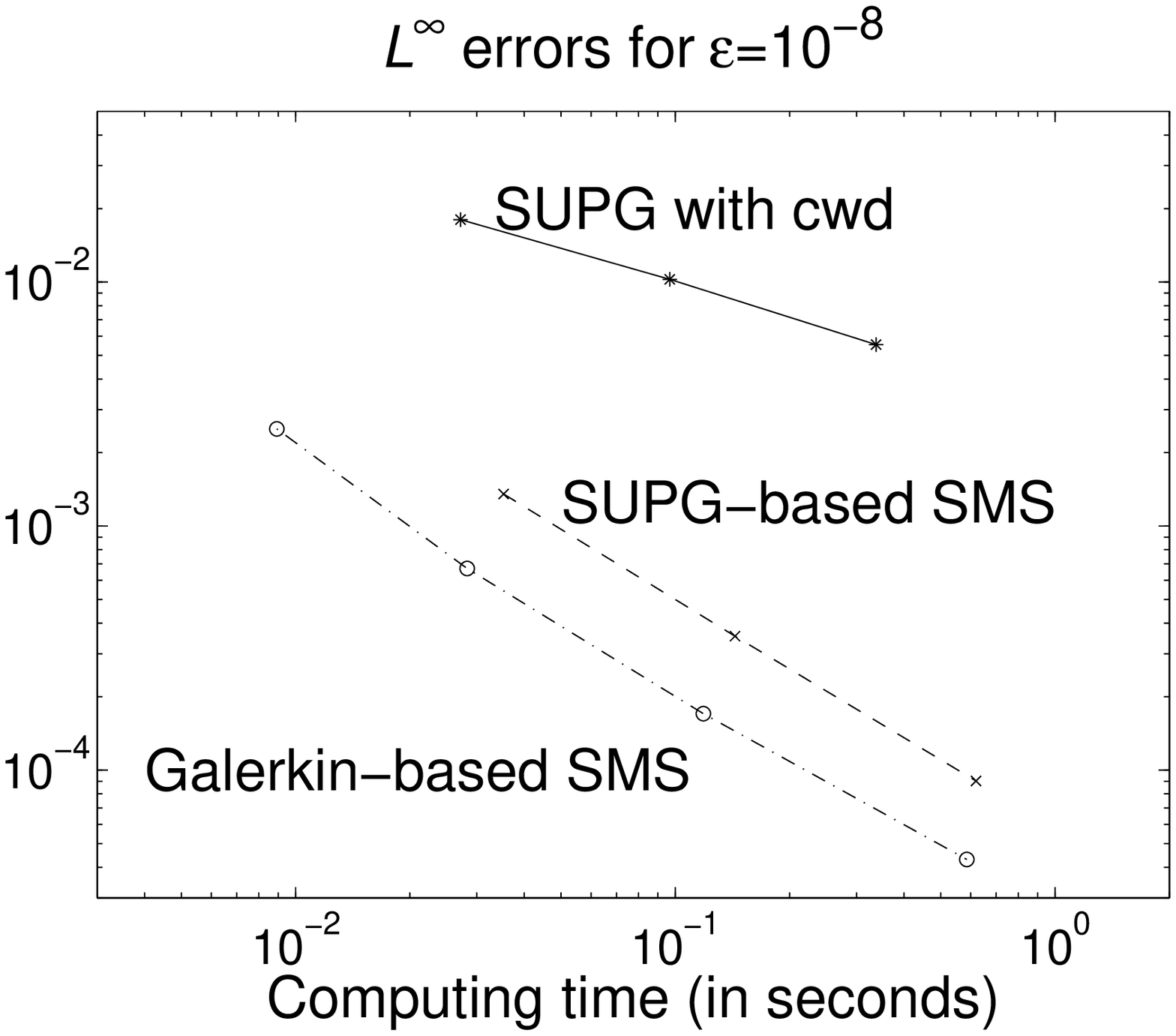}
$$
\caption{Relative efficiency of SMS and SUPG (with crosswind diffusion) methods in Example~\ref{ej:comparem_same}.}
\label{fig_comparemsame}
\end{figure}

Let us comment here on the computational cost of the methods. Recall that the
SMS approximation $\tilde u_h$ is obtained (together with the Lagrange multipliers~$z_h$
and the values $t_j$, $j\in{\mathbb N}_\delta$ as the solution of the optimality
conditions~(\ref{opt2d1}-\ref{opt2d3}).
Let $\varphi_1,\ldots,\varphi_n$ be the nodal basis
of $V_h$ (each $\varphi_i$ takes value 1 on one single node of the triangulation and 0 on the rest of them), and let $A$ and~$S$ be the $n\times n$ matrices with entries
$$
a_{i,j}=a(\varphi_i,\varphi_j),\qquad s_{i,j}=(L(\varphi_i),L(\varphi_j))_{L^2(\Omega\backslash\Omega_h^{+})},\qquad 1\le i,j\le n,
$$
(or $a_{i,j}=a_h(\varphi_i,\varphi_j)$ in the case of the SUPG method) respectively,
and let $E$ be the $n\times m$ matrix whose columns are those of the $n\times n$ identity
matrix corresponding to the indexes in~${\mathbb N}_\delta$. Notice that we may assume
$m\approx\rho n^{(d-1)/d}$ for some
$\rho>0$, ($d$ being the dimension of the euclidean space where the domain~$\Omega$ is). 
Observe also that A, S, and $E$
are typically sparse matrices. The nodal values of~$\tilde u_h$,
the values $t_j$, $j\in {\mathbb N}_\delta$
and the nodal values of~$z_h$ are then obtained by solving a linear system
whose coefficient matrix is
\begin{equation}
M=\left[\begin{array}{ccc} S& 0& -A^T\\ 0& 0& -E^T\\ A& E& 0\end{array}\right].
\end{equation}
This must compared with the SUPG method where the coefficient matrix is~$A$, that is,
a system of order~$n$, whereas the SMS method is of order $2n + \rho n^{(d-1)/d}$
However, as Fig.~\ref{fig_comparemsame} shows, the errors in the SMS method are so small
that greatly compensate for the large computational cost. Also, as Fig.~\ref{fig_comparemv} shows,
the comparison is favourable with
the Galerkin method on Shishkin meshes, where systems of order $2^dn$ have to be solved.

On the other hand notice that
the change of variables $\tilde z_h=-z_h$, changes~$-A^T$ and~$-E^T$ in~$M$
to~$A^T$ and~$E^T$, respectively,
so that the coefficient matrix in the SMS method is symmetric. This allows
to use methods and software for sparse symmetric indefinite matrices which are generally faster than methods
for general sparse matrices (like those in the SUPG or Galerkin methods). A study of the performance
of the direct methods available today for sparse symmetric matrices can be found in~\cite{Gould-Scott-Hu}, where the codes MA57~\cite{MA57} and PARDISO, \cite{PARDISO1}, \cite{PARDISO2}, appear
as best performers. Note however that only serial versions were tried in~\cite{Gould-Scott-Hu} and that
that this is an area of fast development. As for iterative methods, the fact that SMS approximation can
be found by solving a symmetric system allows the use of three term recurrence methods like
the MINRES method of Paige and Saunders~\cite{PS75} (see e.\ g.\
\cite{Saddle-Acta}, \cite{Dollar}, \cite{Silvester-book}, \cite{PARDISO4}, \cite{PARDISO3} and the references cited therein for information on preconditioning this kind of systems).
We remark that for the systems in the SMS method, Matlab's backslash command seems to take no
advantage of the symmetry of the systems,  and, thus, no advantage has resulted from that symmetry
in the experiments reported in the present paper.

\end{example}
\medskip

\begin{example} \label{ej:curved} {\it Irregular grids on curved domains}. We now consider domains
where it may not be easy to set up a Shishkin grid. In particular, we consider the domain~$\Omega$
enclosed by the following curve: 
$$
\gamma(t)=\frac{26+7(1-\sin^9(2t))}{40(2+\sigma)\sqrt{2}}
\left[\begin{array}{rr}1& -1\\ 1& 1\end{array}\right]
\left[\begin{array}{c}2\cos(t)-\sigma\cos(2t)\\
2\sin(t)-\sigma\sin(2t)\end{array}\right],\qquad t\in[0,2\pi],
$$
where $\sigma=0.9$. The curve was created by tilting 45 degrees a centered trochoid and
altering its size in order to have the shape
depicted in~Fig.~\ref{fig_dominio}. In this domain we consider
problem~(\ref{eq:model}--\ref{eq:modelbc}) , with $b=[2,3]^T$ and constant forcing term
$f=1$, with Dirichlet homogeneous boundary conditions.
\begin{figure}[t]
$$
\includegraphics[height=3.5truecm]{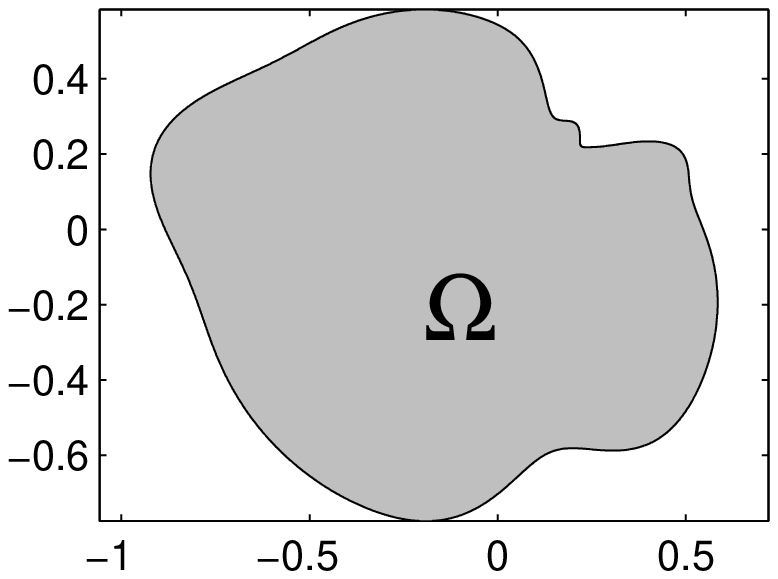}\quad
\includegraphics[height=3.5truecm]{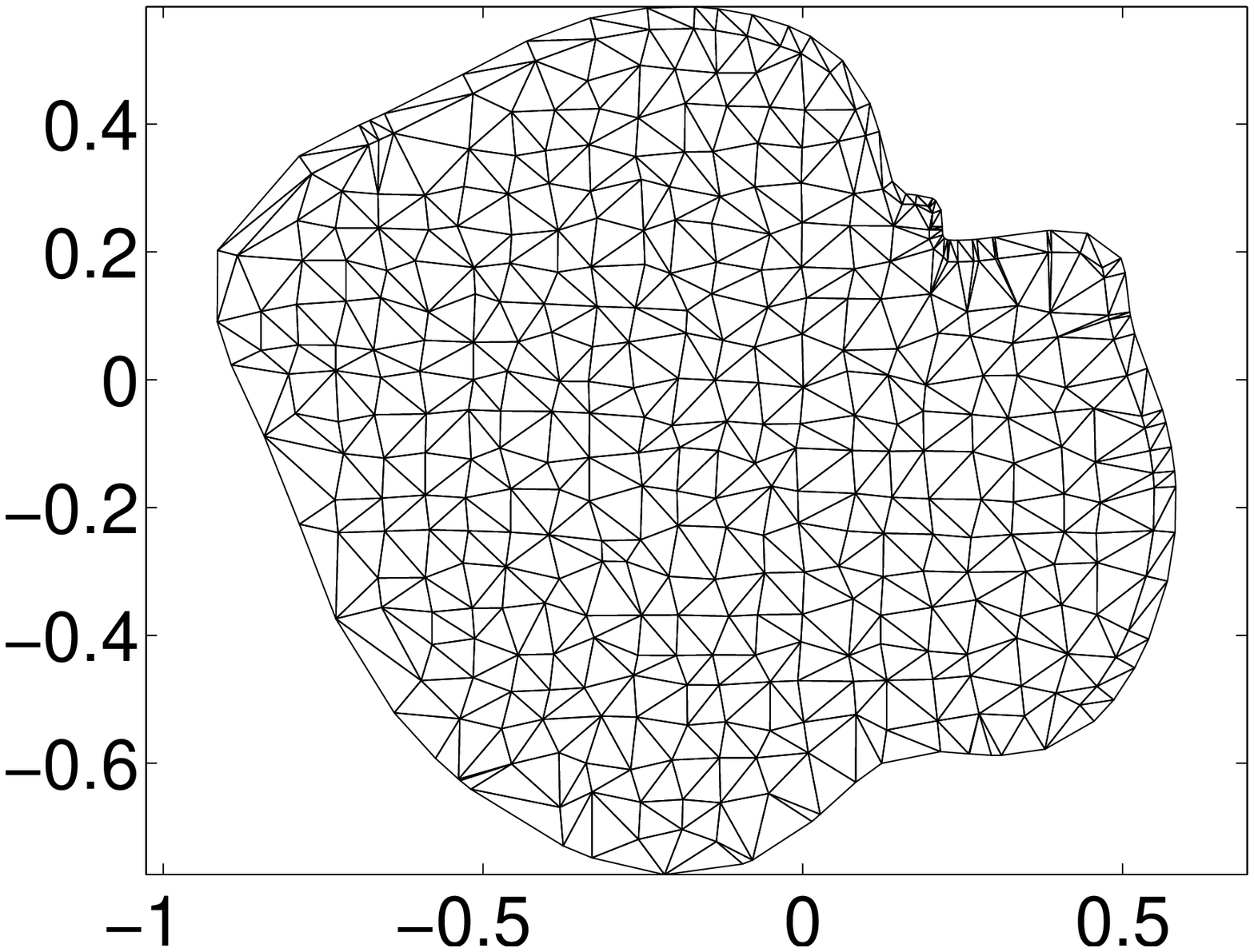}
$$
\caption{The domain of Example~\ref{ej:curved} (left) and a random grid (right).}
\label{fig_dominio}
\end{figure}

We study the behaviour of SMS methods on irregular grids, built as follows.
For a given positive integer $N$, 
we consider Delaunay triangulations with
$(N+1)^2$ points in $\overline\Omega$, where $4N$ of them are
randomly distributed on $\partial\Omega$. If $N_o$ of these are
on~$\Gamma_D^{0+}$, then, 
we built the strip of elements along~$\Gamma_D^{0+}$ as indicated in Fig.~\ref{fig_strip}. The remaining
points are generated by first fitting a uniform grid inside $\Omega$ and the outflow strip,
its diameter~$h$ in the $x$ and~$y$ direction being the value  for which the
number of points is the closest to~$(N-1)^2-N_o$, and then by displacing each point
randomly on the $x$ and~$y$ directions with a uniform distribution on~$[-h/3,h/3]$.
We show one of such grids for~$N=20$ in Fig.~\ref{fig_dominio}.

For $N=40$, we generated 200 random grids and on each of
them we computed the SUPG and SMS approximations and their $L^2$ error in the convective derivative
in~$\hat\Omega_h$, that is
$$
\left\| b\cdot \nabla w_h -1\right\|_{L^2(\hat\Omega_h)},
$$
$w_h$ being each of the three approximations. The errors on the 200 grids
for~$\varepsilon=10^{-4}$ and~$\varepsilon=10^{-8}$
are depicted in Fig.~\ref{fig_random40} (marked differently for each method).
\begin{figure}[h]
$$
\includegraphics[height=4.8truecm]{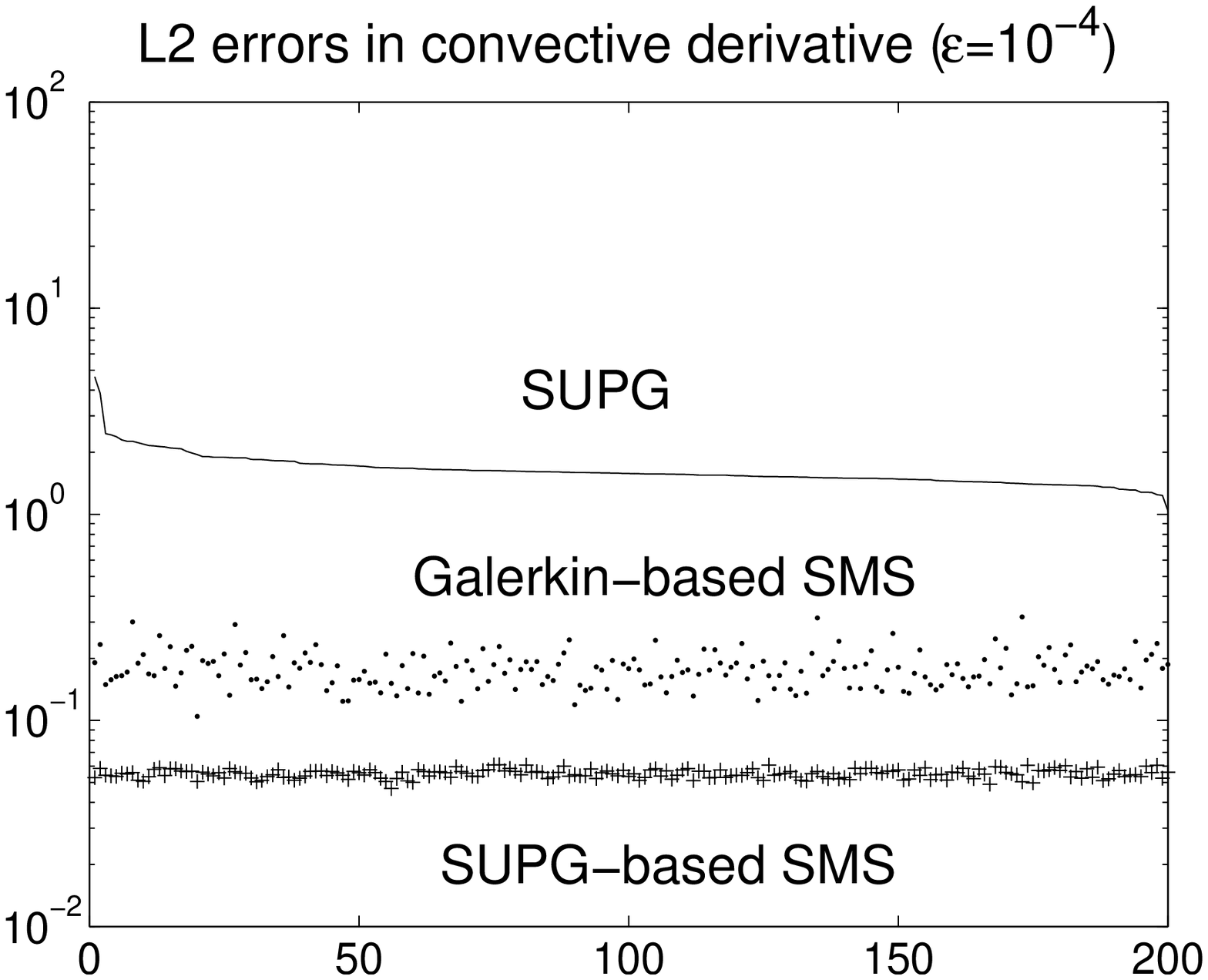}\quad
\includegraphics[height=4.8truecm]{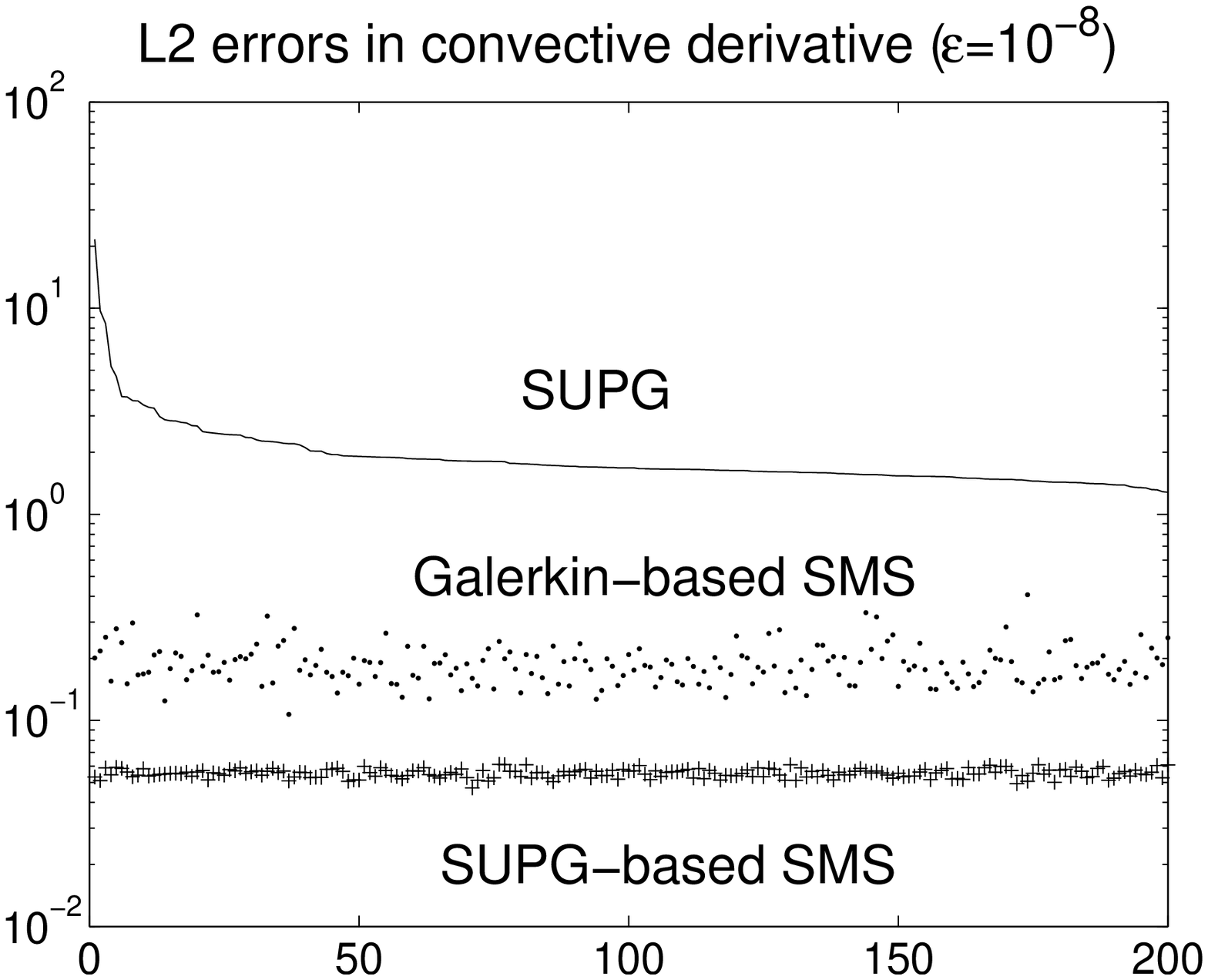}
$$
\caption{Errors in convective derivative on random grids for $\varepsilon=10^{-4}$ (left)
and~$10^{-8}$ (right).}
\label{fig_random40}
\end{figure}
The results have been reordered so that
those of the SUPG method appear in descending order. 
We see that the SMS methods clearly improve the errors of the SUPG method, specially so
in the case of the SUPG-based SMS method. Computing the ratios of the error of the SUPG method
and each of the SMS methods and taking the arithmetic mean, the resulting values for
the Galerkin-based and SUPG-based SMS methods are, respectively, 9.64 and 30.19
for~$\varepsilon=10^{-4}$ and~11.42 and~37.69 for $\varepsilon=10^{-8}$.
That is, the SMS methods commit errors
that are, on average, between 4 and 33 times smaller than those of the SUPG method.
For  $\varepsilon=10^{-8}$ we repeated these computations for growing values of~$N$, and the
rations between the errors of the SUPG method and the SMS method grew with~$N$. For
example, for~$N=320$, these were 46.93 and~380.41 for the Galerkin-based and SUPG-based
SMS methods, respectively.
 
In Fig.~\ref{fig_worstbest40}, where we compare the the SUPG-based SMS approximation that
produced the largest error (left) with that of the SUPG method that produced the smallest error
(right).
\begin{figure}[h]
$$
\includegraphics[height=4truecm]{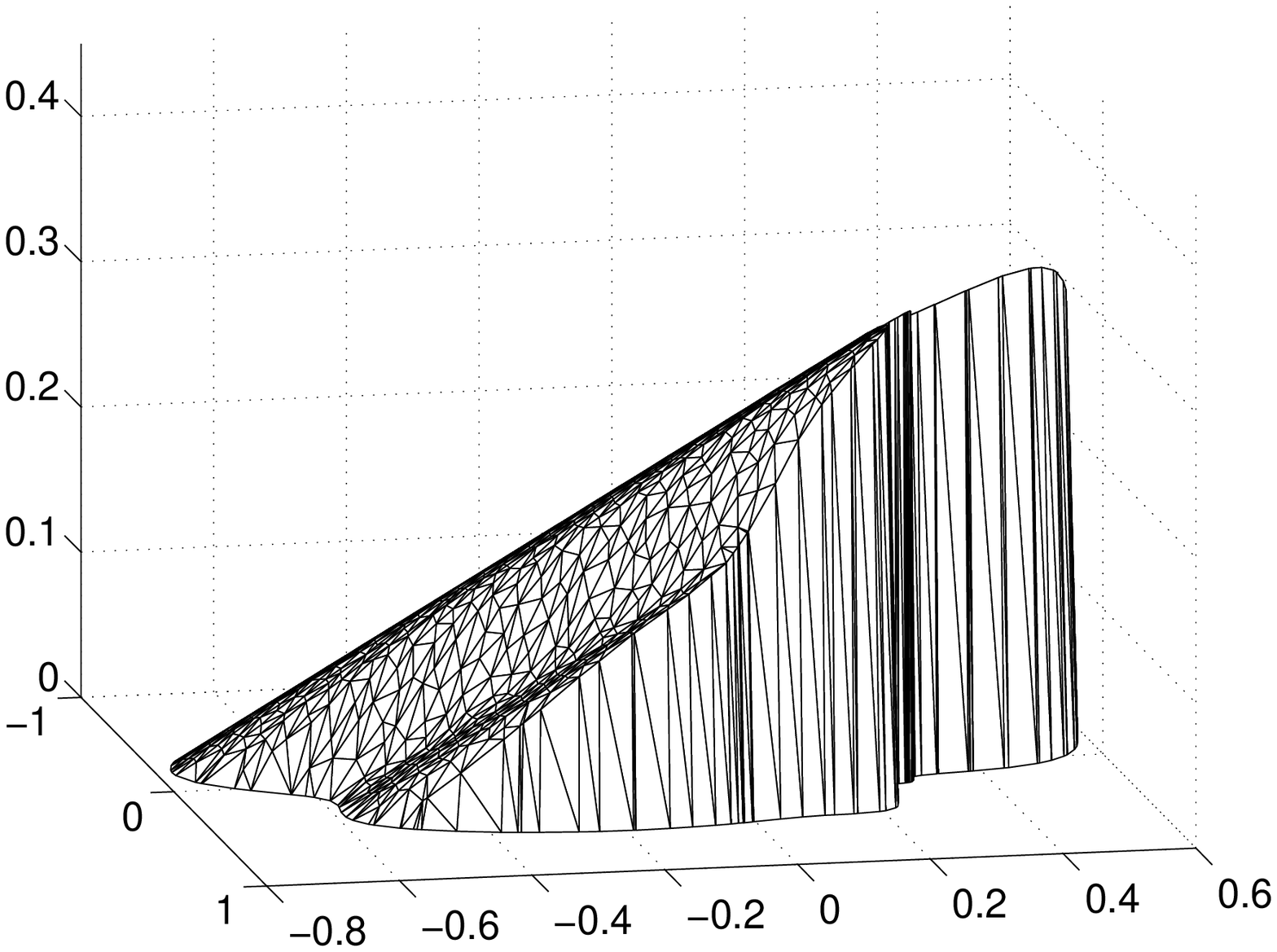}\quad
\includegraphics[height=4truecm]{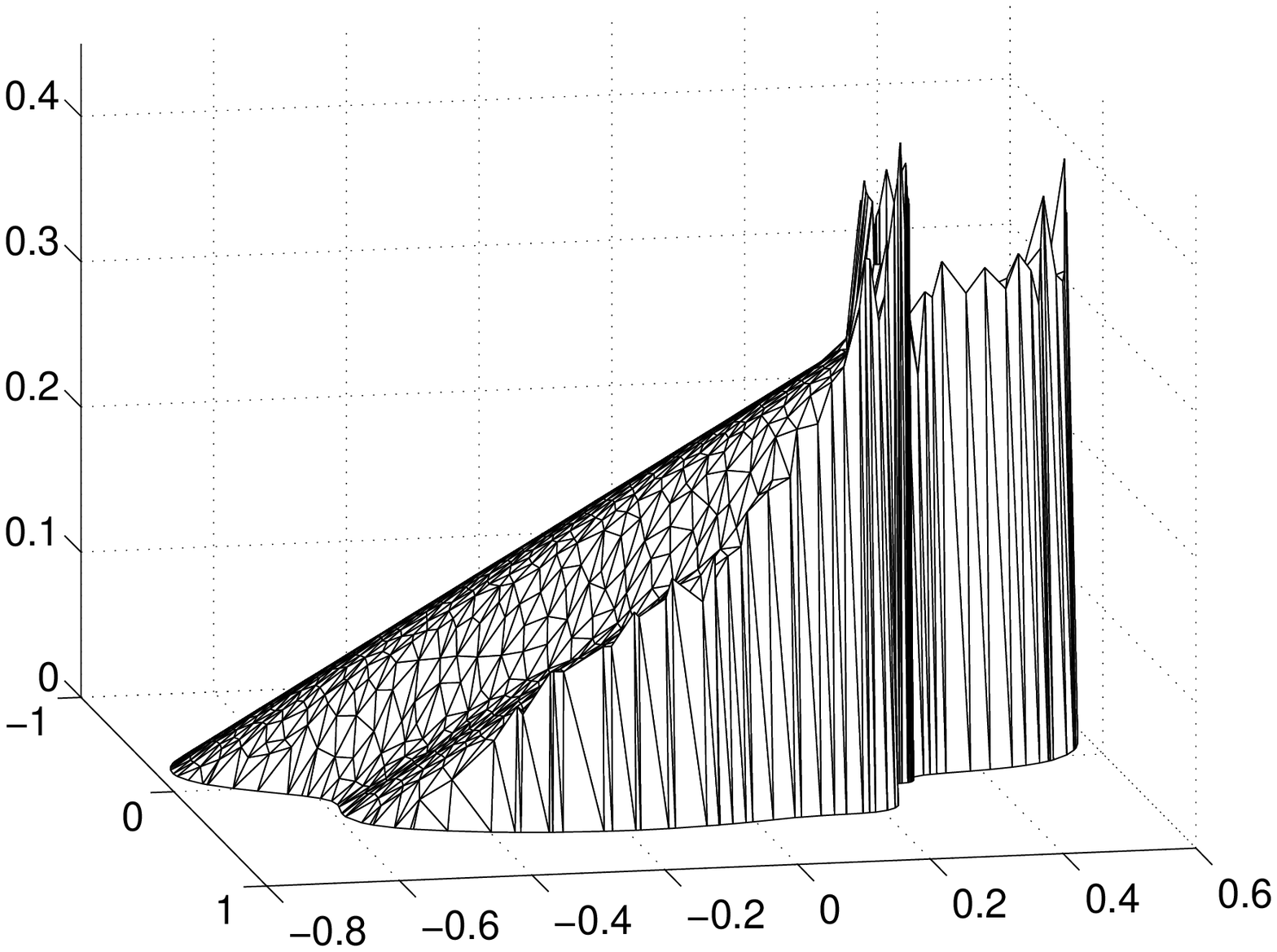}
$$
\caption{The worst case of the SUPG-based SMS method (left) and the best case of the
SUPG method (right).}
\label{fig_worstbest40}
\end{figure}
We notice the typical oscillations of the SUPG method, which are located only near the
outflow boundary but are of considerable amplitude. The
SMS method, on the contrary presents no oscillations. 
Similarly, in Fig.~\ref{fig_worstbest40gs} we show the Galerkin-based SMS approximations that
produced the largest and smallest errors (left and right, respectively).
Both of them compare very favourably with the best case of the SUPG method in
Fig.~\ref{fig_worstbest40}.
\begin{figure}[h]
$$
\includegraphics[height=4truecm]{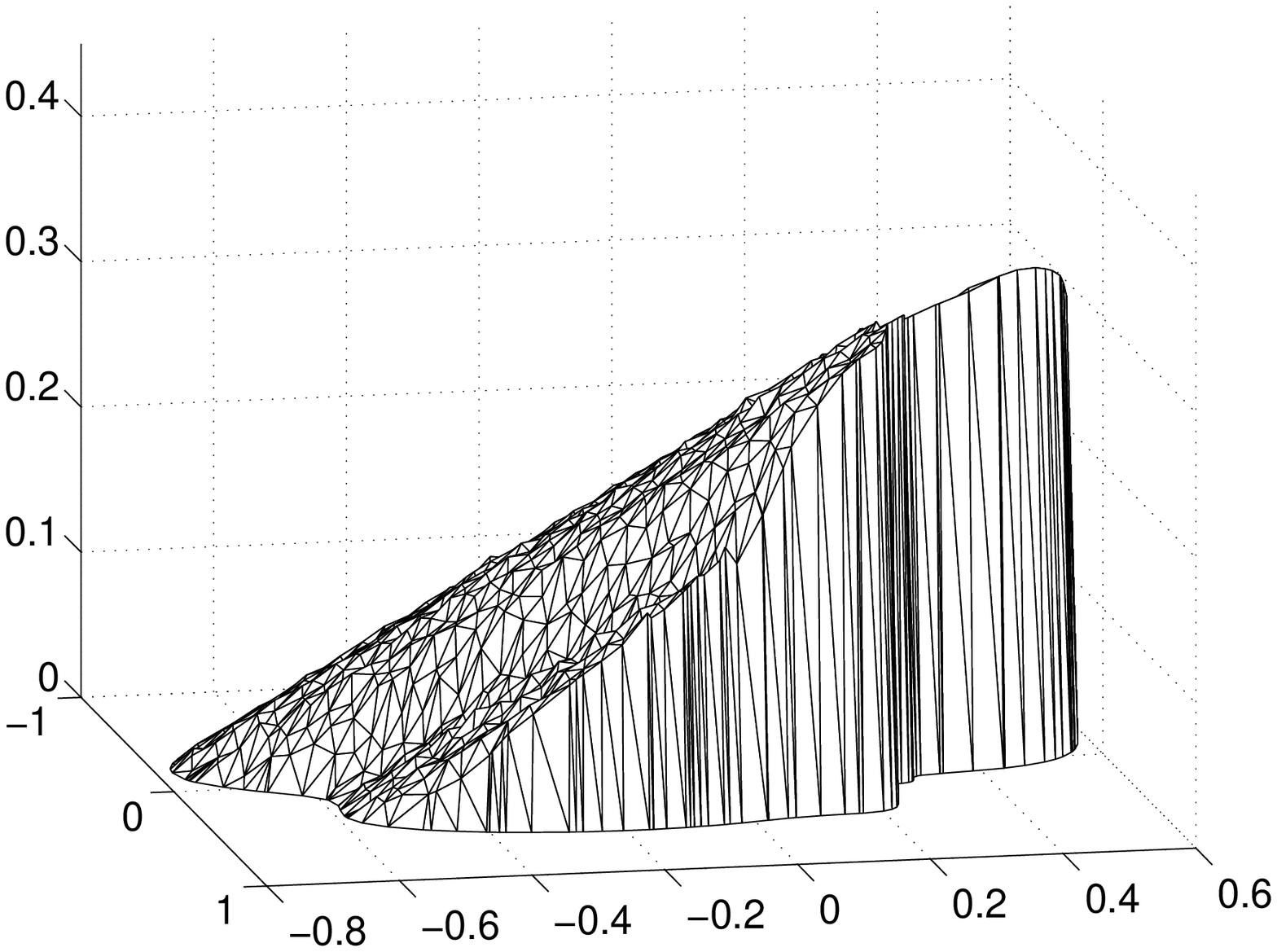}
\includegraphics[height=4truecm]{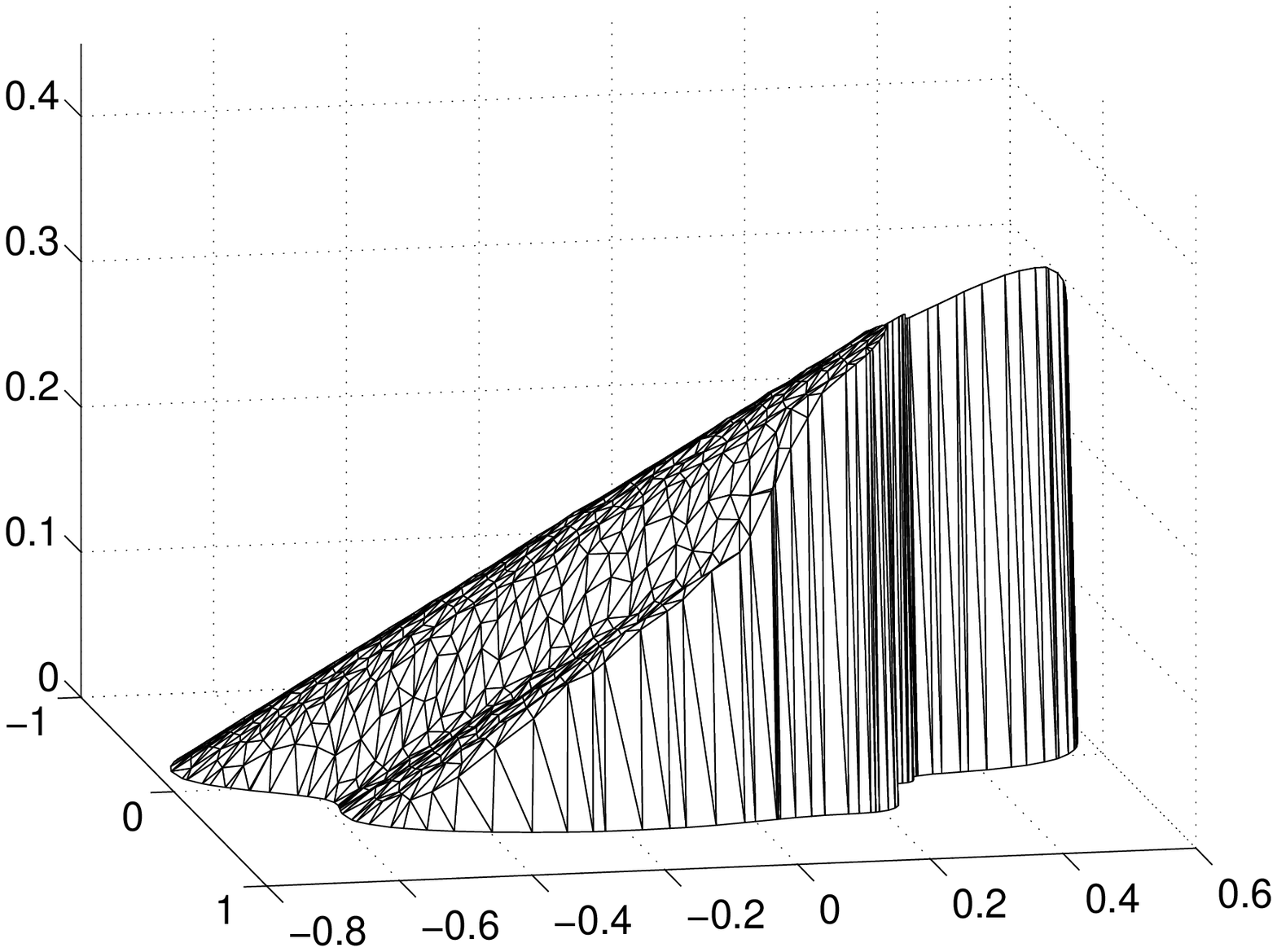}
$$
\caption{The worst and best cases of the Galerkin-based SMS method.}
\label{fig_worstbest40gs}
\end{figure}
We notice however that some small amplitude oscillations can
be observed in the worst case of the Galerkin-based SMS method. We believe that these
are due to the irregularity of the grid, since making the grids less irregular
(i.e., smaller random perturbations) diminishes
these small oscillations in the worst case, whereas making them more irregular increases
them. It is to be remarked, though, that the irregularity of the grid does not affect the
SUPG-based SMS method.
Some results explaining the
degradation of performance of the Galerkin method on irregular grids can be found in~\cite{Chen-Xu} and~\cite{Sun-Chen-Xu}.
Also, further numerical experiments with the SMS method on irregular grids
can be found in~\cite{sms2}.
\end{example}
\medskip

\begin{example} {\it Parabolic Layers}. \label{ej:parabolic}
We solve~(\ref{eq:model}) on $\Omega=(0,1)^2$, with $\varepsilon=10^{-8}$,
$b=[1,0]^T$ and~$f$ constant equal to~$1$ with Dirichlet homogeneous boundary conditions.
This is a well-known test case (see e.g.,~\cite{John-Knobloch-2007}). The solution presents an exponential layer at the outflow boundary at $x=1$ and parabolic
or characteristic layers along $y=0$ and~$y=1$ (see e. g.,~\cite{Roos-Stynes-Tobiska} for a precise
definition of these concepts).
In~Fig.~\ref{fig_parabolic} we show the SUPG and SMS approximations 
on a $20\times 20$ regular grid with Southwest-Northeast diagonals. We notice that whereas
the SUPG approximation suffers from the typical oscillations along the characteristic layers,
the SMS approximation is free of oscillations, both in the exponential layer at $x=1$ and
along the characteristic layers.
\begin{figure}[h]
$$
\includegraphics[height=4truecm]{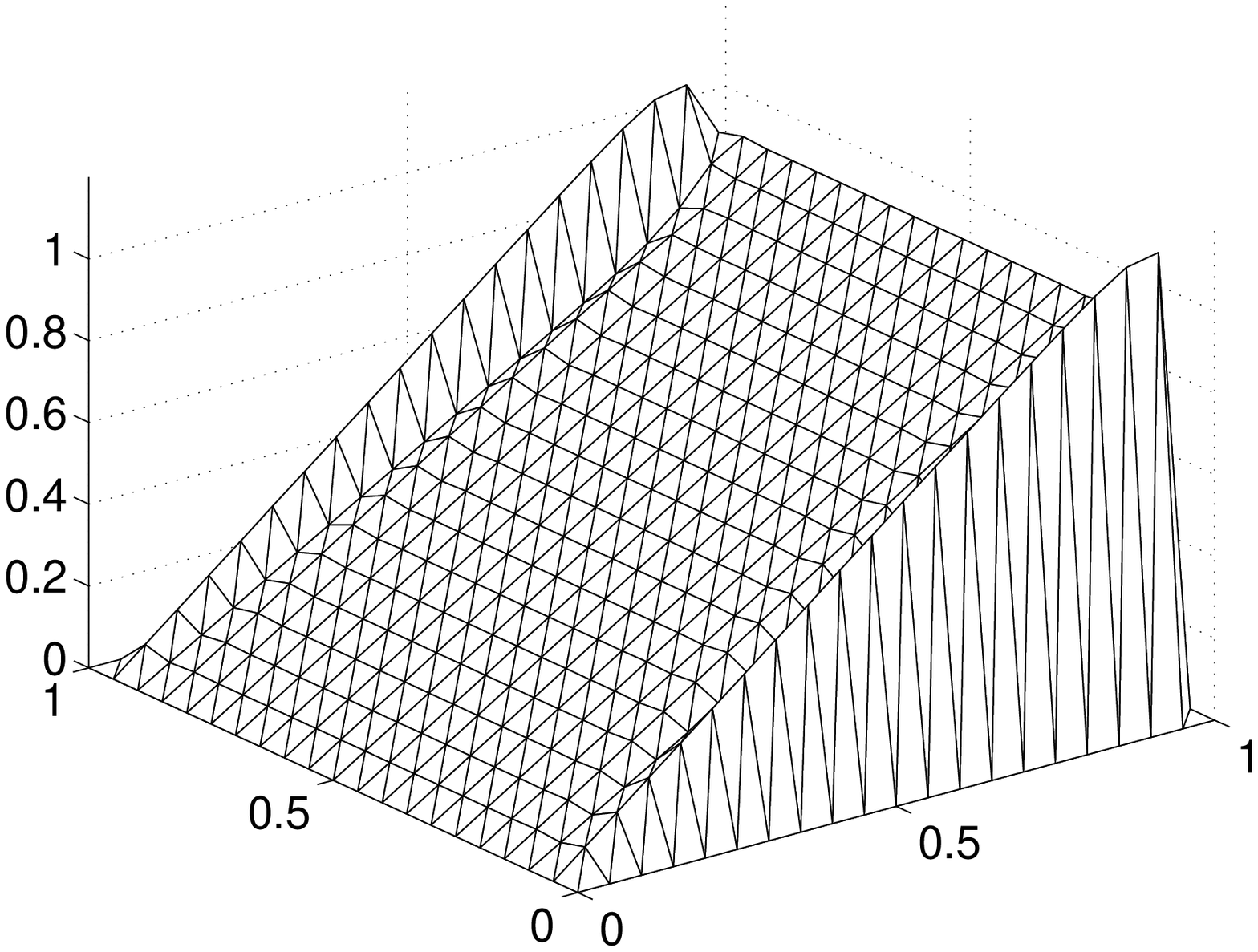}\quad
\includegraphics[height=4truecm]{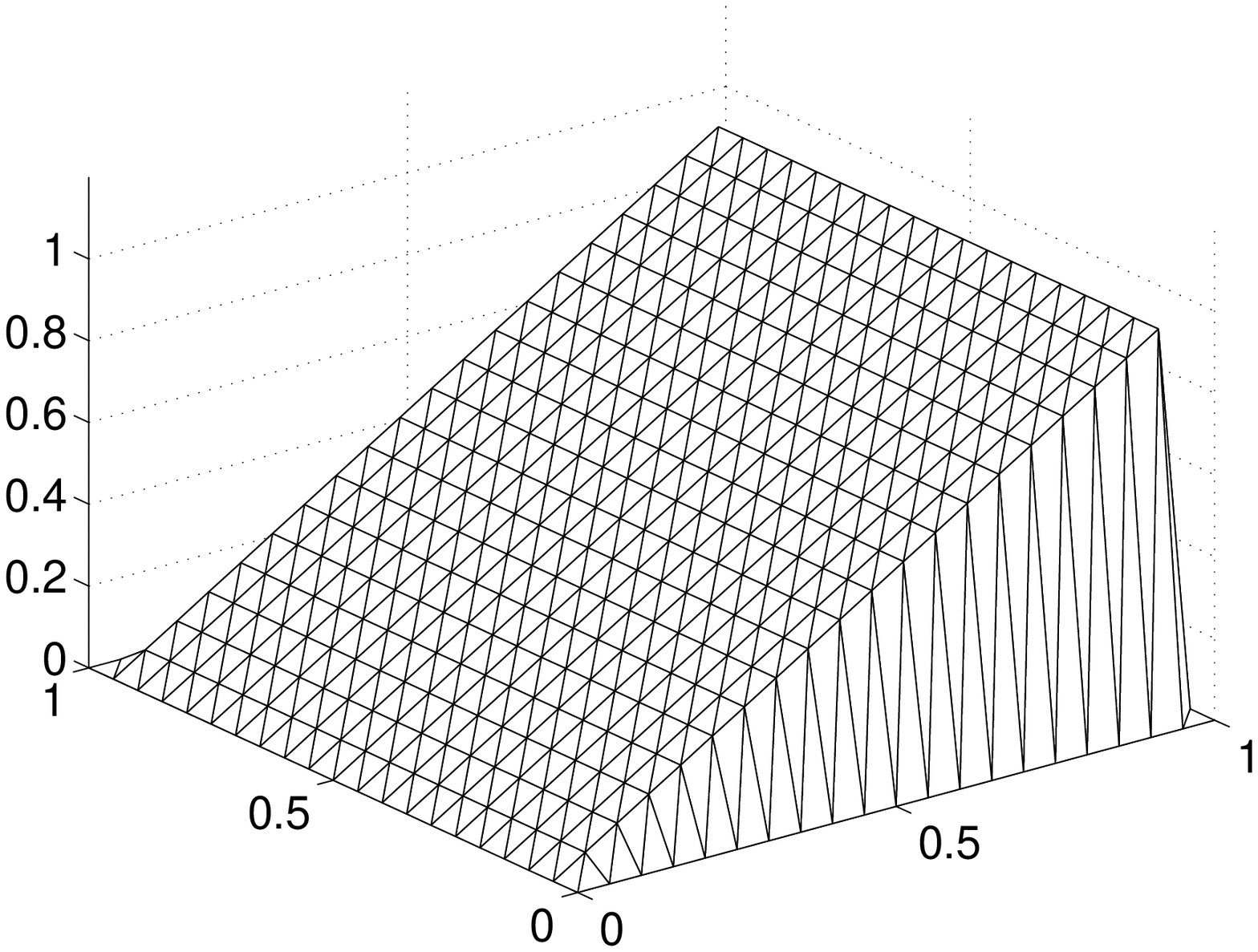}
$$
\caption{The SUPG (left) and the SMS approximation in Example~\ref{ej:parabolic}.}
\label{fig_parabolic}
\end{figure}


In~\cite{John-Knobloch-2007}, several techniques to reduce the oscillations of the SUPG method
are tested. In this example, for approximations~$w_h$ computed on a regular~$64\times 64$ grid,
the following quantities
\begin{align}
\hbox{\rm osc}&:= \max_{y\in\{1/64,\ldots,63/64\}} \{w_h(0.5,y)-w_h(0.5,0.5)
\},\label{osc}
\\
\hbox{\rm smear}&:= \max_{y\in\{1/64,\ldots,63/64\}} \{w_h(0.5,0.5)-w_h(0.5,y)
\},
\label{osc-smear}
\end{align}
are computed in~\cite{John-Knobloch-2007}
as a measure of the oscillations and the smearing along the characteristic
layers, desirable values being, respectively, between~0 and~$10^{-3}$ and
between~0 and~$10^{-4}$.
The value of~osc and~smear for all the methods tested in~\cite{John-Knobloch-2007} except one
was always larger than~$10^{-4}$. The values in the case of SMS methods were below~$10^{-14}$.
The value of osc in the SUPG method in our tests coincided with that in~\cite{John-Knobloch-2007},
0.134 (no value of smear was given in~~\cite{John-Knobloch-2007}).

Similar striking contrast between methods tested in~\cite{John-Knobloch-2007} and the SMS methods
can be found in the experiments on randomly-generated grids for this example in~\cite{sms2}.

\end{example}

\medskip
\begin{example} {\it Interior layers}. \label{ej:interior} This example is also taken
from~\cite{John-Knobloch-2007}. We solve~(\ref{eq:model}) on $\Omega=(0,1)^2$, with
$\varepsilon=10^{-8}$, $b=[\cos(-\pi/3),\sin(-\pi/3)]^T$, $f=0$, and $u=g$ on $\partial\Omega$
where
$$
g(x,y)=\left\{\begin{array}{lll} 0,&\quad &\hbox{\rm if $x=1$ or $y\le 0.7$}\\
1, && \hbox{\rm otherwise.}\end{array}\right.
$$
The solution possesses an interior layer starting at $x=0$ and $y=0.7$, and an
exponential layer on $x=1$ and on the right part of $y=0$.
For SMS methods, the only way we have conceived so far to deal with interior layers is to treat them
as parabolic layers. For this to be possible,
the grid has to include the characteristic curve that starts at the
discontinuity of the boundary data (or a polygonal approximation to it in case of curved
characteristics). We will refer to this characteristic curve as the {\it layer characteristic}. For grids as those depicted on the left of Fig.~\ref{fig_moved},
\begin{figure}[h]
$$
\includegraphics[height=3truecm]{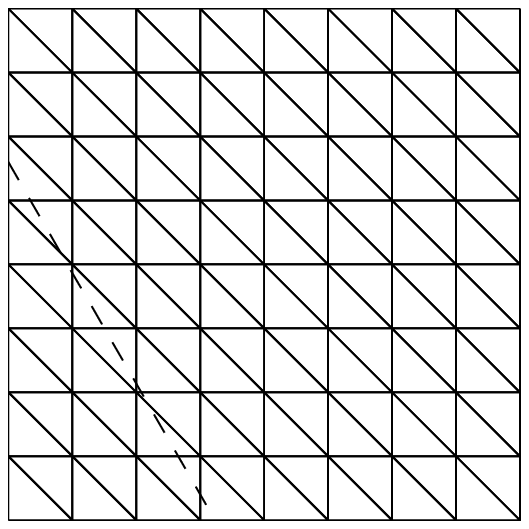}\quad
\includegraphics[height=3truecm]{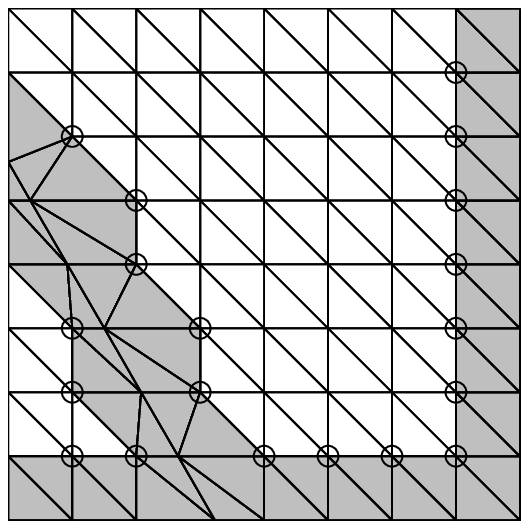}\quad
$$
\caption{Left, a $8\times 8$ uniform grid in Example~\ref{ej:interior}, with the characteristic
curve of the internal layer (dashed line). Right, the grid on the left but moving to the 
internal layer characteristic
its closest points on the grid.}
\label{fig_moved}
\end{figure}
where the layer characteristic is not part
of the grid, one possibility to make it part of it is
to move to the layer characteristic its closest points in the grid,
as we show on the right plot in Fig.~\ref{fig_moved} (see also the next example for an alternative).
Once this is done,  on has to set the value of the solution  along the layer characteristic
as part of Dirichlet
boundary conditions. To do this, we integrate
the reduced problem
(i.e., (\ref{eq:model}--\ref{eq:modelbc}) when $\varepsilon=0$)
along the layer characteristic.
In the present case, this procedure gives $u=0.5$ on the layer characteristic.

With these provisions, we compute the approximations of the
SUPG method SMS methods on a $16\times16$ grid. The results are shown 
in~Fig.~\ref{fig_internalplots} (the SMS approximations were identical and only the
Galerkin-based one is shown).\begin{figure}[h]
$$
\includegraphics[height=4truecm]{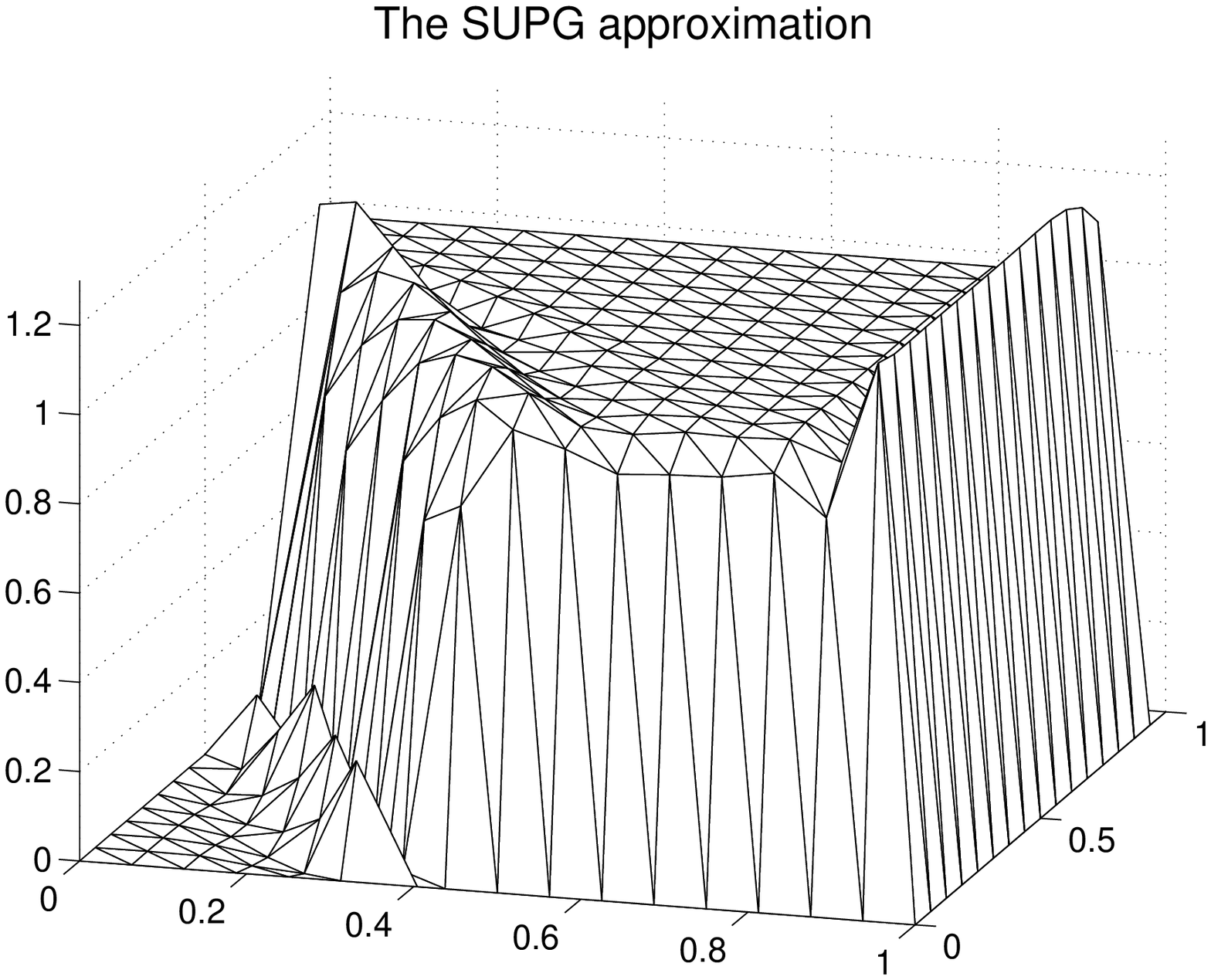}\quad
\includegraphics[height=4truecm]{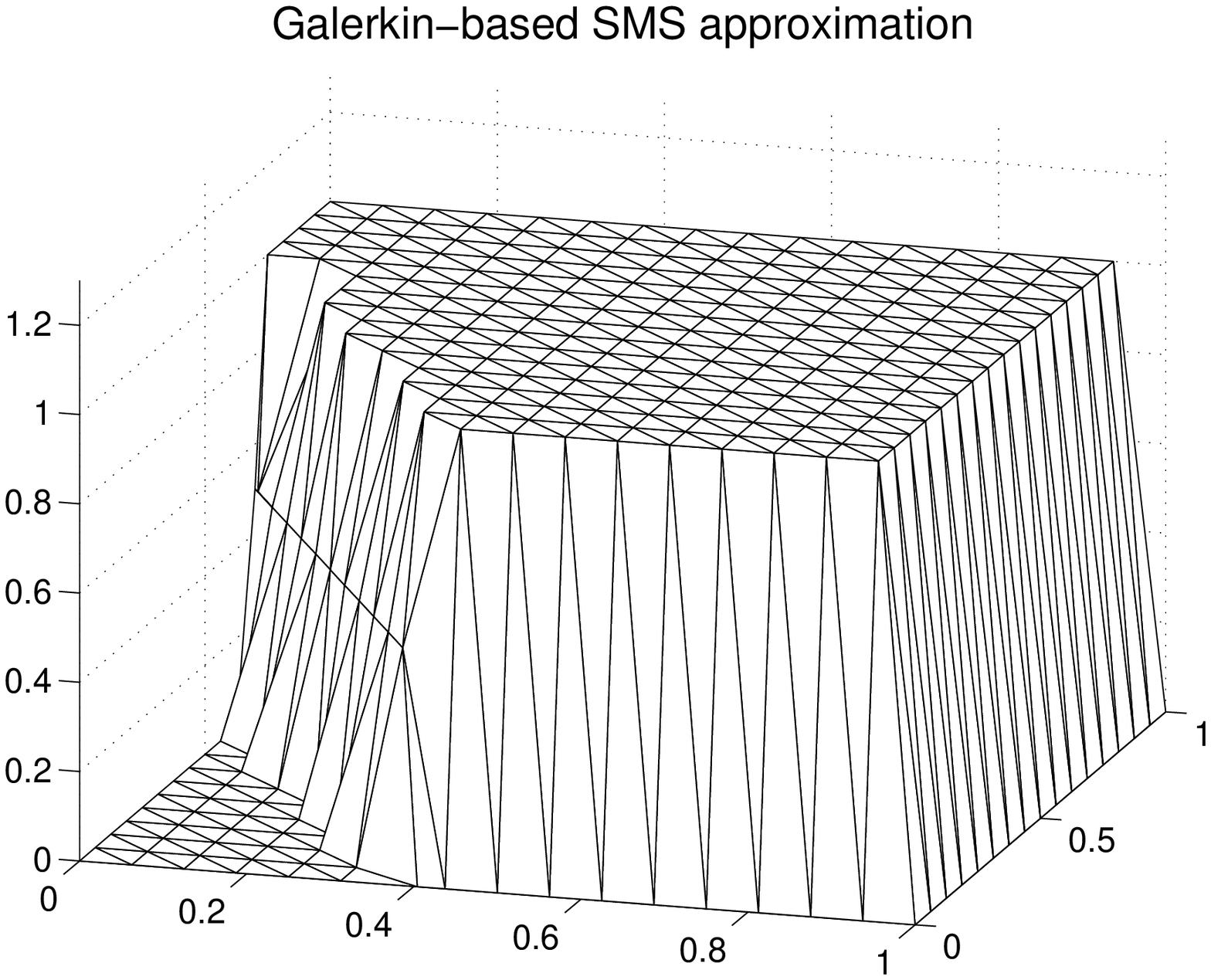}
$$
\caption{The SUPG and SMS (Galerkin-based) approximations on a~$20\times20$ regular grid
in Example~\ref{ej:interior}.}
\label{fig_internalplots}
\end{figure}
 As before, the SMS method produces an approximation with
no oscillations, in sharp contrast with the SUPG method. Comparison with the methods tested
in~\cite{John-Knobloch-2007} can be found in~\cite{sms2}, with results very similar to those
of the previous example.

\end{example}

\medskip
\begin{example}\label{ej:Hemker} {\it Hemker problem}. Here, $\Omega=(-3,9)\times(-3,3)\backslash
\{x^2+y^2\le 1\}$, $b=[1,0]^T$ and $f=0$. The boundary conditions are
\begin{equation}
u(x,y)=\left\{\begin{array}{lll} 0&\quad&\hbox{\rm if $x=-3$,}\\ 1,& &\hbox{\rm if $x^2+y^2=1$,}\\
\varepsilon\nabla u\cdot n=0,&& \hbox{\rm elsewhere.}\end{array}\right.
\label{Hemkerbc}
\end{equation}
This problem, which was originally proposed in~\cite{Hemker}, models a hot column (the unit circle) with normalized temperature equal to~1,
and the heat being transported in the direction of the wind velocity $b$. Thus, a boundary layer
appears in the upwind part of the unit circumference from the lowest to highest point,  and
two internal layers start from these
two points and spread in the direction of~$b$. Notice also that part of the boundary is curved,
a feature which is often encountered in applications.

We are going to present results corresponding
to~$\varepsilon=10^{-8}$ on the grid shown in Fig.~\ref{fig_hemkergrid}, which has 932 triangles and
531 nodes.
\begin{figure}[h]
$$
\includegraphics[height=4truecm]{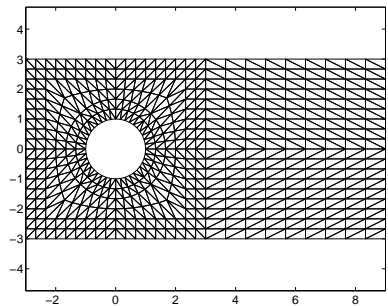}\quad
$$
\caption{The grid in Example~\ref{ej:Hemker}.}
\label{fig_hemkergrid}
\end{figure}
The interior layer characteristics are not part of this grid, as it can be seen in~Fig.~\ref{fig_corte},
where we show one of the layer characteristic in a dashed line. We have seen in the previous
example that interior layers are treated in the SMS method in the same way as characteristic
boundary layers, so that they must be part of the grid. Rather than moving grid points to
the layer as we did in the previous example, in the present one we will enlarge the grid
with more triangles and vertices so that the layer characteristic is part of it (an example
is shown in~Fig.~\ref{fig_corte}). This may be
useful in those cases in practice where, as commented at the end of Section~\ref{Se:elOmega_h} one does not have the freedom to choose or move the grid points.

Thus,
we now describe a general technique to enlarge grids in order to include an inner layer characteristic.
For simplicity we describe it for two-dimensional problems.
Let $\gamma$ be the inner layer characteristic.
We assume that the mesh is sufficiently fine so that $\gamma$
can be well approximated by a straight segment in the interior of every triangle~$\tau$ it intersects.
Let $\tau$ be such a triangle. If $\gamma$ passes through a vertex~$v$, then the $\tau$ is
bisected by~$\gamma$ into two triangles. These two triangles are included in the enlarged triangulation. Otherwise, $\gamma$ intersects two sides, and the element $\tau$ is divided by $\gamma$
into a triangle and a quadrilateral (see an example in~Fig.~\ref{fig_cortedetail}) which in turn is divided into four triangles. 
\begin{figure}[h]
$$
\includegraphics[height=1.9truecm]{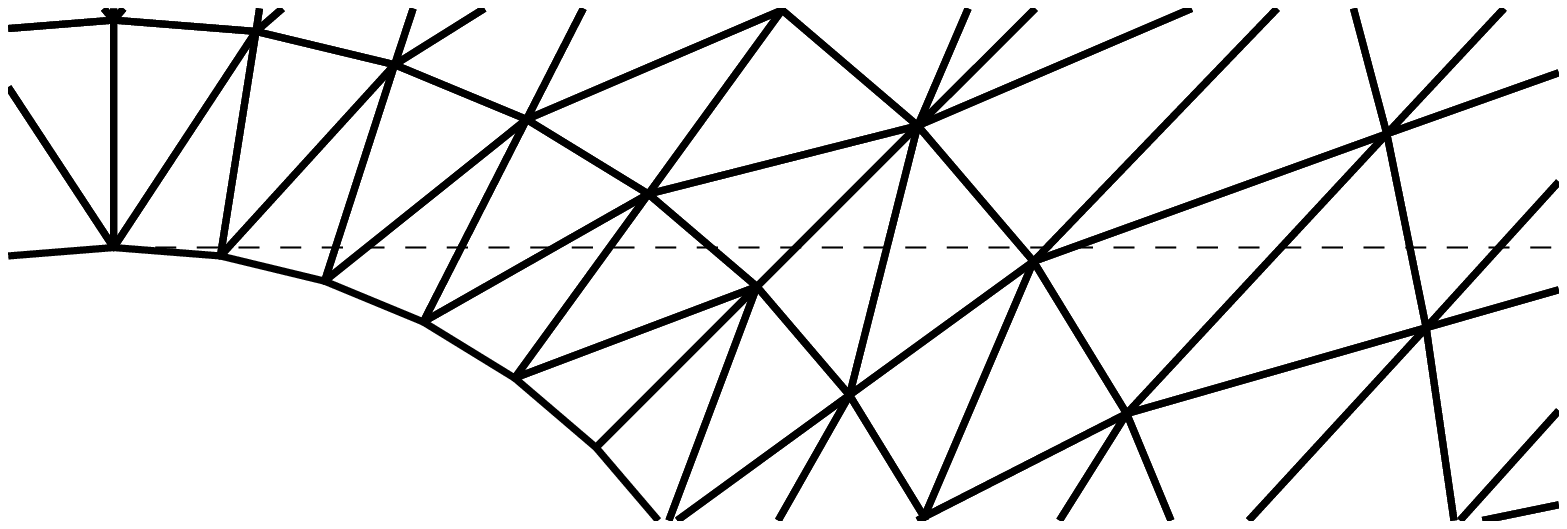}\quad
\includegraphics[height=1.9truecm]{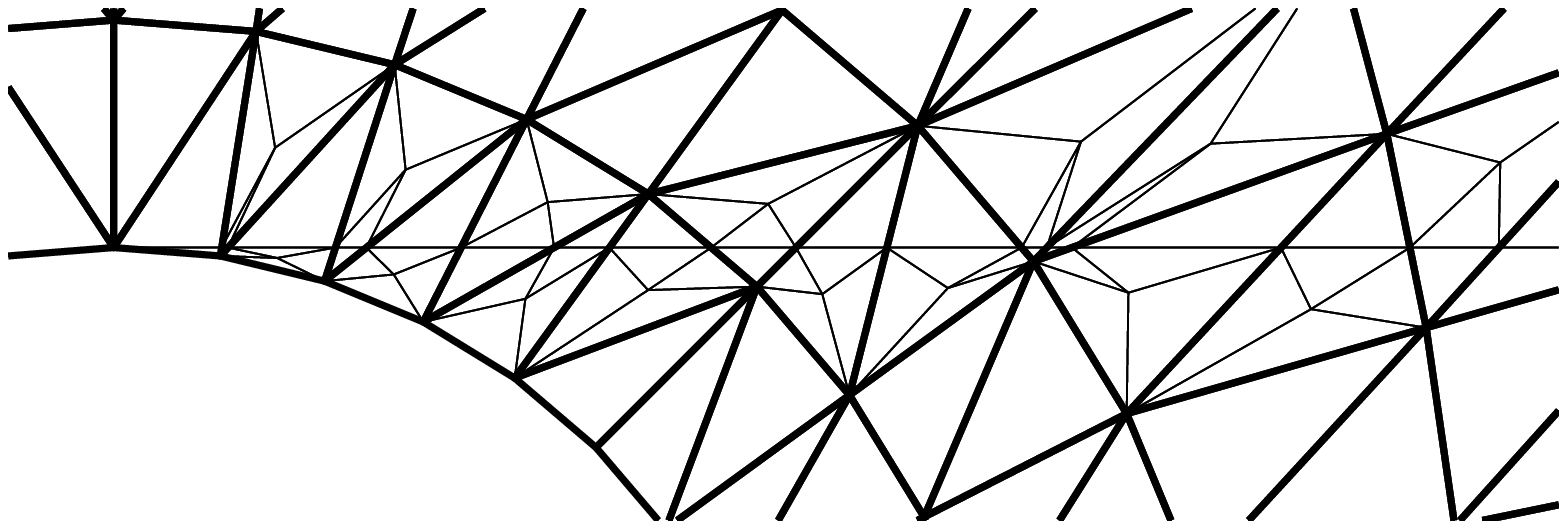}
$$
\caption{Left. Detail on the grid in Fig.~\ref{fig_hemkergrid} showing one of the
inner layer characteristic (dashed line). Right: the same grid enlarged with
more triangles and nodes to include the layer characteristic. New sides are plotted with
thinner lines.}
\label{fig_corte}
\end{figure}
\begin{figure}[h]
$$
\includegraphics[height=2.8truecm]{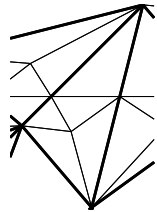}
$$
\caption{Detail of the enlarged grid in~Fig.~\ref{fig_corte} showing a triangle of the original
grid being divided by the layer characteristic in a triangle and a quadrilateral, which we
divide into triangles by joining the vertices with their arithmetic mean.}
\label{fig_cortedetail}
\end{figure}
We remark that it is easy to conceive better strategies to enlarge the triangulation in
order to include the layer characteristic
(techniques that avoid long-shaped triangles, for example), but we have
chosen this one due to its simplicity. Nevertheless, as the experiments below show, its simplicity
does not prevent it from obtaining excellent results.

Fig.~\ref{fig_hemkersols} shows the SUPG and the SMS approximations (Galerkin-based).
Oscillations can be clearly seen in the SUPG approximation, but they are absent in the
SMS approximation. This can be better observed in the bottom plots, where a different
point of view is taken.
\begin{figure}[h]
$$
\displaylines{
\hskip -0.7truecm
\includegraphics[height=4truecm]{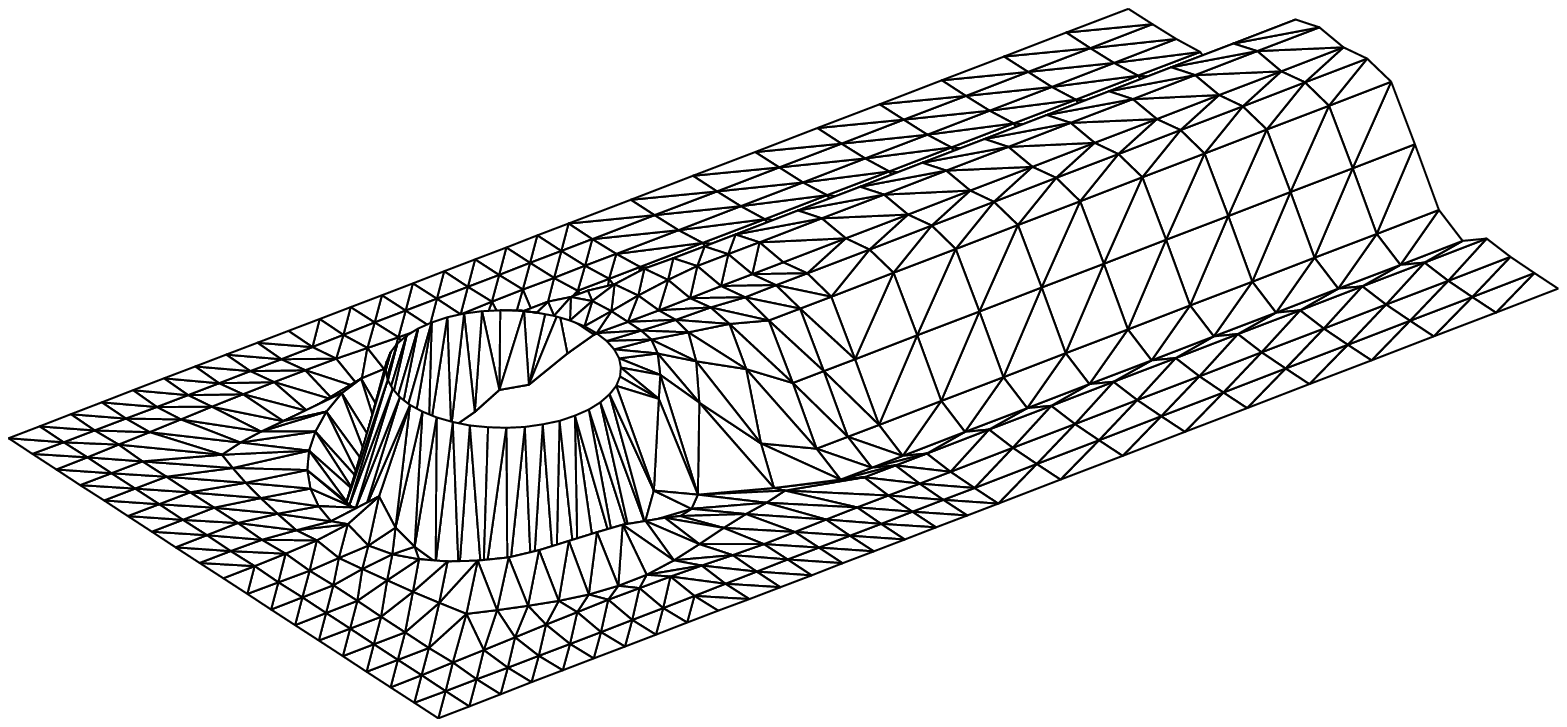}\hskip-0.1truecm
\includegraphics[height=4truecm]{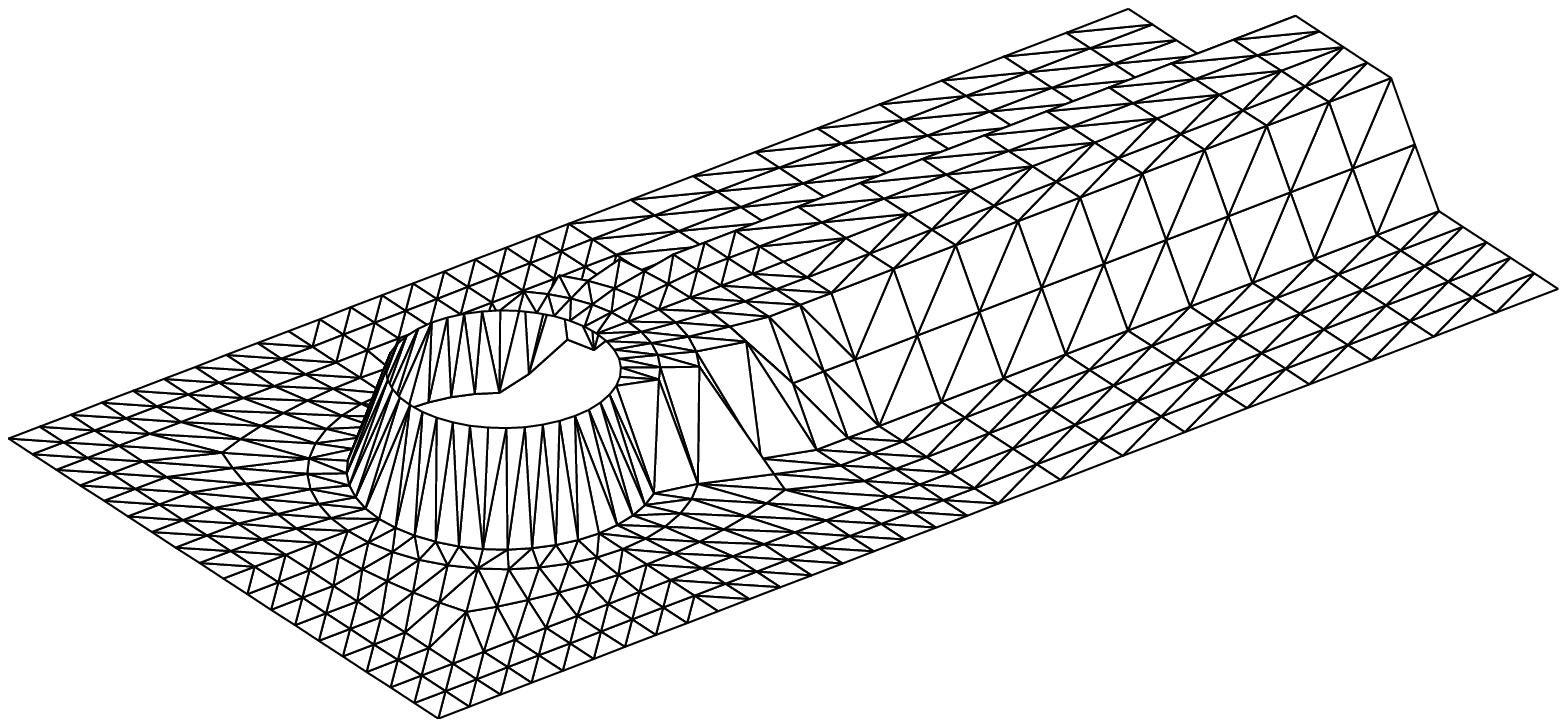}
}
$$
\vskip -1.15truecm

$$
\displaylines{
\hskip -0.7truecm
\includegraphics[height=3truecm]{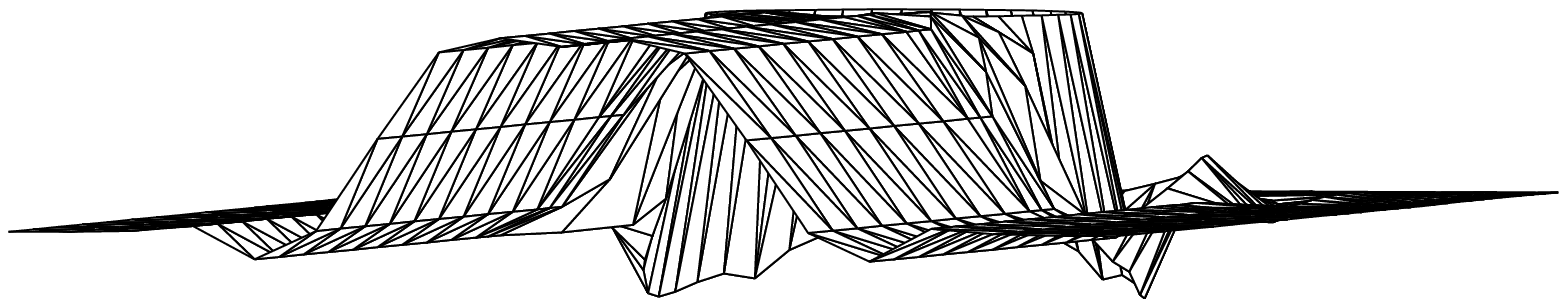}
\includegraphics[height=3truecm]{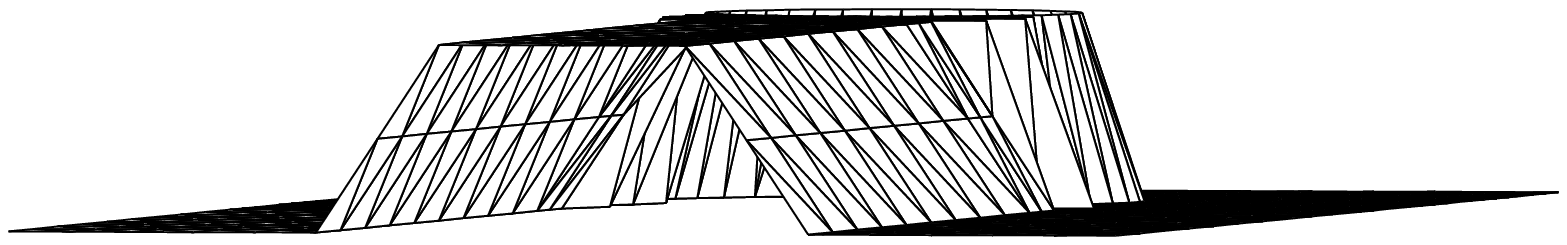}
}
$$

\caption{The SUPG approximation (left) and the Galerkin-based SMS approximation (right)
on the grid depicted in~Fig.\ref{fig_hemkergrid}. A different point of view is taken on the
bottom plots.}
\label{fig_hemkersols}
\end{figure}
Following~\cite{Augustin-et-al-2011}, we measure the undershoots of an approximation~$w_h$
as $\min\left\{w_h\right\}$ and the overshoots as~$\max\left\{w_h-1\right\}$.
The over and undershoots for the SUPG method were $0.04$ and $-0.52$, whereas for
both SMS methods the overshoots where smaller than~$10^{-15}$, and the undershoots
were~0. Similar results were obtained by doubling and multiplying by four the number
of subdivisions in each coordinate direction. The results of the SMS method compare very
favourably not only with the SUPG method, but with most of the methods tested
in~\cite{Augustin-et-al-2011}, which had values very similar to those of the SUPG method.
They also compare very favourably with results in~\cite{TobiskaLP}, where LPS methods present
oscillations of about~$5\%$ of the jump, whereas the SMS methods of less than $10^{-13}\%$.

In \cite{sms2}, the previous experiment is repeated with the vector field~$b$ changed to
$b=[\cos(\theta), \sin(\theta)]^T$, for 100 equidistant values of~$\theta$ between $(0,\pi/4]$,
with results similar to those shown here for all cases except four, where, in order to
get over and undershoots of order $10^{-12}$ it was also necessary for
nodes $O(h_{\min}^2)$ away from interior
layers to be allowed to move (see \cite{sms2} for details).
\end{example}

\medskip
\begin{example}\label{ej:Silvester} {\it Double-glazing test problem\/} \cite{Silvester-book}. This is an example
where the vector field $b$ has vortices. The domain is $\Omega=(-1,1)^2$, $f=0$, and the wind velocity
is $[y(1-x^2),-x(1-y^2)]$, so that the characteristic curves are the closed curves given by
$$
\{(x,y)\mid (1-x^2)(1-y^2)=\hbox{\rm constant}\},
$$
(See Fig.~\ref{fig_silv0}). The Dirichlet boundary conditions are $u=1$ on $x=1$ and $u=0$ otherwise, so that there
are discontinuities in the data at corners in $x=1$, $y=\pm1$. These discontinuities give raise to parabolic
boundary layers which, according to~\cite{Silvester-book}, have a structure difficult to compute
by asymptotic techniques.

\begin{figure}[h]
$$
\includegraphics[height=2.8truecm]{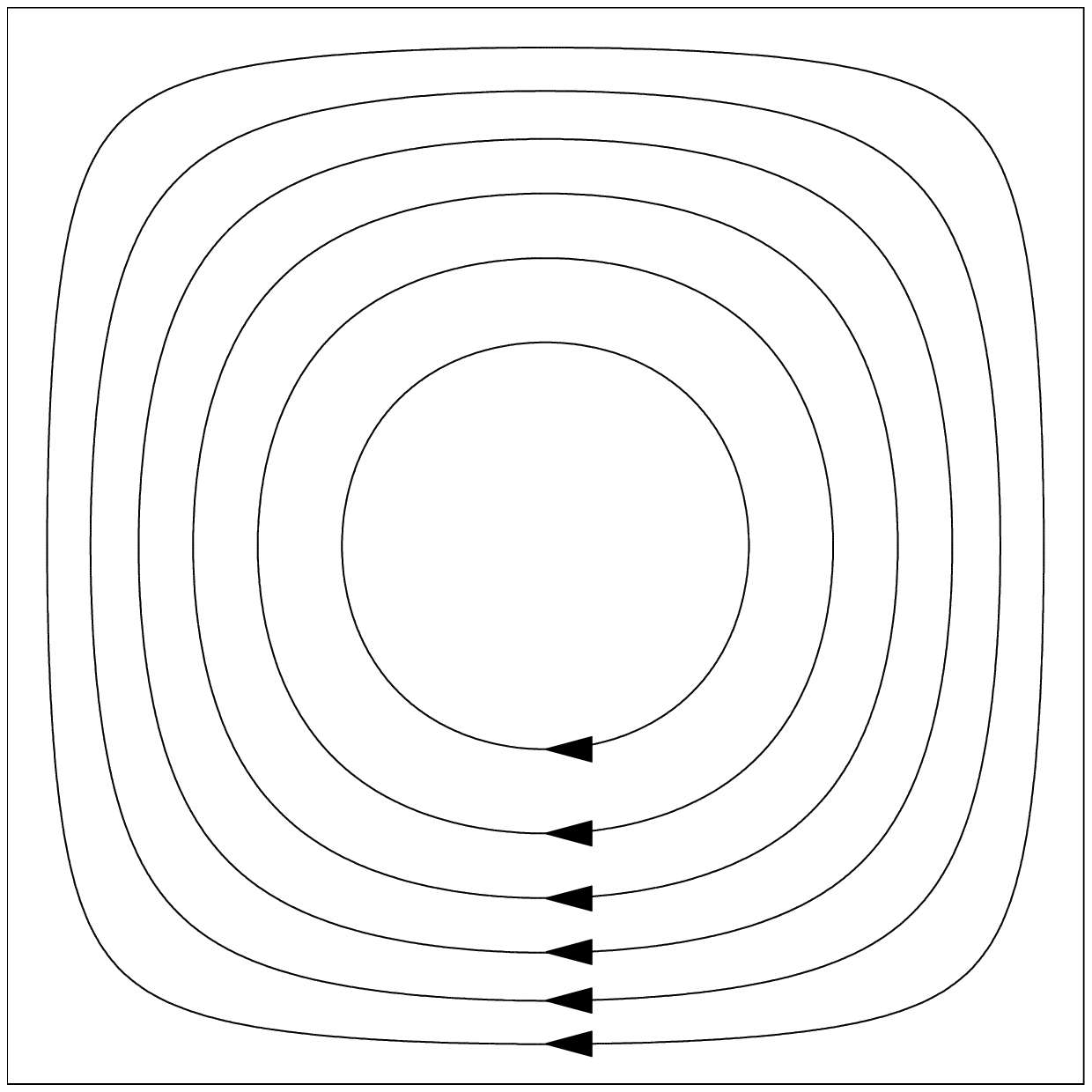}\qquad
\includegraphics[height=3.4truecm]{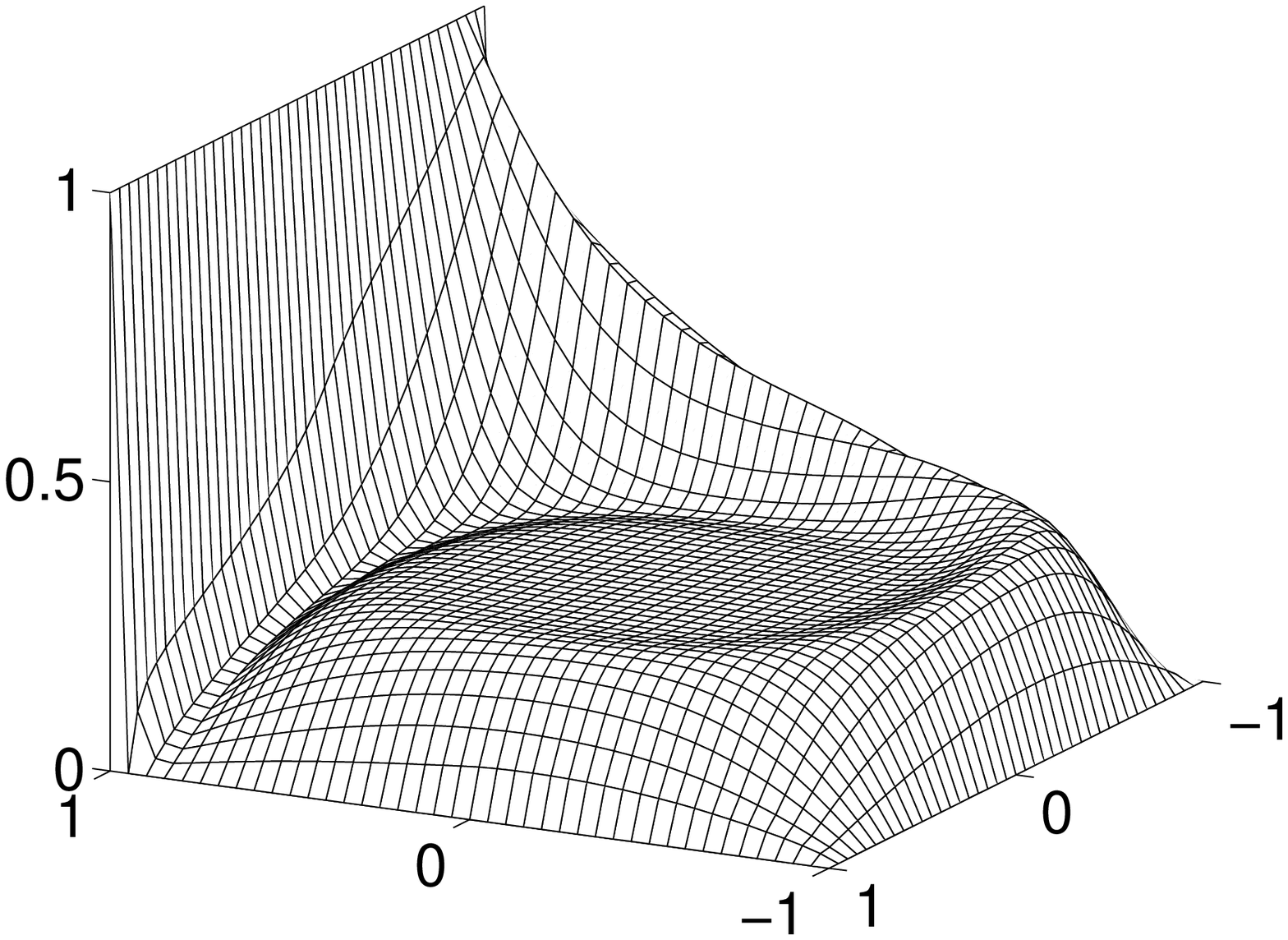}\qquad
\includegraphics[height=2.8truecm]{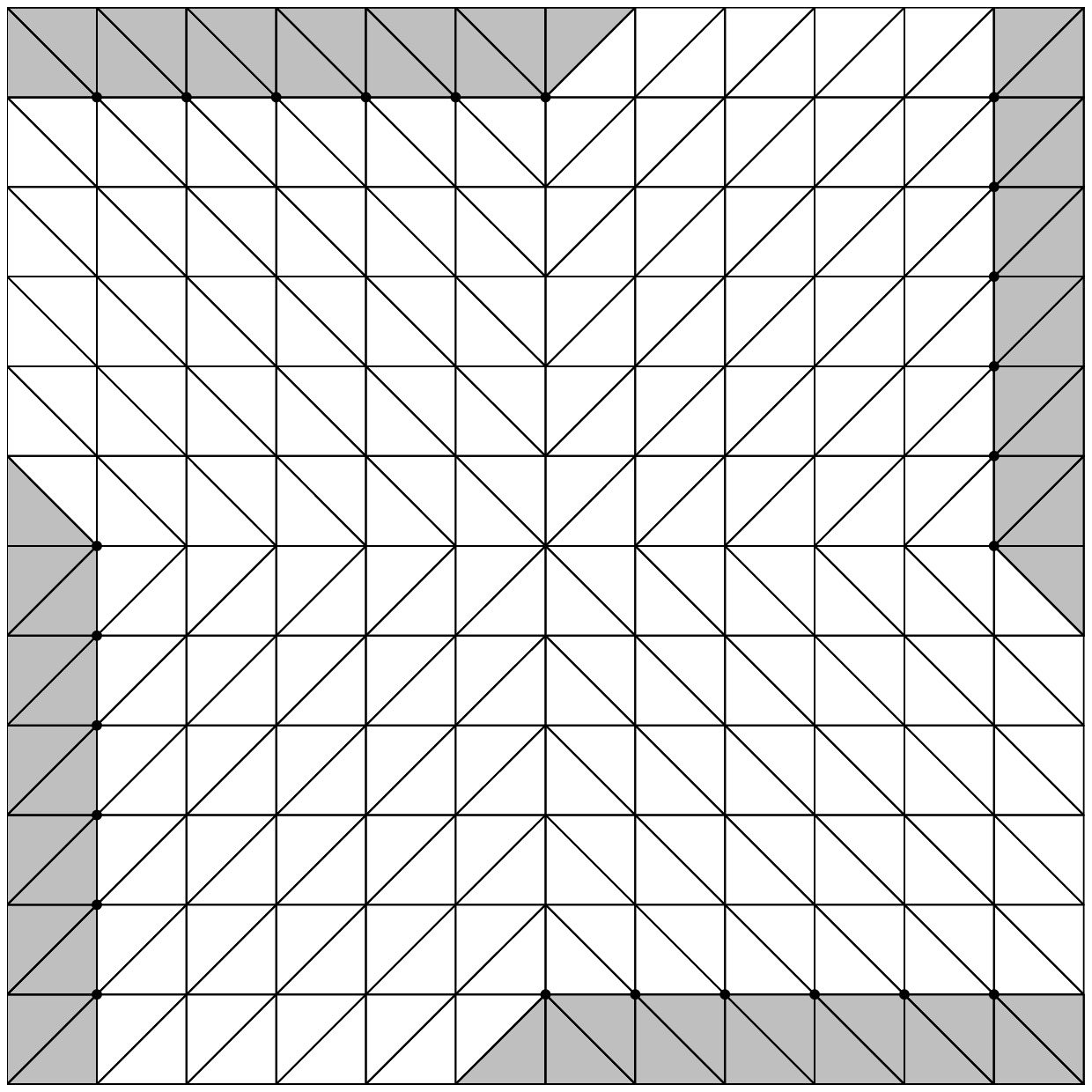}
$$
\caption{Left: The streamlines in Example~\ref{ej:Silvester}.
Center: the solution in Example~\ref{ej:Silvester} for~$\varepsilon=0.005$.
Right: A regular mesh, with the set $\Omega_h^{+}$
shadowed in grey  and points of~${\cal N}_\delta$ are marked with
circles.}
\label{fig_silv0}
\end{figure}

Notice that the hypotheses stated at the beginning of Section~\ref{Se:multidimensional}
($b(x)\ne 0$ and all characteristics entering and leaving the domain in finite time)
do not apply to vector fields with vortices
or closed integral curves.  Nevertheless, as we now show, the results of the SMS method
in the present example are as good as in the previous ones, even though
part of the analysis in Section~\ref{Se:analysis_2d} does not apply to the present case.

Observe that since the four sides of the boundary are themselves characteristic
curves, we have $\Gamma_D^{0+}=\partial
\Omega$, so that building $\Omega_h^{+}$ as described
in~Section~\ref{Se:multidimensional} results in~$\Omega_h^{+}$  being
composed by all elements touching the boundary. In this case
the SMS method produces an approximation equal to~$0$ on all nodes except on those on $x=1$, where it takes value one, and this is only
correct for~$\varepsilon=0$.  For $\varepsilon>0$, 
better results are obtained with the SMS method if,
as we show in~Fig.~\ref{fig_silv0}, the
set $\Omega_h^{+}$ is shrunk so as to correspond to that of a slightly smaller rectangular
domain, that is,
elements are included in $\Omega_h^{+}$ if they
intersect the outflow boundary of $(-1+\delta,1-\delta)^2$ for a small $\delta>0$
(other possibilities to obtain better
results with vector fields with vortices will be reported in future works).
With this selection of~$\Omega_h^{+}$,
the results of the Galerkin-based SMS method for $\varepsilon=10^{-4}$ can be seen
in~Fig.~\ref{fig_silv1}.
\begin{figure}[h]
$$
\includegraphics[height=4truecm]{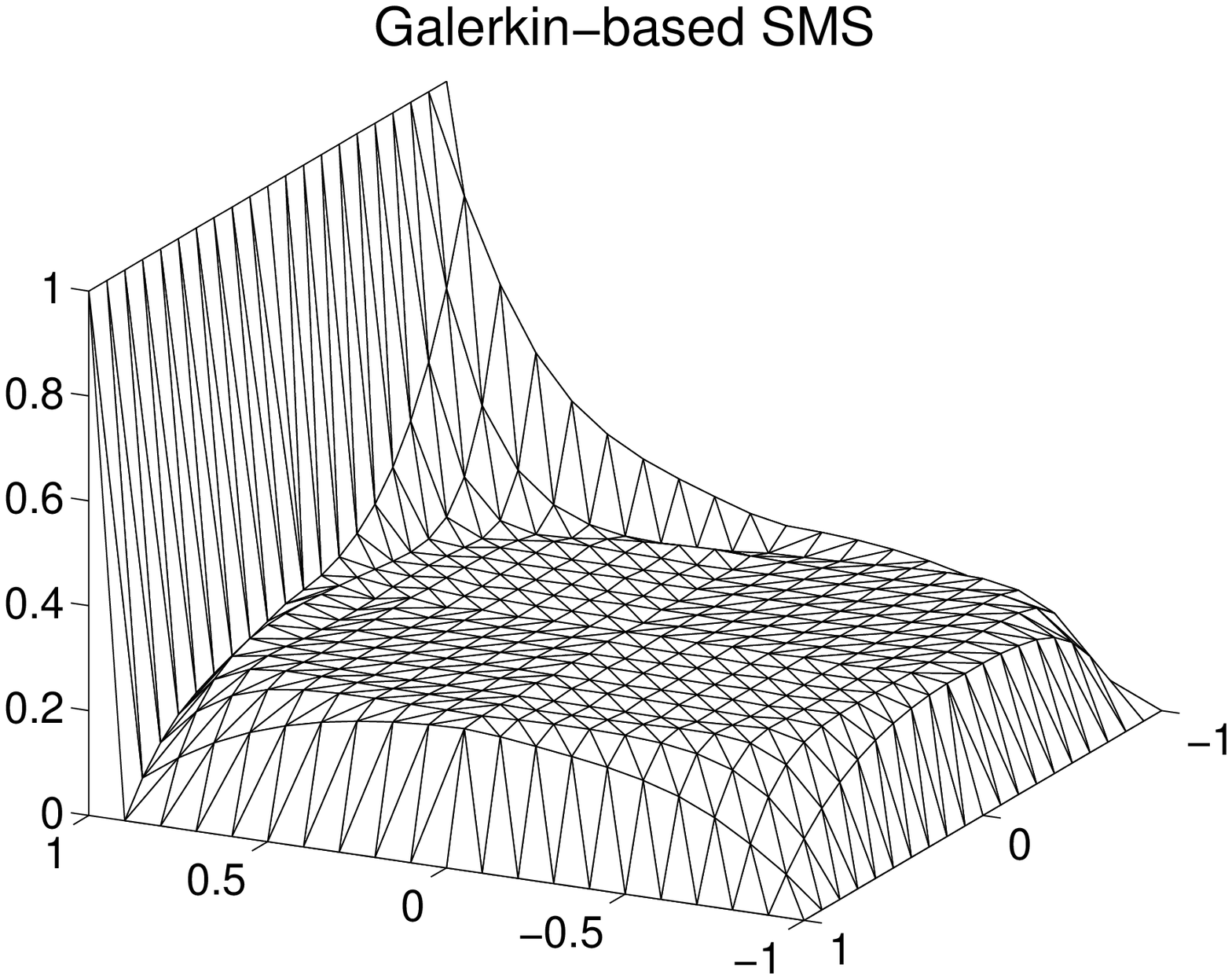}\quad
\includegraphics[height=4truecm]{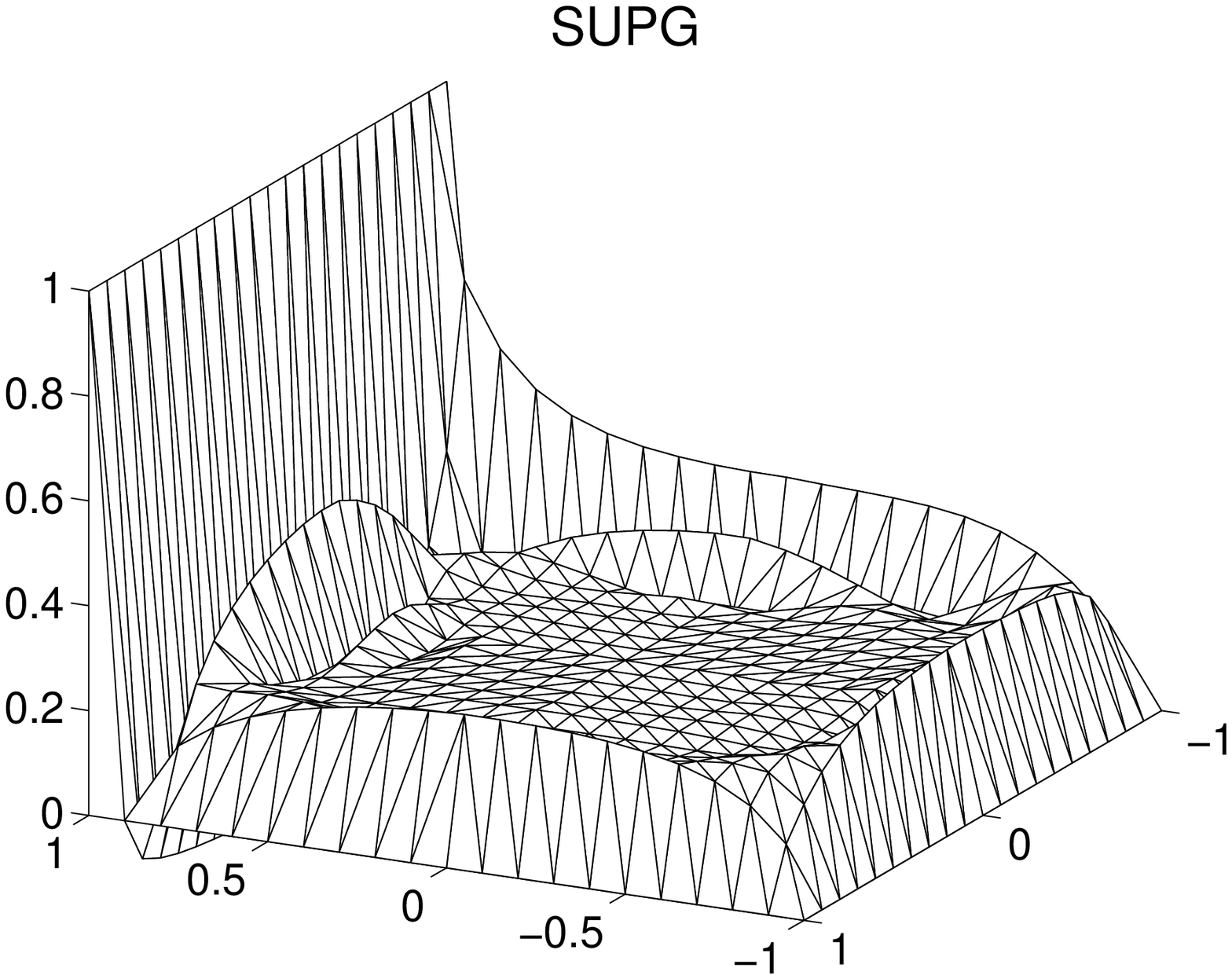}
$$
\caption{The Galerkin based SMS method (left) and the SUPG method on a $20\times 20$ regular
grid in Example~\ref{ej:Silvester} for $\varepsilon=10^{-4}$.}
\label{fig_silv1}
\end{figure}
Also shown in Fig.~\ref{fig_silv1} is the SUPG approximation, where we can observe oscillations
along the characteristic layers, especially at $x=1$ and $y=-1$. In sharp contrast, the SMS method
produces a nonnegative approximation. For this value of~$\varepsilon$, in order to obtain nonnegative approximations on a regular $N\times N$ mesh with the SUPG method one must take $N\ge  144$,
and $N\ge 240$ with the Galerkin method;
for $\varepsilon=10^{-5}$, this is achieved for $N\ge 682$ and $N\ge 800$ with the SUPG and Galerkin methods
respectively, and for $\varepsilon=10^{-6}$, neither method 
was we able to obtain nonnegative approximations with all the meshes we tried up to $N=1800$.  In sharp
contrast, with the SMS method, both Galerkin-based and SUPG-based, nonnegative approximations
were obtained for $N\ge 8$ and $\varepsilon$ as small as $10^{-14}$.
\end{example} 

\section{Concluding remarks}\label{Se:conclusions}

A novel stabilization technique for convection-diffusion problems has been introduced,
tested and partially analyzed. It can be applied to most existing methods based on
conforming piecewise linear finite elements. It consists of, first, adding extra values to
the residual of the method on nodes next to the outflow and characteristic boundaries.
The resulting equations are taken as a restriction on a least-squares problem on the
elementwise residuals of the convection-diffusion operator, where
elements next to the outflow and characteristic boundaries (and internal layers) are left out.

The method has been tried in a series of standard and nonstandard tests, and results
seem to suggests it performs manifestly better than a good deal of the methods of choice today.
The tests include exponential and characteristic layers, and even irregular grids on
domains with nontrivial geometries. This is so in spite of the method being initially
conceived as a simulation on the coarse part of a Shishkin mesh.
The tests also include interior layers and convection with vortices, with similarly outstanding results, although numerical
experiments in the present paper and in~\cite{sms2} suggest that
further research may be  needed on these  topics.

Besides the practical performance shown in tests, it is to be remarked the lack of parameters
in the method (in markedly contrast with most of existing stabilized methods today).
Subject of future research will be extending the method to finite elements of higher
degree, as well as the possible changes when the mesh P\'eclet number tends to one.

%
%
\medskip

\noindent{\bf Acknowledgements}. The author wish to thank Prof.\ Martin Stynes for advice
and helpful discussions in the research summarized in this paper.

\noindent{}
\vskip2.2truecm
\setcounter{page}{1}
\setcounter{equation}{0}
\setcounter{section}{0}
\setcounter{example}{0}
\setcounter{figure}{0}

\begin{center}
{\Large
 Shishkin mesh simulation: Complementary}

{\Large experiments}
\bigskip

{\large Bosco
Garc\'{\i}a-Archilla}\footnote{Departamento de Matem\'{a}tica Aplicada
II, Universidad de Sevilla, Sevilla, Spain. Research supported by
Spanish MEC under grant MTM2009-07849 (bosco@esi.us.es)}

\medskip
{\large
\date{\today}
}
\vskip 1truecm

\begin{abstract}
We present here some of the experiments that, for the sake of brevity were
not included~\cite{sms}. {\it The present paper is neither self-contained nor intended to be published
in any major journal\/}, but is intended to be accessible to anyone wishing to learn more on the
performance of the SMS method. It should be read as a complement to~\cite{sms}.
\end{abstract}

\end{center}
\section{Numerical experiments}

In its present version, this document is not selfcontained, and its contents must
be seen as add ons to~\cite{sms}. Recall that in that paper we are concerned with
\begin{align}
&-\varepsilon \Delta u+b\cdot \nabla u+ cu=f,\quad {\rm in}\quad \Omega,\label{eq:model_inc}\\
&u=g_1,\quad {\rm in}\ \partial \Omega_D,\quad \frac{{\partial
u}}{{\partial n}}=g_2,\quad {\rm in}\ \partial
\Omega_N,\label{eq:modelbc_inc}
\end{align}
where $\Omega$ is a bounded
domain in ${\mathbb R}^d$, $d=1,2,3$,
its boundary $\partial\Omega$ being the disjoint union of~$\Gamma_D$ and~$\Gamma_N$,
$b$ and $c$ are given functions and $\epsilon> 0$ is a
constant diffusion coefficient.

\begin{example}\label{ej:comparem_inc} {\it Simulation of  Shishkin grids}.
Let us also mention that we rerun the SUPG-based SMS manually tunning the streamline
diffusion parameter and found that, in this example, if set 1.5 times larger than the
 value shown in~\cite[Fig.~11]{sms}, the errors improve by a factor of two
approximately, resulting in the most efficient method and at least twice as
efficient as the SUPG on Shishkin grids
\end{example}
\medskip
\begin{example}\label{ej:comparem_same_inc} {\it Comparisons on the same grid}.
 The difference between the SUPG and the SMS methods is even
larger when errors are meassued in the $H^1$ norm, which are shown
in~Fig.~\ref{fig_comparemsameh1}
\begin{figure}[h]
$$
\includegraphics[height=5.3truecm]{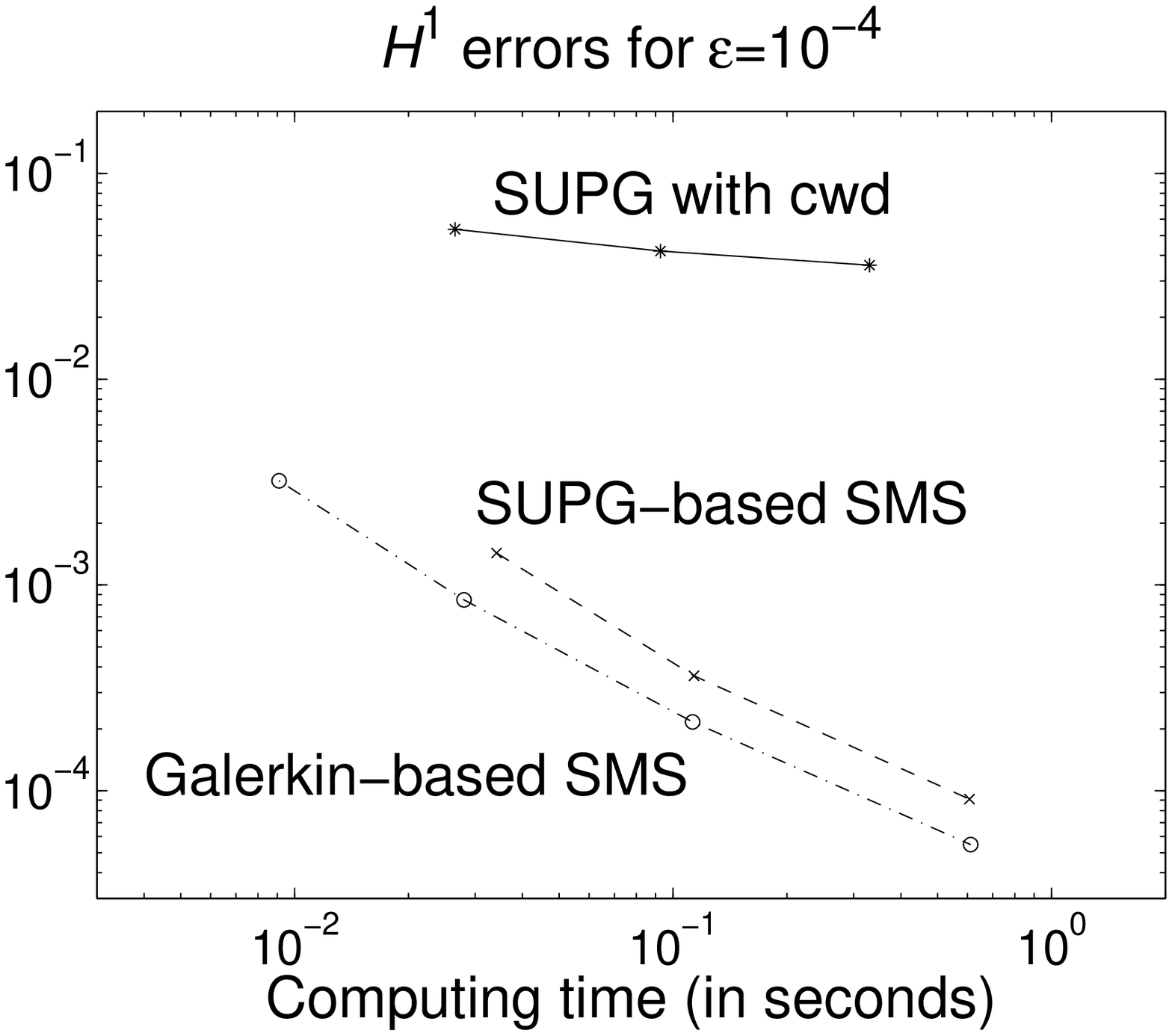}\quad
\includegraphics[height=5.3truecm]{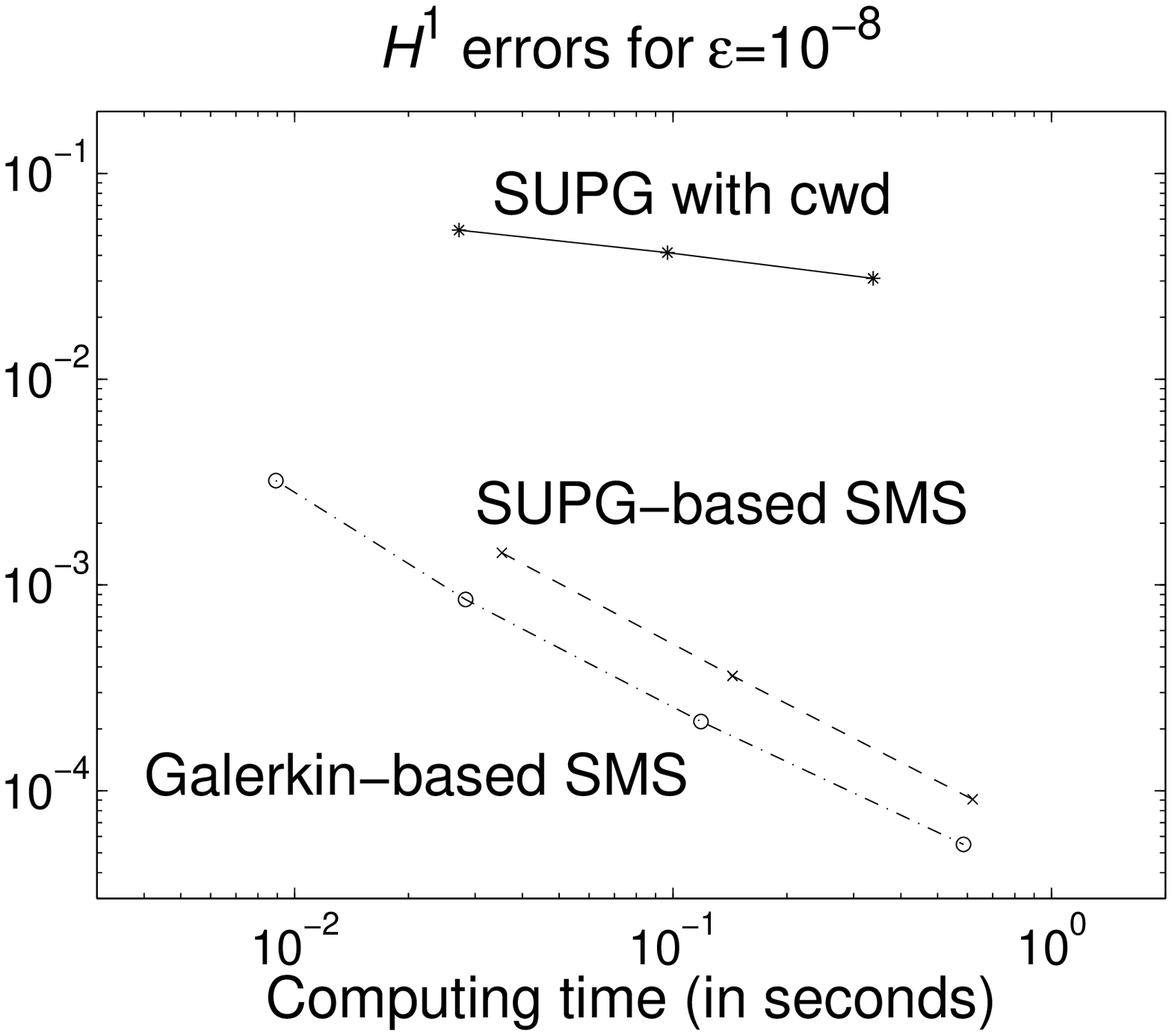}
$$
\caption{Relative efficiency of SMS and SUPG (with crosswind diffusion) methods in the $H^1$norm.}
\label{fig_comparemsameh1}
\end{figure}

\end{example}

\begin{example} \label{ej:curved_inc} {\it Irregular grids on curved domains}.

We study the behaviour of SMS methods on irregular grids. For a given positive integer $N$, we
consider grids with $(N+1)^2$ points in $\overline\Omega$, where $4N$ of them are
randomly distributed on $\partial\Omega$. Let us denote by $N_o$ the number of
these on~$\Gamma_D^{0+}$. Then,
in 
from the remaining~$(N-1)^2$ points, $N_o$ of them are displaced by the normal vector
from those~on $\Gamma_D^{0+}$ (this is done to create a strip of elements on the outflow boundary) and the rest are generated randomly according to the
following two schemes.
\begin{description}
\item{Highly irregular grids} The points are generated following uniform distribution on the square enclosing~$\Omega$, and those enclosed by~$\partial\Omega$ and the outflow strip are included
in the grid until the desired number $(N-1)^2-N_o$ of them is reached.

\item{Mildly irregular grids} A uniform grid is build
first inside $\Omega$ and the outflow strip, its diameter~$h$ in the $x$ and~$y$ direction being the value  for which the
number of points is the closest to~$(N-1)^2-N_o$. Then each point is displaced randomly
on the $x$ and~$y$ directions with a uniform distribution on~$[-h/3,h/3]$. This is the type of
grid used in~the experiments in~\cite{sms}.

\end{description}
An example with the two types of grids can be seen in~Fig.~\ref{fig_gridworst40}.

For different values of~$N$, we generated 200 random grids of both types, and on each of
them we computed the SUPG and SMS approximations and their $L^2$ error in the convective derivative
in~$\hat\Omega_h$, that is
$$
\left\| b\cdot \nabla w_h -1\right\|_{L^2(\hat\Omega_h)},
$$
for $w_h=u_h$ the SUPG approximation and $w_h=\tilde u_h$, the SMS approximation
(both Galerkin and SUPG-based). The errors on the 200 grids for~$N=80$ and~$\varepsilon=10^{-8}$
are depicted in Fig.~\ref{fig_random80}. They are marked with a dot, an asterisks and
a cross for the SUPG, Galerking-based SMS and SUPG-based SMS methods, respectively.
The results have been reordered so that
those of the SUPG method appear in descending order. 
\begin{figure}[h]
$$
\includegraphics[height=5.3truecm]{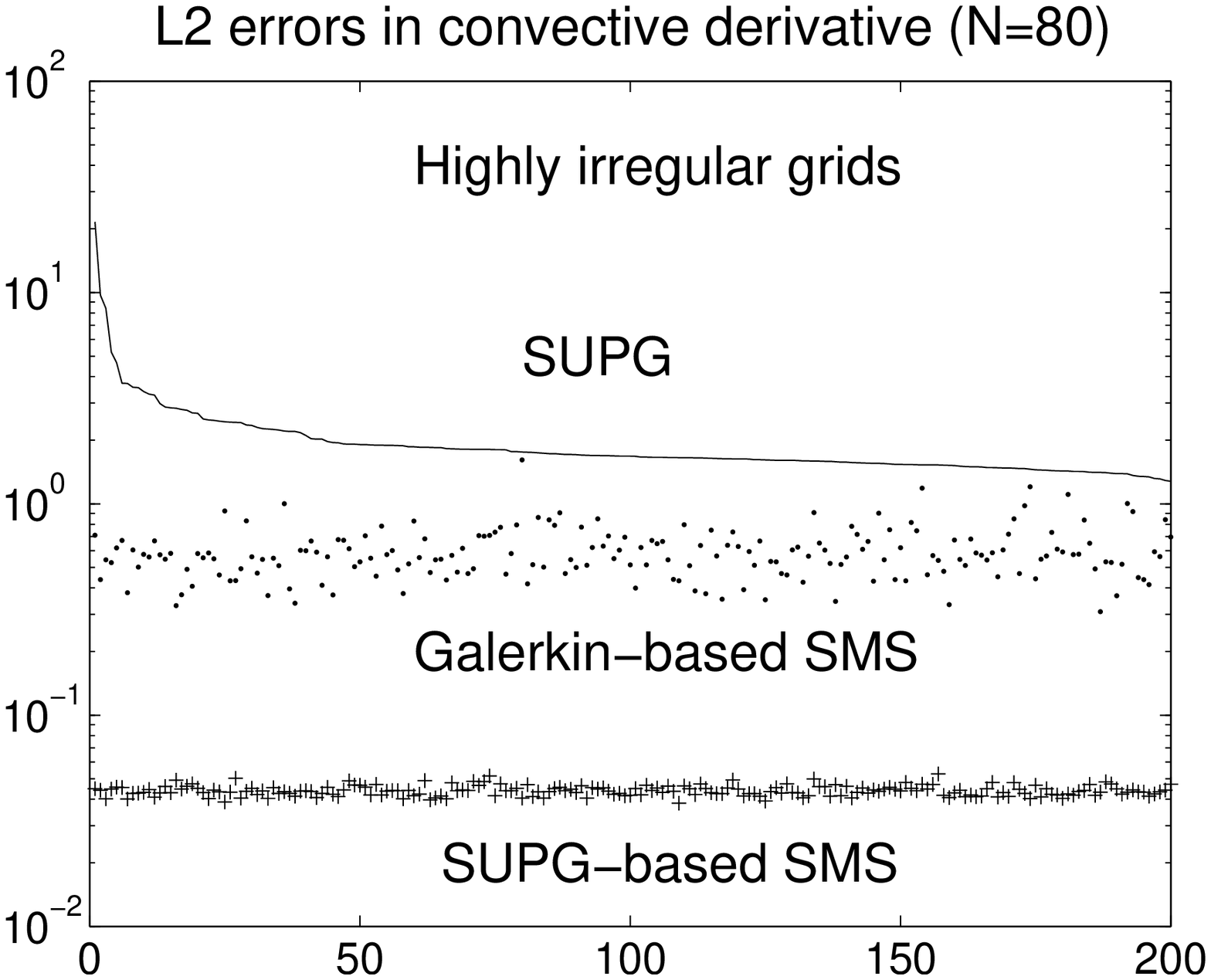}\quad
\includegraphics[height=5.3truecm]{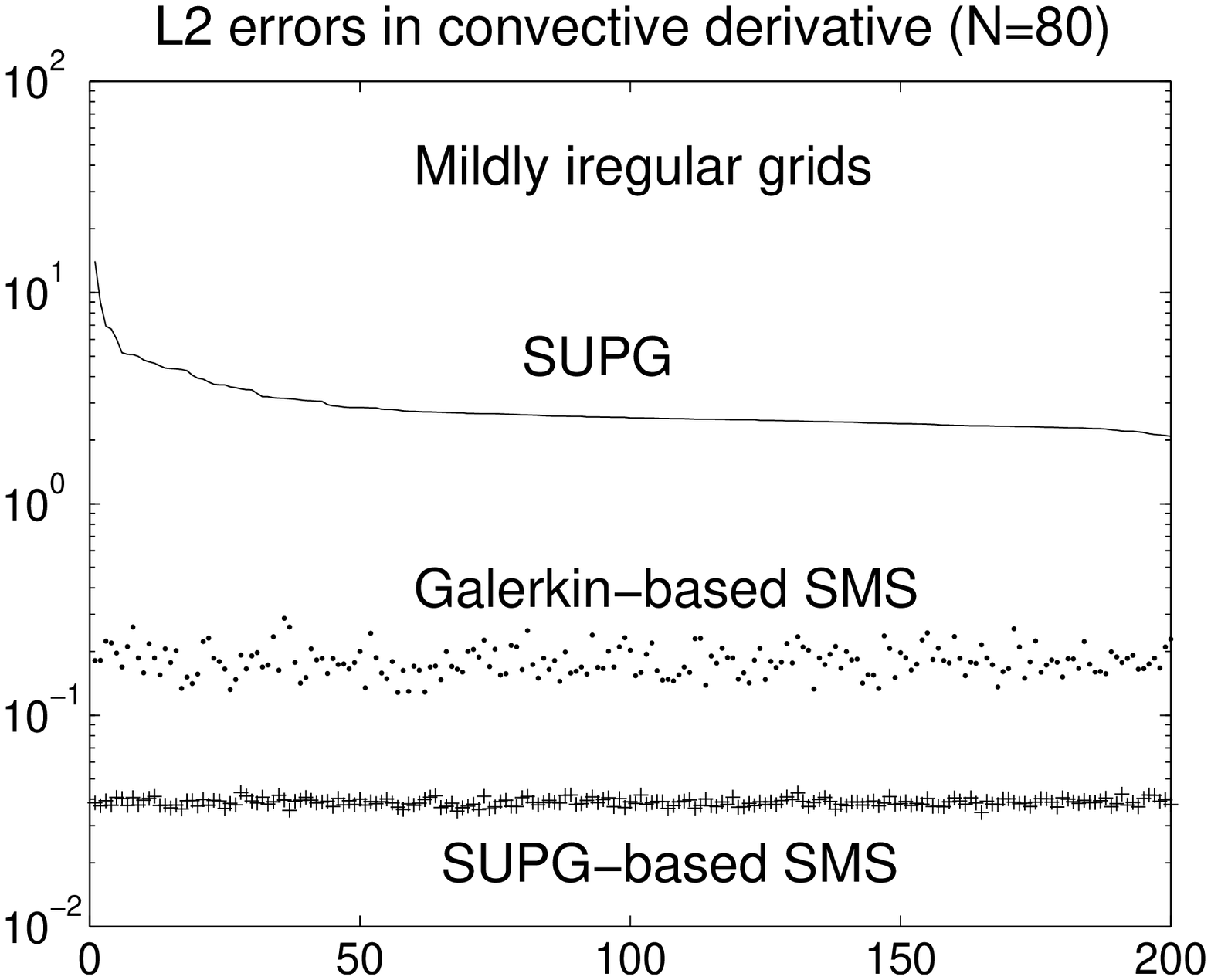}
$$
\caption{Errors in convective derivative on random grids.}
\label{fig_random80}
\end{figure}
As in the results shown in \cite{sms}, the SMS methods clearly improve the error with respect the SUPG method, particularly so
in the case of the SUPG-based SMS method. Computing the ratios of the error of the SUPG method
and each of the SMS methods and taking the arithmetic mean, the resulting values are for
the Galerkin-based and SUPG-based SMS methods, respectively, 5.17 and 65.57 on highly irregular
grids, and 16.30 and~75.16 on mildly irregular grids. That is, the SMS methods commit errors
that are, on average, between 5 and 75 times smaller than those of the SUPG method.
 
In Fig.~\ref{fig_random80} we can also see that while the SUPG and SUPG-based SMS methods are insensitive to the irregularity of the grids, this is not the case of the Galerkin-based SMS method.
 This is
also reflected on Fig.~\ref{fig_worst40}, where we show the worst cases 
of Galerkin-based SMS method on both types
of random grids for~$N=40$: whereas in the highly irregular grids the approximation exhibits
oscillations, this not so noticeable on the mildly irregular grid (the corresponding grids are shown
in~Fig.~\ref{fig_gridworst40}).
\begin{figure}[h]
$$
\includegraphics[height=4.5truecm]{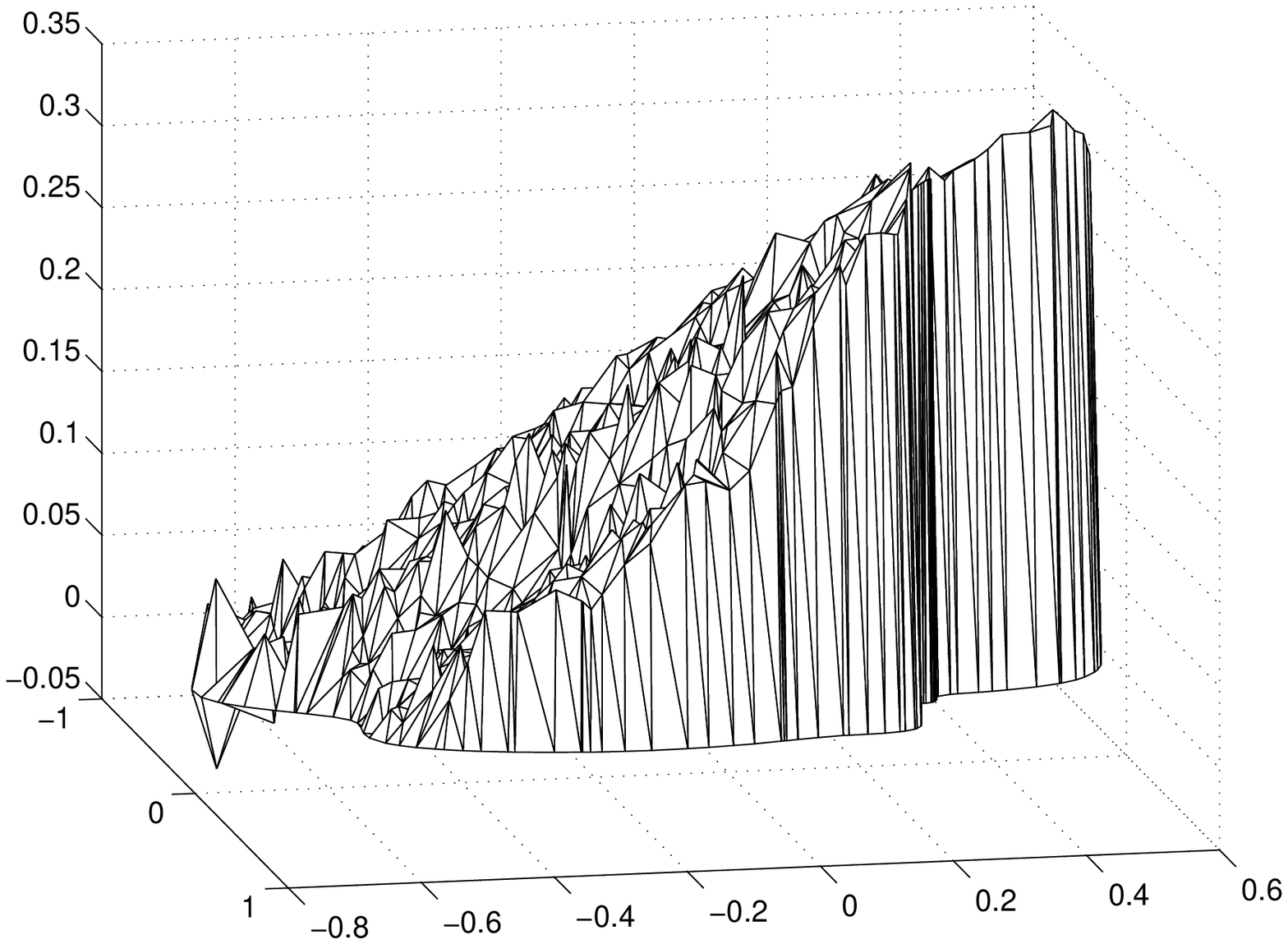}\quad
\includegraphics[height=4.5truecm]{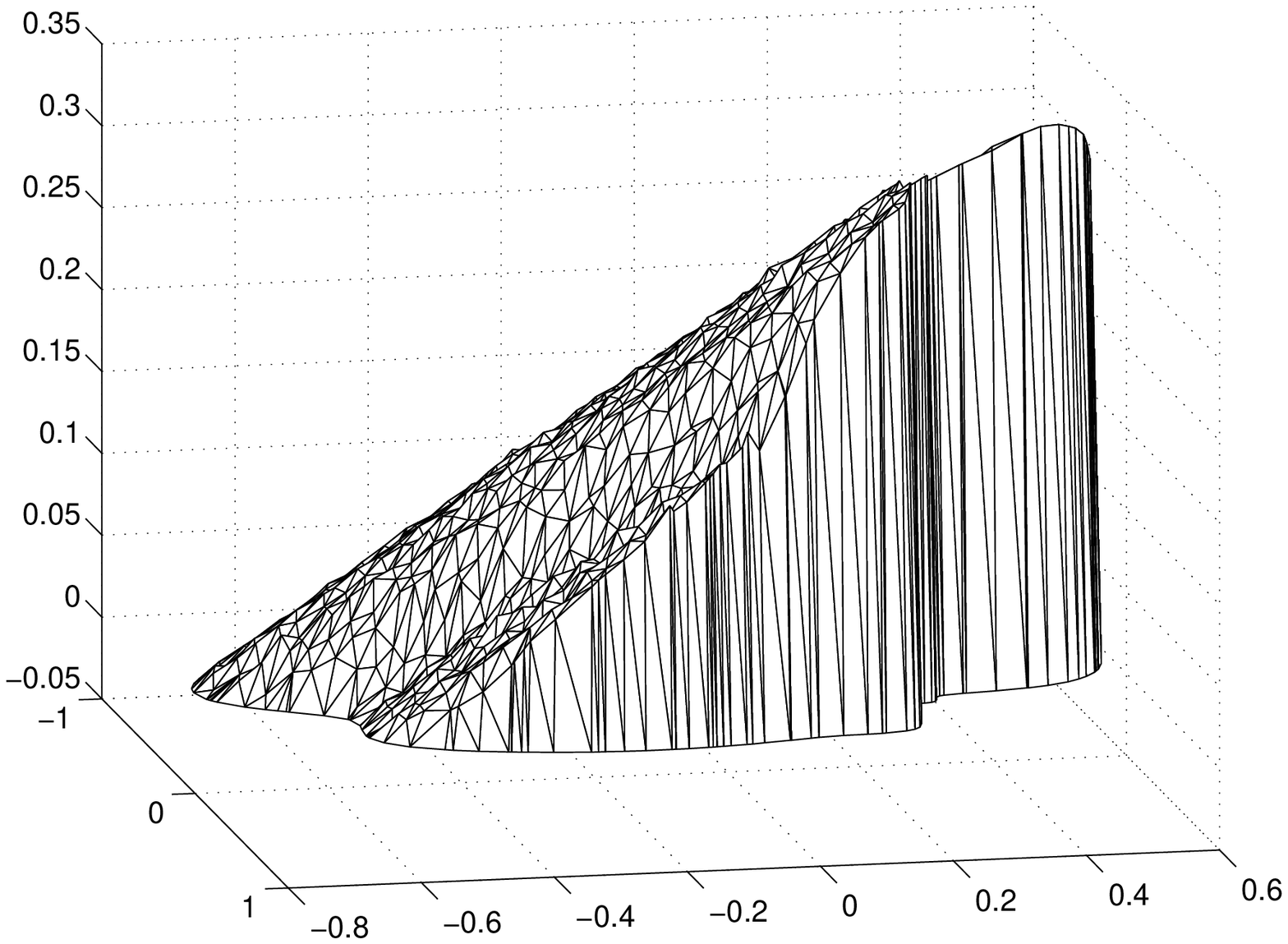}
$$
\caption{The worst cases of the Galerkin-based SMS method for $N=40$ on highly irregular
 and mildly irregular grids (left and right, respectively).}
\label{fig_worst40}
\end{figure}
As commented in~\cite{sms}, some results explaining the
degradation of performance of the Galerkin method on irregular grids can be found in~\cite{Chen-Xu_inc} and~\cite{Sun-Chen-Xu_inc}. However, we must remark that the approximation on the left plot
in~Fig.~\ref{fig_worst40} is the worst case out of 200, and that in most of the cases on
highly-irregular grids the Galerkin-based SMS method did not present spurious oscillations. Nevertheless, due to the better performance of the Galerkin-based SMS method
on mildly-irregular grids, in the rest of the paper, we will only deal with either regular or
mildly-irregular grids. 


\begin{figure}[h]
$$
\includegraphics[height=5.5truecm]{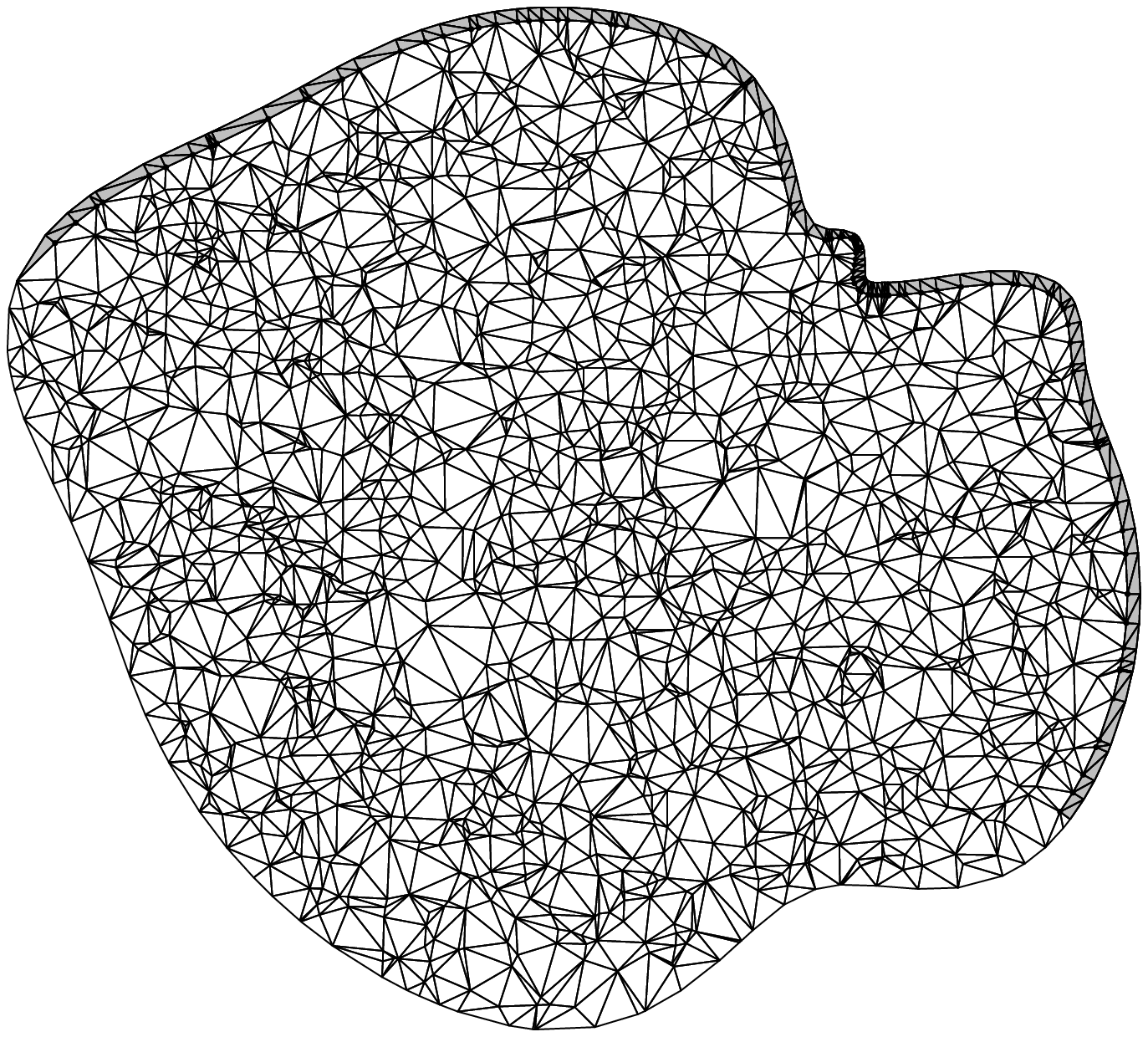}
\includegraphics[height=5.5truecm]{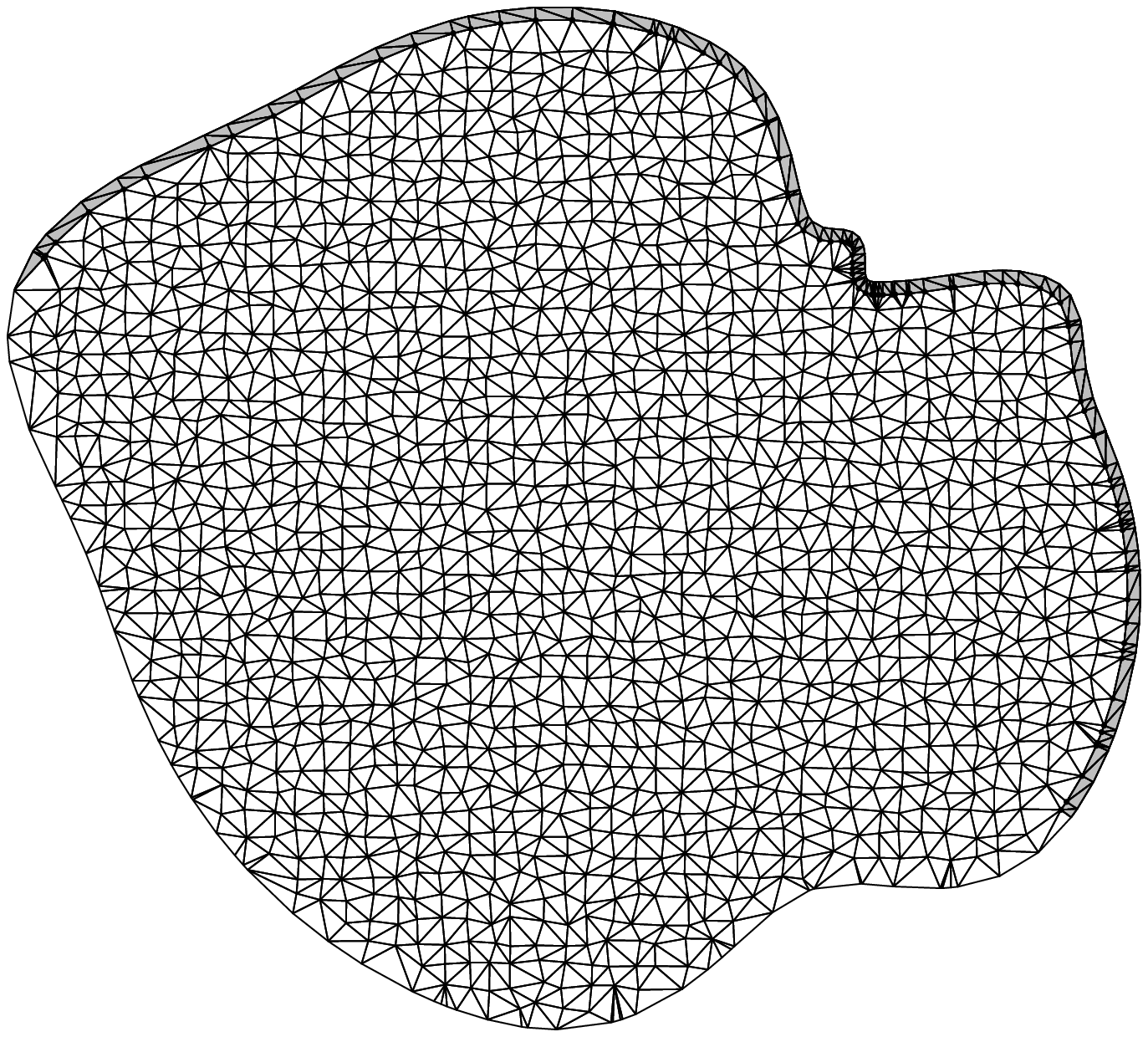}
$$
\caption{The grids of the plots in Fig.~\ref{fig_worst40}. The set~$\Omega_h^{+}$ is shadowed
in grey and the points in~${\cal N}_\delta$ are marked with a dot.}
\label{fig_gridworst40}
\end{figure}
\medskip

We show the arithmetic mean of the errors committed by
the different methods on mildly irregular grids for each value of~$N$ in~Fig.~\ref{fig_dibconv}.
\begin{figure}[h]
$$
\includegraphics[height=5.5truecm]{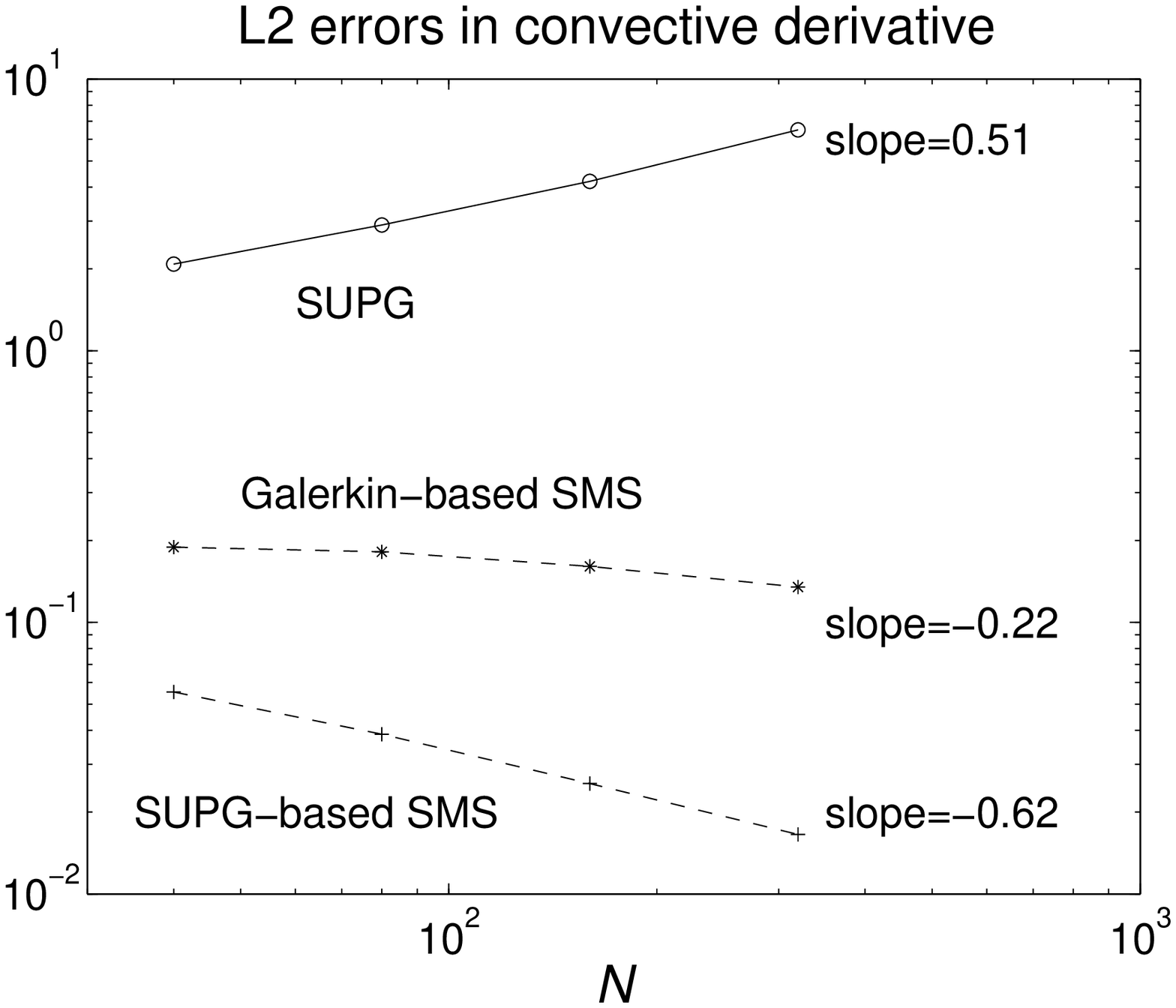}
$$
\caption{Arithemetic means of the errors on mildly irregular grids.}
\label{fig_dibconv}
\end{figure}
We also show the slopes of least-squares fits to the different sets of results. We observe
that the errors in the SUPG-based SMS method decay as the meshes become finer
(at a rate of $N^{-0.6}$) marginally so in the case of the Galerkin-based SMS method and
they grow in the SUPG method. That is, the SUPG method has not reached
the convergence regime yet for these meshes in this problem.  We can have an idea
of why the errors grow
by looking at~Fig.~\ref{fig_worst40supg}, where we show the SUPG approximation on
the grids shown in~Fig.~\ref{fig_gridworst40} (compare with the Galerkin-based SMS
in~Fig.~\ref{fig_worst40}). We notice that large amplitude of
the oscillations near the outflow boundary. The large amplitude of the oscillations
decreases slower than
the grid size, so that the errors in the convective derivative grow.
\begin{figure}[h]
$$
\includegraphics[height=4.5truecm]{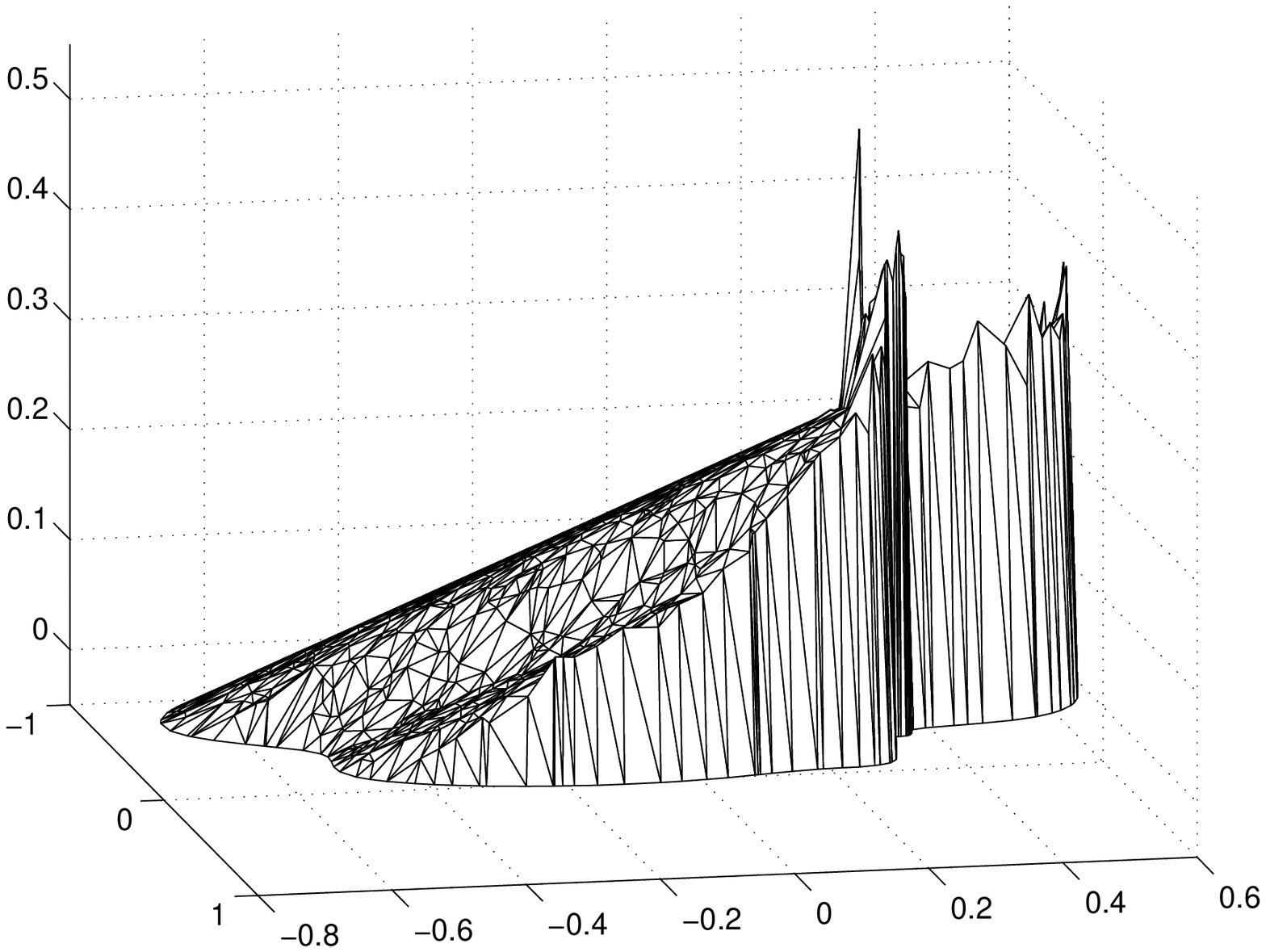}\quad
\includegraphics[height=4.5truecm]{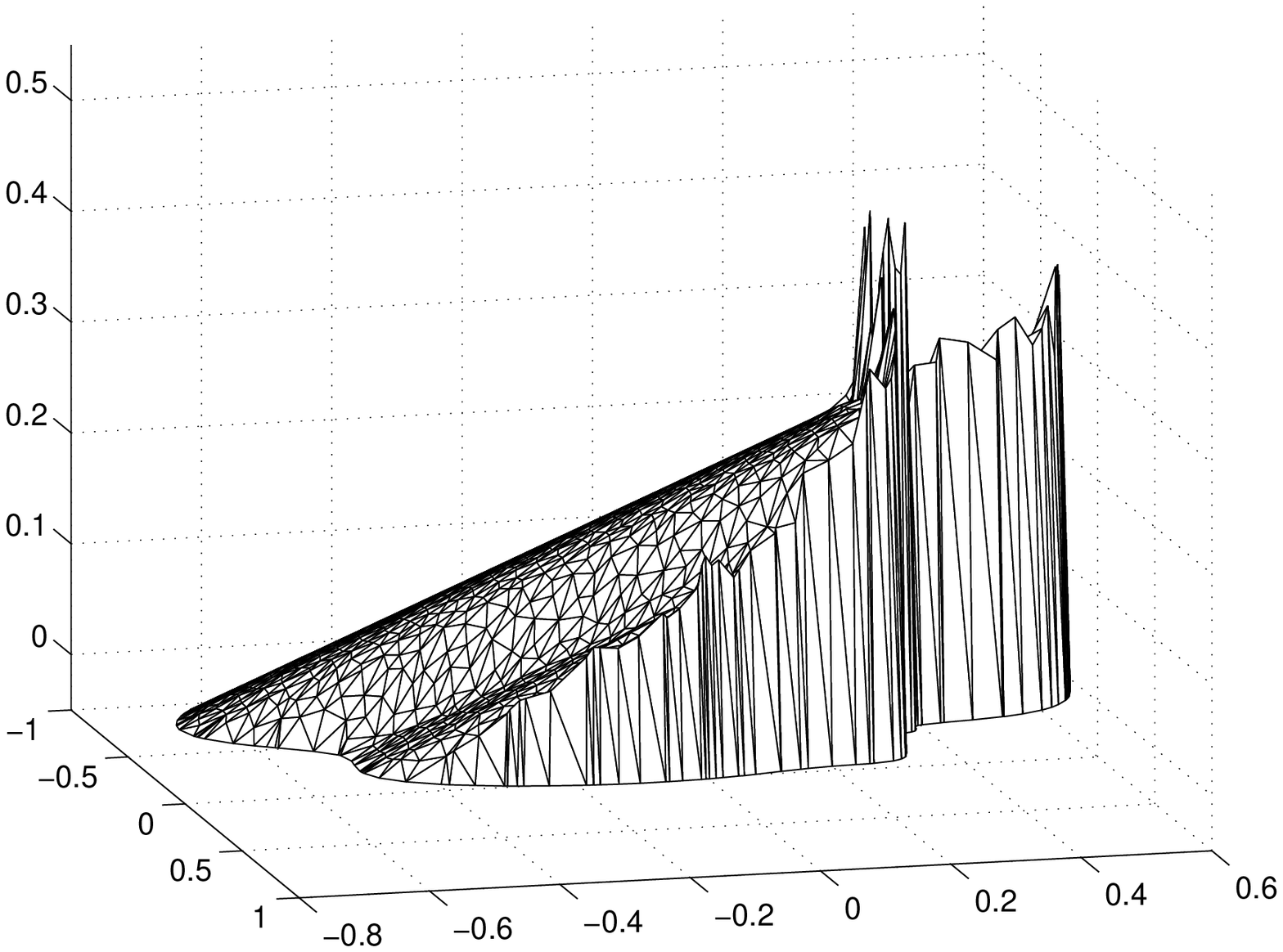}
$$
\caption{The SUPG approximation on the grids shown in~Fig.~\ref{fig_gridworst40}.}
\label{fig_worst40supg}
\end{figure}
\end{example}
\medskip

\begin{example} {\it Parabolic Layers}. \label{ej:parabolic_inc}
Recall that in~\cite{sms} we solve~(\ref{eq:model_inc})
on $\Omega=(0,1)^2$, with $\varepsilon=10^{-8}$,
$b=[1,0]^T$ and~$f$ constant equal to~$1$ with Dirichlet homogeneous boundary conditions.
The solution has one exponential layer at $x=1$ and two parabolic layers along~$y=0$ and~$y=1$.

In~\cite{John-Knobloch-2007_inc}, a random grid is used and the sets
$\Omega_2=(0,0.9)\times(0,0.1]$ and~$\Omega_3=(0,0.9)\times[0.9,1]$, and
$\Omega_4=[0.9,1)\times(0.1,0.9)$ are considered, and the following quantities computed
for numerical approximations~$w_h$:
\begin{align}
\hbox{osc}_{\text{para(2)}}&:=
\max\bigl\{\max_{(x_s,y_s)\in\Omega_2}-\partial_y w_h(x_s,y_s),\max_{(x_s,y_s)\in\Omega_3}\partial_y w_h(x_s,y_s)
\bigr\},
\label{osc_para2}\\
\hbox{osc}_{\text{exp}}&:= \max_{(x_s,y_s)\in\Omega_4} \partial_x w_h(x_s,y_s)
\label{osc_exp}
\end{align}
$(x_s,y_s)$ being the barycenters of the triangles.
The optimal values for~$\hbox{osc}_{\text{para(2)}}$ and~$\hbox{osc}_{\text{exp}}$ are, respectively,
0 and~$1$, the larger these values are, the stronger are the oscillations in the characteristic and
exponential layers, respectively.
None of the methods tested in~\cite{John-Knobloch-2007_inc} presents a value
of~$\hbox{osc}_{\text{para(2)}}$ below~$10^{-2}$, nor a value of~$\hbox{osc}_{\text{exp}}$
that is less than $10^{-2}$ away from the target value~$1$. We try 200 random
irregular grids with $1681$ interior grid points (the grid in~\cite{John-Knobloch-2007_inc} had
1721), which were created as in~Example~\ref{ej:curved_inc}.  On these grids,
the highest values of $\hbox{osc}_{\text{para(2)}}$ and
$\left| \hbox{osc}_{\text{exp}}-1\right|$ in the SMS methods were, respectively, $2.9\times 10^{-13}$
and~$2.8\times 10^{-13}$. Similar striking differences were observed in the rest of the
quantities measured in~\cite{John-Knobloch-2007_inc} in this example.

Also, with highly irregular
grids, as long as they had a strip of elements along~$\Gamma_D^{0+}$, the results were similar to
those of mildly-irregular grids.
\end{example}
\medskip

\medskip
\begin{example} {\it Interior layers}. \label{ej:interior_inc} Nothing to be added here at this stage.
This example is revisited after Exameple~\ref{ej:Hemker_inc} below.
\end{example}

\medskip
\begin{example}\label{ej:Hemker_inc} {\it Hemker problem}.
We test the technique to include the layer characteristic in the grid that was proposed
in~\cite[Example~6]{sms}.
To do this, we change
$b$  and the boundary conditions to
\begin{equation}
\label{hemkertheta}
b=\left[\begin{array}{c}\cos(\theta)\\ \sin(\theta)\end{array}\right],
\qquad
u(x,y)=\left\{\begin{array}{lll} 0,&\ &\hbox{\rm if $x=-3$ or $y=-3$,}\\ 1,& &\hbox{\rm if $x^2+y^2=1$,}\\
\varepsilon\nabla u\cdot n=0,&& \hbox{\rm elsewhere}\end{array}\right.
\end{equation}
for $\quad\theta\in(0,\pi/4]$.
(The boundary conditions are changed so that the inflow boundary $\Gamma^{-}$ is part of the
Dirichlet boundary~$\Gamma_D$). The idea is to try different combinations of characteristic layers
and grid alignments. In the following computations, the enlarged grids and the
corresponding sets~$\Omega_h^{+}$ can be as irregular as those depicted
in~Fig.~\ref{fig_spikes}.
\begin{figure}[h]
$$
\includegraphics[height=4truecm]{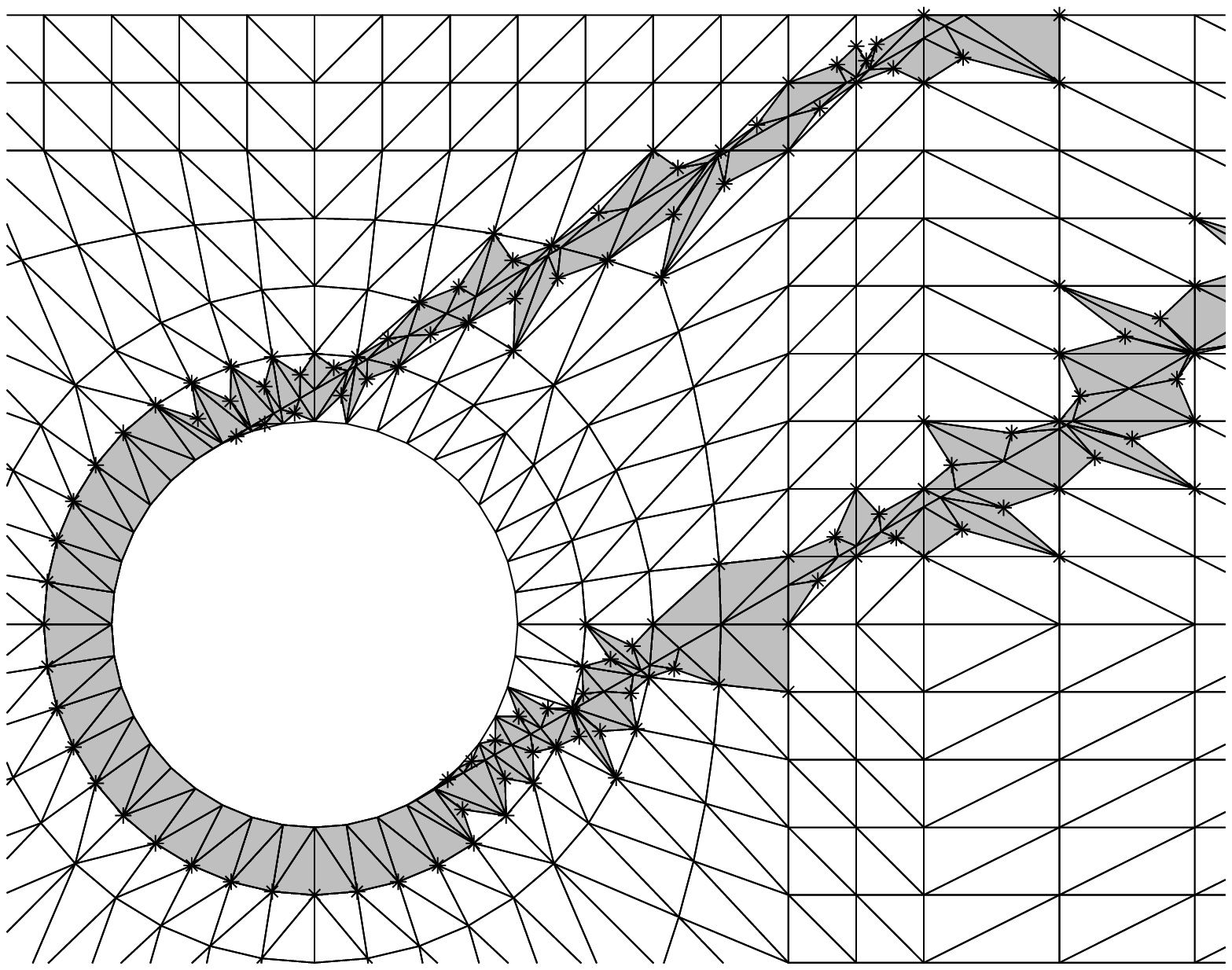}
$$
\caption{Detail of the enlarged grid, the set $\Omega_h^{+}$ (shadowed in grey) and the
points in~${\cal N}_\delta$ (marked with asterisks) in Example~\ref{ej:Hemker_inc} for $\theta=\pi/6$ in~(\ref{hemkertheta}).}
\label{fig_spikes}
\end{figure}
To prevent elements with one side parallel to~$b$ being downwind of~$\Omega_h^{+}$
(see~\cite[Section~3.2.2]{sms} for the possible adverse effects of this)
when $\theta=\pi/4$ we red-refined all triangles in~$\Omega_h^{+}$ along the interior layers.
Otherwise, no extra provisions where taken in the computations that follow, except removing isolated components in  $\mathring{\hat\Omega{}}_h$ (when present) as explained
in~\cite[Section~3.2.2]{sms}.

We considered 100 equidistant values of~$\theta$ in~$(0,\pi/4]$. For each of them we computed
the overshoots and undershoots of the SUPG and SMS approximations
on the grid depicted in~\cite[Fig.~20]{sms}.
For the SUPG and the Galerkin-based SMS approximations,
Fig.~\ref{fig_hemkerovers} shows the difference between the overshoots and undershoots,
as a measure of the oscillations in each method
(the results of the SUPG-based SMS method were similar to those of the Galerkin-based method
and are not shown).
\begin{figure}[h]
$$
\includegraphics[height=3truecm]{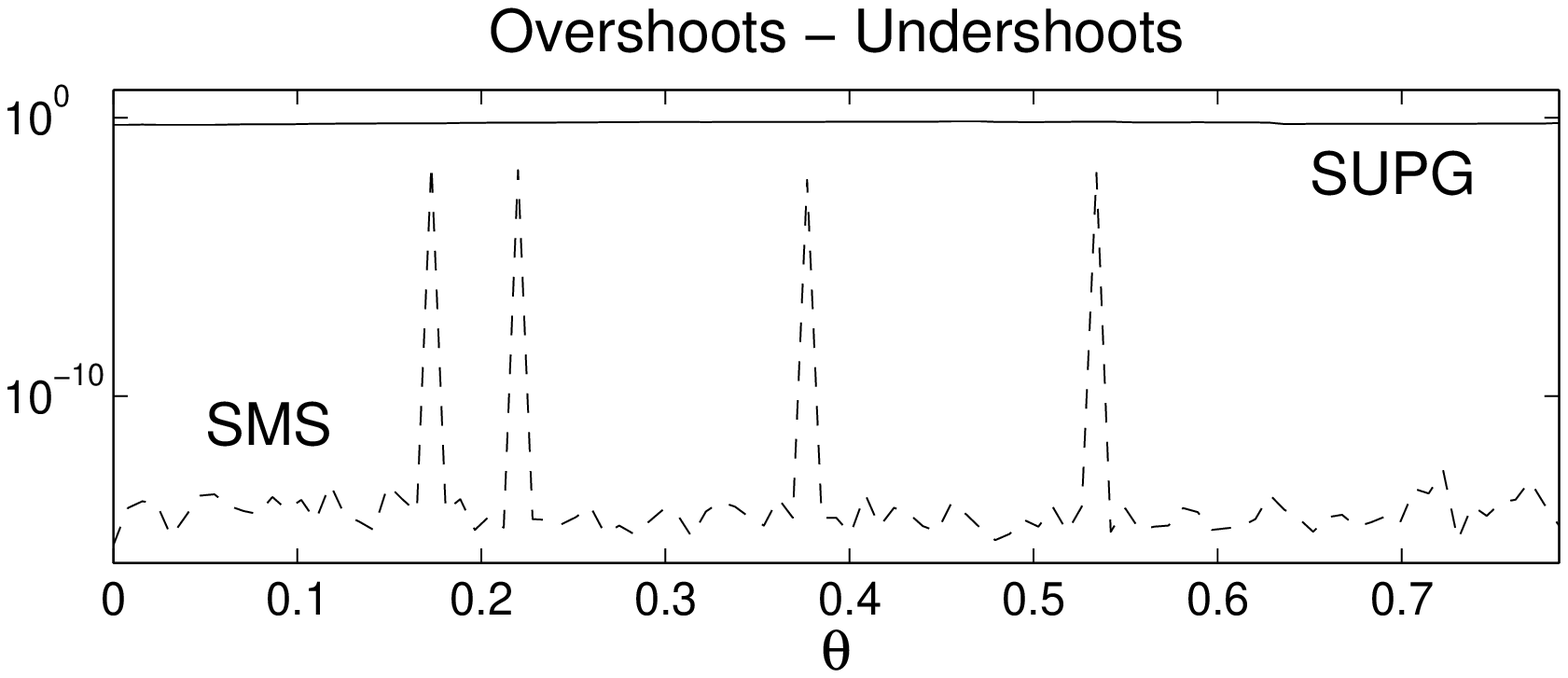}
$$
\caption{Difference between the over and undershoots for
the approximations on the grid depicted in~\cite[Fig.~20]{sms} in Example~\ref{ej:Hemker_inc},
with $b$ and boundary conditions as specified in~(\ref{hemkertheta}).}
\label{fig_hemkerovers}
\end{figure}
We see that out of the 100 values of~$\theta$ tried, only in four of them are the oscillation on
the SMS method
larger than~$10^{-12}$, the arithmetic mean of these four values being $0.013$, yet well below the
average value in the SUPG method, $0.65$.

We investigate the reason of those four cases presenting values of over and undershoots so different
to the rest of the cases. The reason is that some mesh points of the enlarged grid are too close to the interior layer, and this
results in the set~$\Omega_h^{+}$ being too thin in the neighbourhood of those points.
This can be seen in Fig.~\ref{fig_thin}, where on the left plot we show a detail of the
enlarged grid (with $\Omega_h^{+}$ shadowed in grey) for $\theta=23\pi/400$, which is the
first value in~Fig.~\ref{fig_hemkerovers} where the overshoots or undershoots are above~$10^{-12}$
in the SMS method. The center plot is a magnification of the left one around the point too close
to the layer characteristic, where we can see that very stretched triangles have been created
on the enlarged grid.
\begin{figure}[h]
$$
\includegraphics[height=2.5truecm]{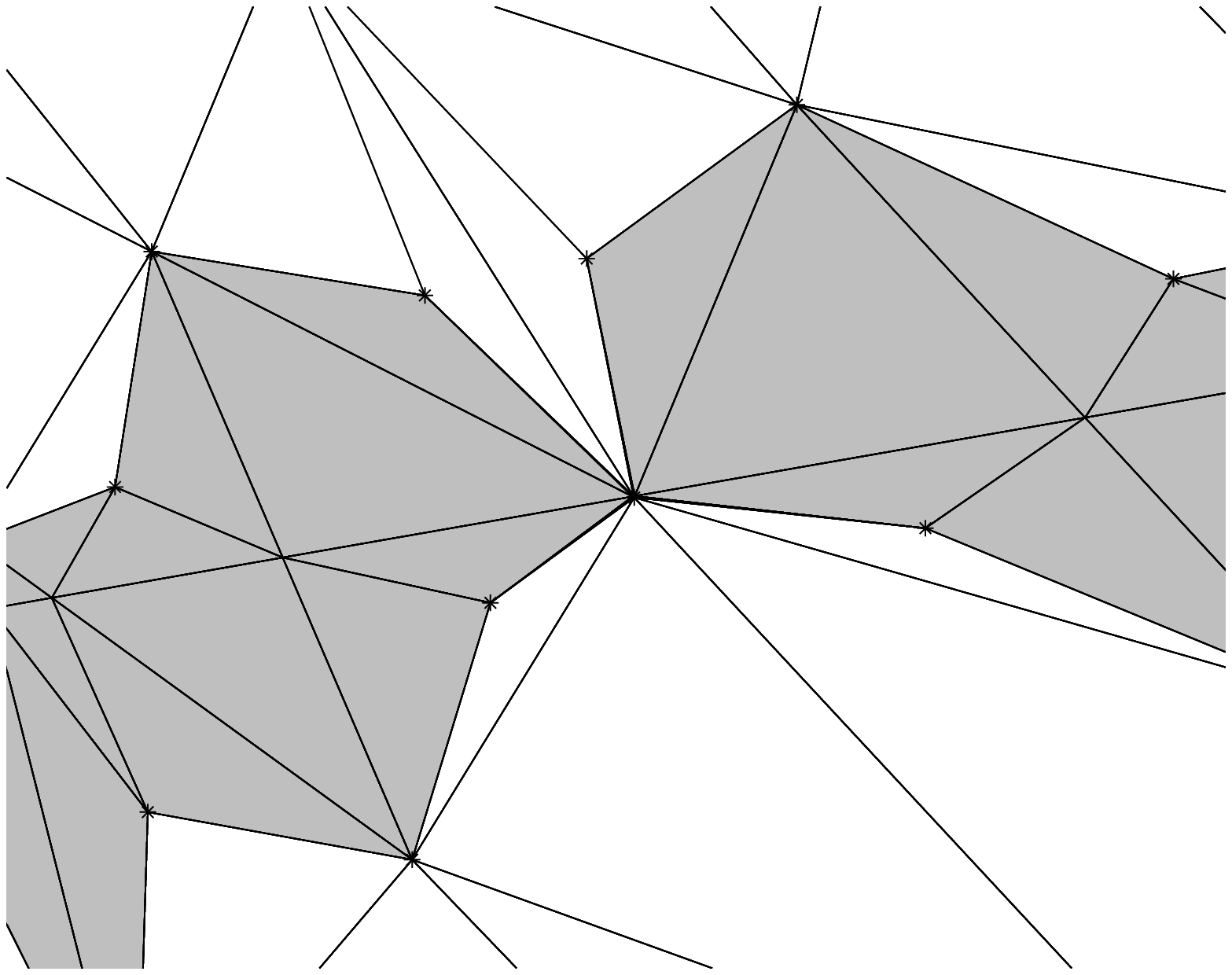}\quad
\includegraphics[height=2.5truecm]{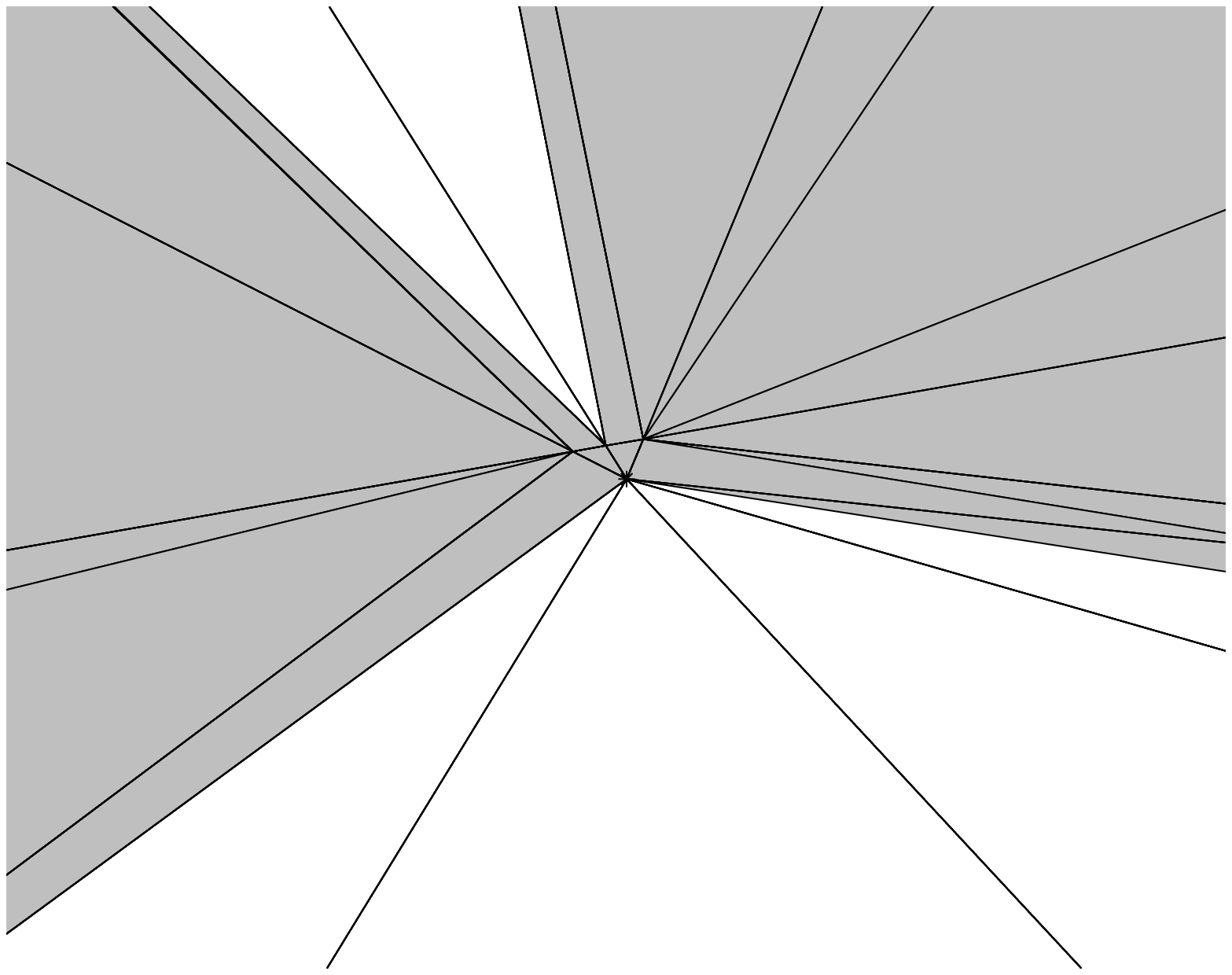}\quad
\includegraphics[height=2.5truecm]{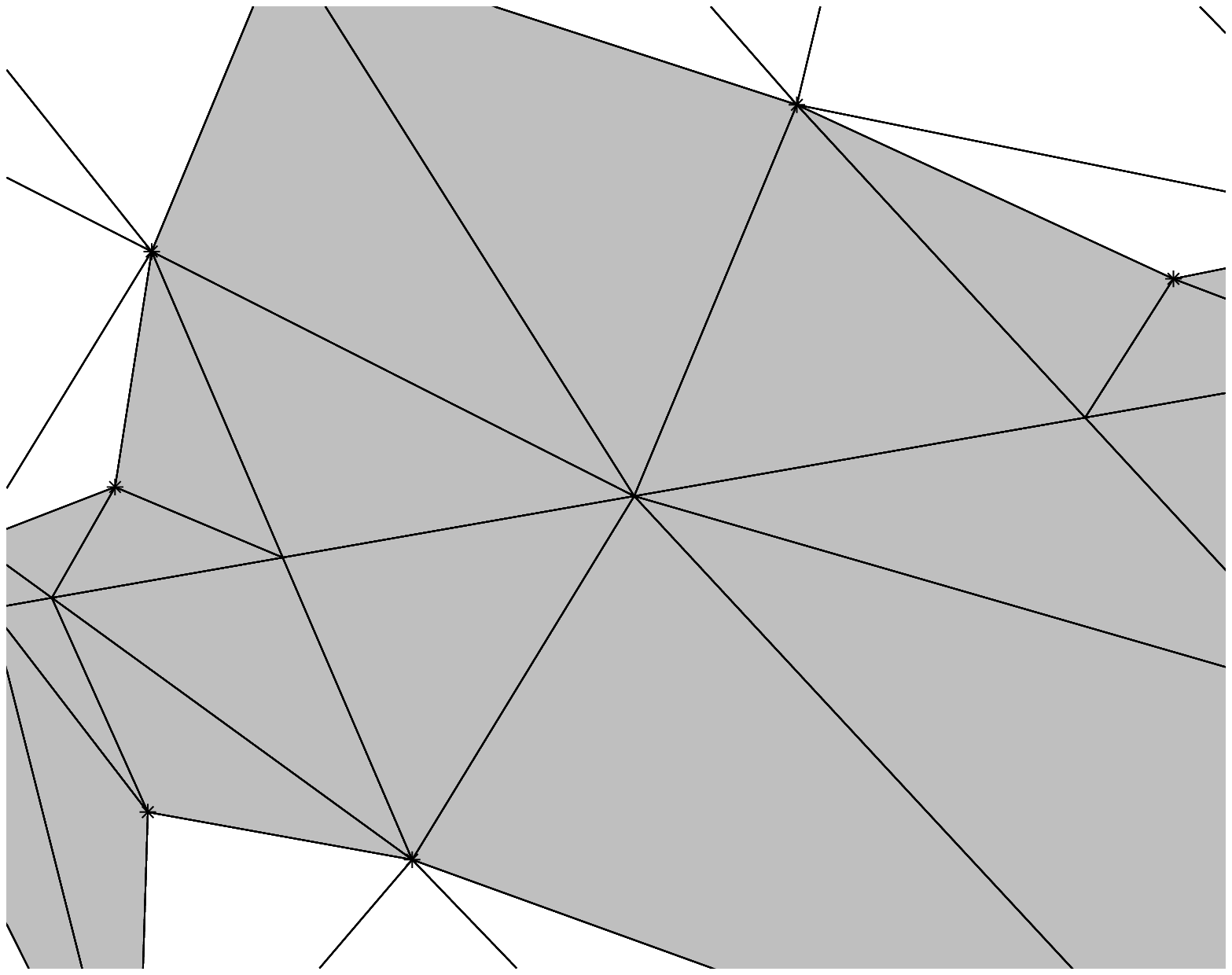}
$$
\caption{Left: detail of the enlarged grid and $\Omega_h^{+}$ for $\theta=23\pi/400$
in~Example~\ref{ej:Hemker_inc} with~(\ref{hemkertheta}). Center: magnfication of the left
plot around a point too close to the characteristic layer. Right:
the set $\Omega_h$ after moving that point to the layer characteristic.}
\label{fig_thin}
\end{figure}
If we move the point to the layer characteristic, the resulting set $\Omega_h^{+}$ is shown
in~the right plot of~Fig.~\ref{fig_thin}, and we can see it is much wider than that in the
left plot.

We repeated the four cases in the SMS method were oscillations were much larger
than the rest of the cases,
but this time allowing 
a vertex $v$ of a triangle~$\tau$ to move
to the layer characteristic if its distance to it was less than~$(h_{\min}(\tau))^2/10$, were
$h_{\min}(\tau)$ is the minimum distance between the vertices of~$\tau$. 
Once the SMS approximations
were computed on the moved grids, values on the unmoved grid were computed
by interpolation and over and undershoots were again computed. The overshoots in all cased
turned out to be below $10^{-14}$, and the undershoots were zero in all cases.
Notice that when not altering the grid (not even by less that $h_{\min}^2/10$) is an important issue,
then, another option to avoid mesh points too close to the layer characteristic is to slightly perturb
the wind velocity~$b$ in the neighbourhood of those points so that the layer characteristic passes through them.
\end{example}
\medskip

\noindent{\bf Example 5 (revisited).} As in example~\ref{ej:parabolic_inc}, we compute with our method
some of the quantities used in~\cite{John-Knobloch-2007_inc} to measure the quality of the
approximations. 
For an approximation $w_h\in V_h$, these are
\begin{align}
\label{osc_int}
\hbox{\rm osc}_{\text{int}}&:=\biggl(\sum_{(x,y)\in\Omega_1}
\bigl(\min\{0,w_h(x,y)\}\bigr)^2+\bigl(\max\{0,w_h(x,y)-1\}\bigr)^2\biggr)^{1/2},
\\
\hbox{\rm smear}_\text{int}&:=x_2-x_1,
\label{smear_int_inc}
\end{align}
where $\Omega_1=\{(x,y)\in \Omega \mid x\le 0.5,\quad y\ge 0.1\}$, and $x_1$ and~$x_2$ are
the first points on the line $y=0.25$ satisfying, respectively, $w_h(x_1,0.25)\ge 0.1$
and~$w_h(x_2,0.25)\ge 0.9$. It is argued in~\cite{John-Knobloch-2007_inc} that~(\ref{smear_int_inc})
is a measure of the thickness of the interior layer.  Observe also that the target value for
$\hbox{\rm osc}_{\text{int}}$ is 0.

We computed the value of these quantities in the case of the SMS approximations on a
regular $64\times 64$ grid with Southwest-Northeast diagonals, using the technique of
the previous example to include the layer characteristic in the grid.
On a similar mesh, out of 18 methods methods tested in~\cite{John-Knobloch-2007_inc} 
the value of~$\hbox{\rm osc}_{\text{int}}$ was larger than~$0.1$ in~14 of them, and larger than
$0.001$ in three more of them. In contrast,
the value of~$\hbox{\rm osc}_{\text{int}}$ in the SMS methods was
$1.8\times 10^{-14}$ (Galerkin-based) and~$2.2\times10^{-13}$ (SUPG-based). Also, the
vale of~$\hbox{\rm smear}_\text{int}$ was larger than  $6.2\times 10^{-2}$ in all the methods
tested in~\cite{John-Knobloch-2007_inc}, whereas for both SMS methods was~$1.3\times 10^{-2}$.
For the SUPG method, the values of~$\hbox{\rm osc}_{\text{int}}$
and~$\hbox{\rm smear}_\text{int}$ were, respectively, $0.59$ and~$0.062$, the same values
as in~\cite{John-Knobloch-2007_inc}. 

On the irregular grid in~\cite{John-Knobloch-2007_inc}, they obtain results were very similar to those
of the regular grid. We, on our part, tried the 200 random grids of Example~4 in the present one.
The mean value~$\hbox{\rm smear}_\text{int}$ in the SMS methods was $2.7\times10^{-2}$.
With respect the value
of~$\hbox{\rm osc}_{\text{int}}$ in the SMS methods, out of the 200 runs, in only 20 of them
was this value larger than $10^{-8}$, being the mean value in these 20 cases of 0.1822, which still compares very favourable with most of the methods
tested in~\cite{John-Knobloch-2007_inc}.  We rerun these 20 cases allowing a grid points to move
to the layer characteristic if the distance to it was less $2h_{\min}^2$, and in 19 of them the value
of~$\hbox{\rm osc}_{\text{int}}$ move to under~$10^{-10}$. In the remaining case this was also
achieve if we red-refining the triangles in~$\Omega_h^{+}$ next to the inner layer characteristic.

\medskip
\begin{example}\label{ej:Silvester_inc} {\it Double-glazing test problem\/} \cite{Silvester-book_inc}.
Whatever we may have to add in this case will be reported elsewhere.
\end{example}


\end{document}